\documentclass{amsart}

\usepackage{amssymb}
\usepackage{graphicx}
\usepackage{epstopdf}
\usepackage{amscd}
\usepackage{graphics}
\usepackage{young}
\newtheorem{Theorem}{Theorem}[section]
\newtheorem{Proposition}[Theorem]{Proposition}
\newtheorem{Corollary}[Theorem]{Corollary}
\newtheorem{Lemma}[Theorem]{Lemma}

\def\proof{\par{\it Proof}. \ignorespaces}
\def\endproof{{\ \vbox{\hrule\hbox{%
     \vrule height1.3ex\hskip0.8ex\vrule}\hrule }}\par}
\newenvironment{Proof}{\proof}{\endproof}

\theoremstyle{definition}
\newtheorem{Definition}{Definition}[section]
\newtheorem{Example}{Example}[section]
\newtheorem{Notation}{Notation}[section]

\theoremstyle{remark}
\newtheorem{Remark}{Remark}[section]

\numberwithin{equation}{section}

\let\trueint=\int
\let\truesum=\sum
\def\int{\mathop{\textstyle\trueint}\limits}
\def\sum{\mathop{\textstyle\truesum}\limits}

\let\<=\langle
\let\>=\rangle
\def\Real{{\mathbb{R}}}
\def\Complex{{\mathbb{C}}}
\def\F{{\mathbb{F}}}

\def\G{{\mathcal G}}
\def\S{{\mathcal S}}

\def\half{{\textstyle\frac12}}

\def\Gr{{\rm Gr}}
\def\SL{{\rm SL}}
\def\SO{{\rm SO}}

\def\mapright#1{\smash{\mathop{\longrightarrow}\limits^{#1}}}
\def\mapdown#1{\Big\downarrow\rlap{$\vcenter{\hbox{$\scriptstyle#1$}}$}}

\begin{document}
 
\title[Cohomology of real Grassmannian and KP flow]
{Cohomology of real Grassmann manifold and
KP flow}

\author{Luis Casian}
\address{Department of Mathematics, Ohio State University, Columbus,
OH 43210}
\email{casian@math.ohio-state.edu}
\author{Yuji Kodama$^*$}
\thanks{$^*$Partially
supported by NSF grant DMS0806219}

\address{Department of Mathematics, Ohio State University,
Columbus, OH 43210}
\email{kodama@math.ohio-state.edu}


\begin{abstract}
We consider a realization of the real Grassmann manifold $\Gr(k,n)$ based on a particular flow
defined by the corresponding (singular) solution of the KP equation. Then 
we show that the KP flow can provide an explicit and simple construction of the incidence graph for the integral cohomology of $\Gr(k,n)$.  
It turns out that there are two types of graphs, one for the trivial coefficients and other for the twisted
coefficients, and they correspond to the homology groups of the orientable and non-orientable cases
of $\Gr(k,n)$ via the Poincar\'e-Lefschetz duality.
We also derive  an explicit formula of the Poincar\'e polynomial
for $\Gr(k,n)$ and show that the Poincar\'e polynomial is also related to the number of points
on a suitable version of  $\Gr(k,n)$ over a finite field $\mathbb{F}_q$ with $q$ being a power of a prime. In particular, 
we find that the number of $\mathbb{F}_q$ points on $\Gr(k,n)$ can be computed by counting the number of singularities along the KP flow.
\end{abstract}

\maketitle

\thispagestyle{empty}
\pagenumbering{arabic}\setcounter{page}{1}
\tableofcontents
%\clearpage

%%%

\section{Introduction}

This paper attempts to  extract the topology of  real Grassmannians based on
a {\it realization} of the manifolds related to the KP  hierarchy.  As in the case of the real flag manifolds
discussed in \cite{casian:06} using the Toda lattice hierarchy, we here show that the KP
equation can be used to obtain  similar results for the real Grassmannians. 

The KP equation is a two-dimensional extension of the well-known KdV equation, and it was
introduced by Kadomtsev and Petviashvili in 1970 to study the stability of one KdV soliton under the
influence of weak two-dimensional perturbations \cite{KP:70}. The equation provides also a model to describe a two-dimensional shallow water wave phenomena (see for example \cite{AS:81, YLK:10, K:10, CK:09}.
The KP equation is a dispersive wave equation for the scaler function $u=u(x,y,t)$
with spatial variables $(x,y)$ and time variable $t$, and is given by
\[
\frac{\partial}{\partial x}\left(-4\frac{\partial u}{\partial t}+6u\frac{\partial u}{\partial x}+
\frac{\partial^3 u}{\partial x^3}\right)+3\frac{\partial^2 u}{\partial y^2}=0.
\]
We write the solution in the form 
\[
u(x,y,t)=2\frac{\partial^2}{\partial x^2}\,\ln\tau(x,y,t),
\]
where the function $\tau$ is called the {\it tau} function of the KP equation \cite{S:81, MJD:00}.
It is also well-known that some of the exact solutions can be written in the Wronskian determinant form
(see for example \cite{H:04}),
\begin{equation}\label{tau1}
\tau={\rm Wr}(f_1,f_2,\ldots,f_k):=\left|\begin{matrix}
f_1 & f_2 & \cdots & f_k \\
f_1'&f_2'&\cdots &f_k'\\
\vdots&\ddots&\ddots&\vdots\\
f^{(k-1)}_1&f_2^{(k-1)}&\cdots&f_k^{(k-1)}
\end{matrix}\right|,
\end{equation}
where $f^{(l)}_j:=\partial^lf_j/\partial x^l$, and each of the functions $\{f_i(x,y,t):i=1,2,\ldots,k\}$ satisfies the simple linear equations,
\[
\frac{\partial f_i}{\partial y}=\frac{\partial^2f_i}{\partial x^2},\qquad \frac{\partial f_i}{\partial t}=\frac{\partial^3 f_i}{\partial x^3}.
\]
We then take the following finite dimensional solutions, (i.e.  finite Fourie series),
\[
f_i=\sum_{j=1}^na_{i,j}E_j,\qquad{\rm with}\quad E_j:=\exp\left(\lambda_jx+\lambda_j^2y+\lambda^3_jt\right),
\]
where $A:=(a_{i,j})$ is an $k\times n$ real matrix of rank $k$, and $\lambda_j$'s are real and distinct numbers, say
\[
\lambda_1~<~\lambda_2~<~\cdots~<~\lambda_n.
\]
Since $\lambda_j$'s are all distinct, the set of exponential functions $\{E_j:j=1,\ldots,n\}$ is linearly independent and forms a basis of $\Real^n$.  Then the set of $\{f_i:i=1,\ldots,k\}$ defines a $k$-dimensional
subspace in $\Real^n$.  This implies that each solution defined by the $\tau$-function \eqref{tau1}
can be characterized by the $A$-matrix, hence the solution parametrizes a point of the real Grassmannian.
Note here that the Grassmannian denoted by $\Gr(k,n,\Real)$ can be described as
\[
\Gr(k,n,\Real)={\rm GL}_k(\Real)\setminus M_{k\times n}(\Real),
\]
where $M_{k\times n}(\Real)$ is the set of all $k\times n$ matrices of rank $k$, and ${\rm GL}_k(\Real)$
is the general linear group of $k\times k$ matrices.  We then expect that some particular solutions
of the KP equation contain certain information on the cohomology of $\Gr(k,n,\Real)$.  This is our main
motivation for the present study.  In the previous paper \cite{casian:06}, we found the similar results
for the case of real flag variety using the Toda lattice equation.

\begin{figure}[t!]
\centering
\includegraphics[scale=0.5]{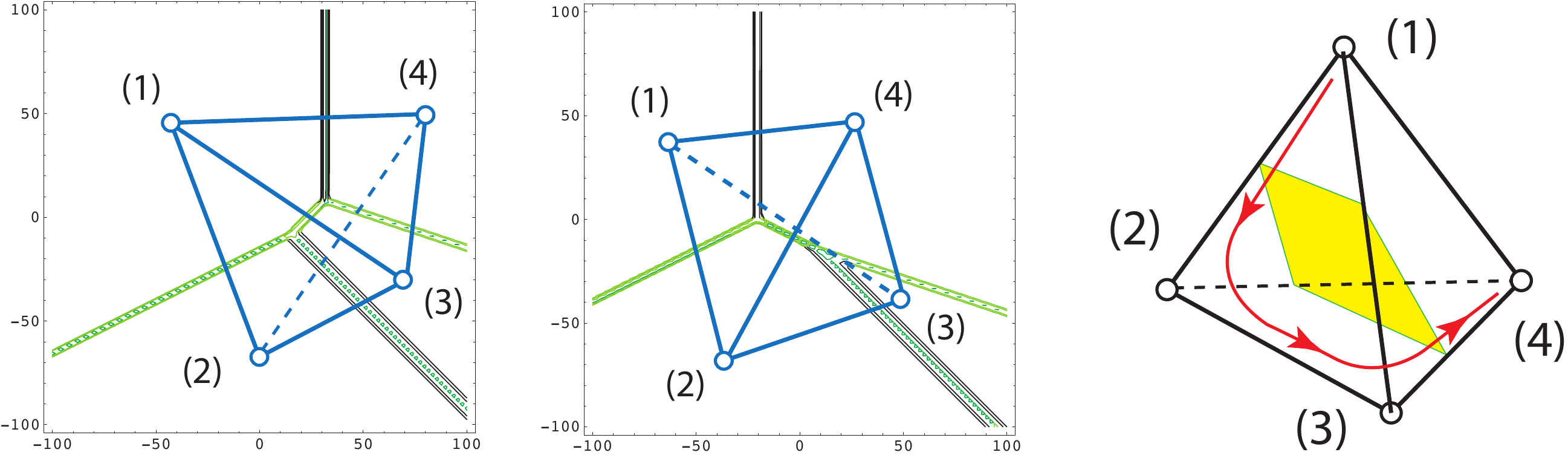}
\caption{The contour plot of the KP solution with the $\tau$-function given by \eqref{tauG14}.
The signs are chosen as $(\epsilon_1,\ldots,\epsilon_4)=(+,-,-,+)$.  In the left two figures,
the number $(i)$ of each region represents the dominant exponential $E_i$, and
the dark (light) colored contour shows the regular (singular) part of the solution.
The left figure shows the contour plot of the solution at $t=-5$ and 
the middle one is at $t=5$.  The tetrahedron in the figures  indicates as a dual graph of the solution pattern, that is, the vertices of the polytope are labeled by
the dominant exponentials. The directed curve in the right figure shows a flow of the KP equation, and
this flow blows up twice by crossing the boundary between the regions with different signs,
which corresponds to $\tau=0$.}
\label{fig:G14}
\end{figure}

Let us explain how some solutions contain  information of the cohomology of $\Gr(k,n,\Real)$ by taking
simple example with $k=1$ and $n=4$.  That is, we recover the cohomology of real variety 
$\Real P^3\cong \Gr(1,4,\Real)$ from a KP flow.  We consider the $\tau$-function in
the form,
\begin{equation}\label{tauG14}
\tau=\epsilon_1\,E_1+\epsilon_2\,E_2+\epsilon_3\,E_3+\epsilon_4\,E_4,
\end{equation}
where $\epsilon_j$ are arbitrary non-zero real constants. Since $|\epsilon_j|$ can be absorbed in the exponential term, we only consider their signs, and here we assume $\epsilon_j\in\{\pm\}$ (which we refer to as {\em KP signs}).
Let us consider the case with the specific signs  $(\epsilon_1,\ldots,\epsilon_4)=(+,-,-,+)$.
Using simple asymptotic argument, one can easily find that each $E_j(x,y,t)$ becomes dominant in certain
region of the $xy$-plane for a fixed $t$. Figure \ref{fig:G14} illustrates the contour plots of the solution $u(x,y,t)$ at $t>0$ and $t<0$.  
 Each region marked by the number $(i)$ indicates the dominant exponential $E_i$, and at the boundary of the region, the $\tau$-function can be expressed by two exponential
functions $E_i$ and $E_j$ (which is dominant in the adjacent region), i.e.  writing $E_i=e^{\theta_i}$
with $\theta_i=\lambda_ix+\lambda_i^2y+\lambda_i^3t$, 
\begin{align*}
\tau&\approx \epsilon_i\,E_i+\epsilon_j\,E_j \\
&=\pm 2 e^{\frac{1}{2}(\theta_i+\theta_j)}\,\left\{\begin{array}{lll}
\displaystyle{\cosh\half(\theta_i-\theta_j)}\quad&{\rm if}\quad \epsilon_i\epsilon_j=+,\\[1.5ex]
\displaystyle{\sinh\half(\theta_i-\theta_j)}\quad&{\rm if}\quad\epsilon_i\epsilon_j=-.
\end{array}\right.
\end{align*}
Then at the boundary, the solution has the following form depending on the sign $\epsilon_i\epsilon_j$,
\begin{align*}
u=\frac{1}{2}(\lambda_i-\lambda_j)^2\,\left\{
\begin{array}{lll}
\displaystyle{{\rm sech}^2\half(\theta_i-\theta_j)}\quad&{\rm if}\quad \epsilon_i\epsilon_j=+,\\[1.5ex]
\displaystyle{{\rm csch}^2\half(\theta_i-\theta_j)}\quad&{\rm if}\quad\epsilon_i\epsilon_j=-.
\end{array}\right.
\end{align*}
In Figure \ref{fig:G14},  the sech-part (regular) of the solution is shown by the dark colored contour, and 
the csch-part (singular) is shown in the light colored contour.
%%%%%%%%%%%%%%%%%%%%%%%%%%%%%%%%%%%%%%%%%%%%%%%%
%%%%%%%%%%%%%%%%%%%%%%%%%%%
We construct a polytope whose vertices are labeled by the dominant exponentials.
The polytope for our example is a tetrahedron as shown in Figure \ref{fig:G14}.
We then consider a KP flow choosing a particular parameter set $(x,y,t)$, so that the flow 
is passing near the edges of the polytope. The flow might be expressed as a graph,
\begin{equation}\label{RP3}
(1)~\longrightarrow~(2)~\Longrightarrow~(3)~\longrightarrow~(4),
\end{equation}
where $\rightarrow$ shows the singular flow, and the $\Rightarrow$ shows the regular one. It is then interesting to note that the graph is the incidence graph of the real projective space $\Real P^3$,
if the arrows are identified as the coboundary operators and their incidence numbers are assigned as $0$ for $\rightarrow$ and $\pm 2$ for $\Rightarrow$. 

%\begin{Remark} 

This kind of graph, exemplified by  graph \eqref {RP3}, was introduced in \cite{casian:06} in an initial attempt to compute the number of connected components carved 
up by the zero divisors of the $\tau$-function, i.e. $\tau_i=0$ for some $i$,  within a moment polytope for the Toda flow. In the case of the KP flow Figure \ref{fig:G14} shows  that  this graph does not contain enough information to attempt to do this (one would have to connect (1) and (4) with an edge $\Rightarrow$).
Still the main results in \cite{casian:06},  which  then motivated  this paper,  were the result of an observation
  that  such a graph, exclusively defined in terms of the flow crossing or not crossing singularities, agrees, in the case of the Toda flow, with the incidence graph of real flag manifolds.   The results
  in this paper extend this type of results  to the KP flow.
 
%\end{Remark}

As in  \cite{casian:06} a  polynomial can be defined  by counting singularities crossed and taking an alternating sum. In this case we
obtain  $1-q^2$, and now
 $q(q^2-1)$, with $q$ a power of a prime number,  counts the number of points over a finite field, ${\mathbb F}_q$, of a variety
closely related to $\Gr(1,4,\Real)$ . This is explained in more detail below and generalized. The cohomology  calculation can be
 thought of as  actually taking place over a field of positive characteristic
(etale cohomology) so that these $q^i$, $i$ counting singularties crossed,   appear as Frobenius eigenvalues in cohomology.

One should note in the graph that the numbers $(i)$ representing the dominant exponential $E_i$  label
the Schubert cells of the variety $\Real P^3$, i.e.
\[
\Real P^3=X(1)\sqcup X(2)\sqcup X(3)\sqcup X(4),
\]
where $X(i)$ is the Schubert cell of dimension $i-1$ (see Section \ref{gras}).  In this paper, we intend to
show that this view of the KP flow can be used to construct the incidence graph of the real
Grassmann variety $\Gr(k,n,\Real)$ in general.

    We simply denote $\Gr(k,n, \Real)$ by $\Gr(k,n)$.  Here we present an overview of the paper.
    Let $\S_n$ be the symmetric group with the simple reflections $s_j$ for
    $j=1,\ldots,n-1$, and $\mathcal{P}_k$ be the parabolic subgroup of $\S_n$ generated by
    $\{s_j:j=1,\ldots,n, j\ne k\}$.  We then have the Schubert decomposition of $\Gr(k,n)$,
    \[
    \Gr(k,n)=\bigsqcup_{w\in \S_n^{(k)}}X_w,
    \]
 where $\S_n^{(k)}$ is the set of minimal coset representatives of $\S_n/ {\mathcal P}_k$ (\cite{bjorner:05}, see also Section \ref{action}).  In this paper, we are interested in the integral cohomology of $\Gr(k,n)$ which we find
 by constructing the corresponding incidence graph (Section \ref{grafica}).
 The co-chain complex may be expressed as follows:  The set of the co-chains are defined by
 \[
 \mathcal{C}^*=\bigoplus_{j=0}^{k(n-k)}\mathcal{C}^j,\qquad{\rm with}\quad
 \mathcal{C}^j=\bigoplus_{j=l(w)}\mathbb{Z}\langle w\rangle,
 \]
 where $l(w)$ is the length of $w\in \S_n^{(k)}$, and $\langle w\rangle:=\overline{X_w}$, the Schubert cycle.
For each $j$, we define an operator $\delta_j:\mathcal{C}^j\to\mathcal{C}^{j+1}$ in the form,
for $\langle w\rangle\in\mathcal{C}^j$ and $\langle w'\rangle\in\mathcal{C}^{j+1}$,
\[
\delta_j\langle w\rangle=\sum_{l(w')=l(w)+1}[w;w']\,\langle w'\rangle.
\]
where $w'=s_iw$ for all $s_i$ with $l(s_iw)=l(w)+1$.  The coefficient $[w,w']$ takes the value either
$0$ or $\pm 2$ (\cite{koch}, see also Section \ref{grafica}).
Then we define a graph $\G(k,n)$, which can be described along the  KP flow, and
which consists of the vertices given by $\langle w\rangle$ for $w\in\S_n^{(k)}$ and
the edges given by the (double) arrows between $\langle w\rangle$ and $\langle w'\rangle$ when $[w,w']\ne0$, that is, there is no blow up.  For example, in the case of $\Gr(1,4)\cong \Real P^3$, we have
$\mathcal{C}^*=\mathbb{Z}[\langle e\rangle]\oplus\mathbb{Z}[\langle s_1\rangle]\oplus\mathbb{Z}[\langle s_2s_1\rangle]\oplus\mathbb{Z}[\langle s_3s_2s_1\rangle]$, and the graph $\G(1,4)$ is given by
\[
\langle e\rangle ~\longrightarrow~ \langle s_1\rangle ~\Longrightarrow~\langle s_2s_1\rangle~\longrightarrow~\langle s_3s_2s_1\rangle.
\]
This graph is the incidence graph for $\Real P^3$ when we assign the coefficients as
$[e;s_1]=0, [s_1;s_2s_1]=\pm 2$ and $[s_2s_1;s_3s_2s_1]=0$ (cf. \eqref{RP3}). This can be extended to the general case, and
 in Section \ref{teorema}, we prove the main theorem,
\begin{Theorem}\label{main1}
The graph $\G(k,n)$ is the incidence graph for the real Grassmannian $\Gr(k,n)$.
\end{Theorem}

The proof of this theorem is based on \cite{casian99}  and  a  simplified description of these results for the real split case
which are  described in terms of the Toda lattice in   \cite{casian:06};  finally an  adaptation of   \cite{casian:06} so that 
incidence graphs are expressed exclusively in terms of the   KP flow,  which, it turns out,  further
simplifies the case of real Grassmanians allowing easy explicit calculation of Betti numbers.

Based on the graph $\G(k,n)$, we then  derive  
explicit formulas for Betti numbers of real Grassmanians, and discuss relations of these  Betti numbers  to  the number of points over a finite field of certain related varieties.  
For this computation, we  first replace  the real Grassmanian manifold $\Gr(k,n)$ 
    with  a complex manifold,   a Zariski open subset  $\Gr(k,n)_{\mathbb C} \subset\Gr(k,n,{\mathbb C})$ with the same homotopy type (see subsection \ref{Fqpoints} for the details).
     The  variety $\Gr(k,n)_{\Complex}$ can be considered over other fields, for example,  a finite field $\F_q$ with $q$ elements. The number of   $\F_q$ points
i.e.  the number of points on $\Gr(k,n)_{\mathbb{F}_q}$, is given by a  polynomial of $q$  denoted by $|\Gr(k,n)_{\mathbb{F}_q}|$.
 In section \ref{Ppolynomial}, we show that the polynomial $|\Gr(k,n)_{\mathbb{F}_q}|$ is related to
 the Poincar\'e polynomial of the real variety $\Gr(k,n)$.
   The  cohomology of  $\Gr(k,n)$
 is much more complicated than the cohomology of the complex variety $\Gr(k,n,{\mathbb C})$,  since it involves  torsion.   However we find that,  the Betti numbers of $\Gr(k,n,{\mathbb C})$ and $\Gr(k',n')$
 for some $(k',n')$ determined from $(k,n),$  agree  after taking into account a change in degrees.
    Recall that the Poincar'e polynomial of the complex variety $\Gr(k,n,\Complex)$, denoted by $P_{(k,n)}^{\Complex}(t)$, 
    is simply obtained  by the number of $\mathbb{F}_q$ points, $|\Gr(k,n, \Complex)|$.  That is, we have 
  \[
P_{(k,n)}^{\Complex}(t)= |\Gr(k,n,\F_q)|_{q=t^2}=\sum_{w\in\S_n^{(k)}}t^{2l(w)}=\left[\begin{matrix} n \\ k \end{matrix}\right]_{t^2},
 \]
Here the polynomial ${\left[\begin{matrix}n\\k\end{matrix}\right]_q}$ is a $q$-analog of the binomial coefficient $\binom{n}{k}$ (see section \ref{gras}).

  The simplest example is the case of $\Gr(1,2)\cong\Real P^1$ which
is a circle  $S^1$  and  can be described with an equation $x^2+y^2=1$.  The number of ${\mathbb F}_q$ points of $\Gr(1,2)$, that is the number of the solutions of $x^2+y^2=1$ for $x,y\in\mathbb{F}_q$,  is $q-1$ (see Remark \ref{points} below).  We thus have a polynomial in $q$, namely,
\[
|\Gr(1,2)_{\mathbb{F}_q}|=q-1, 
\]
and we note that
  this  gives  the Betti numbers $\beta_k$   as coefficients of  $(-1)^{1-j}q^j$, $j=0,1$,  in the polynomial, i.e. $|\Gr(1,2)_{\mathbb{F}_q}|=\sum_{j=0}^n(-1)^{n-j}\beta_jq^j$.
  Namely the Poincar\'e polynomial for $\Gr(1,2)$ denoted by $P_{(1,2)}(t)$ is given by
  \[
  P_{(1,2)}(t)=\sum_{j=0}^1\beta_jt^j=1+t.
  \]
     \begin{Remark} \label{points} In fact, there are two formulas  for the number of points of $x^2+y^2=1$ over ${\mathbb F}_q$ depending on $q$ ($q$ not a power of $2$).  For instance over $\mathbb F_3$ we have $|S^1 ( \mathbb F_3)|=q+1=4$;  however over $F_5$ there  $|S^1;(\mathbb F_5)|=q-1=4$.  By extending $\mathbb F_3$ to  $\mathbb F_{3^2}$ (adding $\sqrt{-1}$),   the  number of ${\mathbb F}_q$ points 
with $q=3^2=9$ becomes $q-1$.  We assume in this paper $\sqrt{-1}\in {\mathbb F}_q$.    
 \end{Remark}

 In a similar computation (see the subsection \ref{Fqpoints}), we obtain the number of $\mathbb{F}_q$ points of $\Real P^3\cong\Gr(1,4)$  as
 \[
 |\Gr(1,4)_{\mathbb{F}_q}|=q(q^2-1).
 \]
The Poincar\'e polynomial for this case is given by
\[
P_{(1,4)}(t)=1+t^3,
\]
and the relation with $|\Gr(1,4)_{\mathbb{F}_q}|$ is not obvious, but we note that the replacement
$q^2-1$ to $1+t^3$ gives the relation between those polynomials.  It turns out that the similar replacement 
can be done for the general cases of $\Gr(2j+1,2m)$ (see the subsection \ref{power}).

  A more interesting example is $\Gr(2,4)$which is discussed in Example \ref {isotropic} below.  
 The polynomial that results in this case is  $q^2(1+q^2)$ giving the number of ${\mathbb F}_q$ points
 on Gr$(2,4)$.  Notice that the polynomial $1+q^2$ is also $\Gr(1,2, {\mathbb F}_{q^2})$  i.e., we have 
 \[
 |\Gr(2,4)_ {{\mathbb F}_{q}}|=q^2|\Gr(1,2, {\mathbb F}_{q^2})| =q^2P_{(1,2)}^{\Complex}(q)=q^2\left[\begin{matrix}2\\1\end{matrix}\right]_{q^2}.
 \]
 The Poincar\'e polynomial is obtained  from $1+q^2$ by  replacing $q$  with $t^2$, and is related to
 $P_{(1,2)}^{\Complex}(q)$  i.e. 
 \[
 P_{(2,4)}(t)= P_{(1,2)}^{\Complex}(t^2)=1+t^4.
\] 
\begin{Remark}
This relation of the Poincar'e polynomials suggests that
the ring of polynomials for the real cohomology for $\Gr(2,4)$ be expressed by
\[
H^*(\Gr(2,4),\Real)\cong \frac{\Real[p]}{\{p^2\}},
\]
where $p$ is the Pontrjagin class $p\in H^4(\Gr(2,4),\Real)$. We discuss this more in section \ref{recursive}.
\end{Remark}

 This paper generalizes  those low dimensional examples  to any real Grassmannians $\Gr(k,n)$  introducing explicit polynomials $p_{(k,n)}(q)$ in $q$ (see the subsection \ref{power}). The polynomial $p_{(k,n)}(q)$ has an interesting connection with the KP flow.  For example,
 in the case of $\Gr(1,4)$ discussed above, the flow blows up twice along the orbit considered (see
 Figure \ref{fig:G14}).  Then we consider the following polynomial
 \[
 p_{(1,4)}(q)=\sum_{w\in\S_4^{(1)}}(-1)^{l(w)}\,q^{\eta(w)}=1-q^2,
 \]
 where $\eta(w)$ counts the number of blow-ups from $e=(1)$ to $w=(i)=s_{i-1}\cdots s_1$, i.e. $\eta(e)=0, ~\eta(s_1)=\eta(s_2s_1)=1$ and $\eta(s_3s_2s_1)=2$. 
Namely,  we also  assign
powers $q^{\eta(w)}$ to the vertex labeled by $w\in\S_n^{(k)}$  which  keep track of the number of  blow-ups
from $e$ to $w$.  Then each vertex in the graph $\G(k,n)$ would have additional label $q^{\eta(w)}$,
which defines the {\it weighted} Schubert cell (see subsection \ref{power}).  In the case of $\Gr(1,4)$, we have the graph with the weighted Schbert cells $(\sigma_i)=(i)$ (see Figure \ref{fig:G14}),
\[
((1); q^0)~\longrightarrow~((2);q^1)~\Longrightarrow~((3);q^1)~\longrightarrow~((4);q^2).
\]
The polynomial $p_{(k,n)}(q)$ appears in $|\Gr(k,n)_{\mathbb{F}_q}|$ given above, and in section \ref{Fqpoints}, we show that this is true in general.
  Without technicalities, the main result is that these polynomials $p_{(k,n)}(q)$ times  certain power of 
 $q$ compute the number of  ${\mathbb F}_q$ points.  The polynomial in $t$
 given by   $p_{(k,n)}(t^2)$ is  the corresponding Poincar\'e polynomial. 
The cases of the form $(k,n)=(2j+1, 2m)$ turn out to be  exceptions obeying  a slightly different   formula. Excluding the exceptional cases  $(k,n)=(2j+1, 2m)$,  
 the relation between the complex Grassmannians $\Gr(j,m,{\mathbb C})$ and the
 real Grassmannians $\Gr(k,n)$ with  $k=2j, 2j+1$, $n=2m,2m+1$ is given by     $p_{(k,n)}(q^{1/2})$  vs.   $p_{(k,n)}(q)$. 
 
 %Strictly speaking the {\em real} cases are just Zariski open subsets of the {\em complex} cases  $\Gr(k,n,{\mathbb C})$ and the  cohomology of  the real Grassmannians
 %is much more complicated involving torsion but the Betti numbers of $\Gr(j,m,{\mathbb C})$ and $\Gr(k,n)$ have a very simple relationship.

 %The following theorem then gives  the number of ${\mathbb F}_q$ points in the subvariety $\Gr(k,n)_{\mathbb F}$ of $ \Gr(k,n, {\mathbb F})$.  

We then obtain the following Theorems in Section \ref{Ppolynomial}:
 
\begin{Theorem}  \label{A} The number of ${\mathbb F}_q$ points of $\Gr(k,n)$ is  given as follows:  
 
 \begin{enumerate}
 
 \item [(i)]  If   $(k,n)$ equals $(2j, 2m), (2j, 2m+1)$ or $ (2j+1, 2m+1)$, we have 
 \[
   | \Gr(k, n)_{{\mathbb F}_q} |= q^r  \left[\begin{matrix} m \\ j \end{matrix}\right]_{q^2}\quad {\rm with}\quad r=k(n-k)-2j(m-j).
   \]
 \item [(ii)]   If $(k,n)=(2j+1,2m)$, then we have
 \[
 | \Gr(k,n)_{{\mathbb F}_q} |= q^r  \left[\begin{matrix} m-1 \\ j \end{matrix}\right]_{q^2}(q^m-1)
 \quad {\rm with}\quad r=k(n-k)-2j(m-j-1)-m.
 \]  
 \end{enumerate}
  \end{Theorem}
 
One should note here that  $ \left[\begin{matrix} m \\ j \end{matrix}\right]_{q^2}=| \Gr(j, m , {\mathbb F}_{q^2} )|$ for any $m\ge j\ge1$.
Let $P_{(k,n)}(t)$ denote the Poincar\'e polynomial of $\Gr(k,n)$, i.e.
\[
P_{(k,n)}(t)= \sum_{i=0}^{k(n-k)} \beta_i\,t^i\quad{\rm with}\quad \beta_i:={\rm dim}\, H^i(\Gr(k,n), {\mathbb Q} ).
\] 
We also show that those $\mathbb{F}_q$ points are  related to the Poincar\'e polynomials of $\Gr(k,n)$:
 \begin{Theorem} \label{B}  The Poincar\'e polynomial of $\Gr(k,n)$ is given as follows:
 \begin{itemize}
   \item [(i)]  If   $(k,n)$ equals $(2j, 2m), (2j, 2m+1)$ or $(2j+1, 2m+1)$, we have 
   \[
   P_{(k,n)}(t)= \left[\begin{matrix} m \\ j \end{matrix}\right]_{t^4} .
 \]
   \item [(ii)]  If $(k,n)=(2j+1,2m)$, we have
 \[
 P_{(k,n)}(t)= (1+t^{2m-1} ) \left[\begin{matrix} m-1 \\ j \end{matrix}\right]_{t^4}.
 \]
\end{itemize}  
  \end{Theorem}
This then gives the well-known formulas on the Euler characteristic for $\Gr(k,n)$.
\begin{Corollary}
The Euler characteristic $\chi(\Gr(k,n))$ has the following form.
\begin{itemize}
\item[(i)] If $(k,n)$ equals $(2j,2m), (2j,2m+1)$ or $(2j+1,2m+1)$, we have
\[
\chi(\Gr(k,n))=P_{(k,n)}(-1)=\binom{m}{j}.
\]
\item[(ii)] If $(k,n)=(2j+1,2m)$, then we have
\[
\chi(Gr(k,n))=P_{(k,n)}(-1)=0.
\]
\end{itemize}
\end{Corollary}
It is interesting to note in particular that for the case (i), we have 
  \[
  P_{(k,n)}(t)=P^{\Complex}_{(\lfloor{k}/{2}\rfloor,\lfloor{n}/{2}\rfloor)}(t^2),
  \]
which gives an information on the structure of the cohomology ring $H^*(\Gr(k,n),\Real)$
(see \cite{BT:82, MS:74}).

  Let us summarize the results for low dimensional cases, $\Gr(1,2), \Gr(1,3)$ and  $\Gr(2,4)$.
\begin{itemize}
\item[(a)]  The case $\Gr(1,2)\cong S^1$:  From Theorem \ref{A} with $(k,m)=(1,2)$, 
 we recover the polynomial
 \[
 |\Gr(1,2)_{\F_q}|=(q-1)\left[\begin{matrix} 1 \\ 0 \end{matrix}\right]_{q^2}=q-1.
 \]   
 The Poincar\'e polynomial is then obtained from Theorem \ref{B} as
 \[
 P_{(1,2)}(t)=(1+t)\left[\begin{matrix} 1 \\ 0 \end{matrix}\right]_{t^4}=1+t.
 \]
 Notice in general that $\Gr(1,n)\cong \Real P^{n-1}$ and $P_{(1,n)}(t)=1+t^{2m-1}$ when $n=2m$.
\item[(b)]  The case $\Gr(1,3)\cong \Real P^2$:  We have
\[
|\Gr(1,3)_{\F_q}|=q^2\left[\begin{matrix} 1 \\ 0 \end{matrix}\right]_{q^2}=q^2.
\]
 The Poincar\'e polynomial is
 \[
 P_{(1,3)}(t)=1.
 \]
 that is, the $\Gr(1,3)$ is connected but not orientable.  In general, we have $P_{(1,n)}(t)=1$ when $n$ is odd.
 \item[(c)] The case $\Gr(2,4)$:  We have
\[
|\Gr(2,4)_{\F_q}|=q^2\left[\begin{matrix} 2 \\ 1 \end{matrix}\right]_{q^2}=q^2(1+q^2).
\]
 The  Poincar\'e polynomial is then given by
 \[
 P_{(2,4)}(t)=\left[\begin{matrix}2\\1\end{matrix}\right]_{t^4}=1+t^4,
 \]
 which shows that the manifold is connected and orientable.
 \end{itemize}

   %%%%%%%%%%
 
\section{The real Grassmann manifold $\Gr(k,n)$}  \label{gras}

Let $G=\SL_n(\Real)^{\pm}$, the set of $n\times n$ real matrices of determinant $\pm1$.
%We denote by $\Gr(k,n, \Real)$ or just $\Gr(k,n)$,  the real Grassmannian consisting of $k$-dimensional subspaces of $\Real^n$. 
We have $\Gr(k,n)=G/P_k$ for a maximal parabolic subgroup $P_k$ stabilizing the $k$-dimensional subspace $V_0(k)={\rm Span}_{\Real}\{E_1,\ldots,E_k\}$ with the standard basis $\{E_j:j=1,\ldots,n\}$. In this section, we give an elementary introduction of the Grassmann manifolds
and the Schubert decompositions.

The vector spaces  $\Real^n= V(E_1,\dots E_n)$ and $\bigwedge^k \Real^n=\bigwedge^k V(E_1,\dots, E_n)$  for $k=1,\ldots, n$ are also the {\it fundamental} representations of $G$ with the standard action of $G$ on $V(E_1,\dots E_n)$.  Note that the real Grassmannian can be also expressed as a (unoriented) homogeneous space, 
\[
\Gr(k,n)\cong \frac{{\rm O}_n(\Real)}{{\rm O}_k(\Real)\times {\rm O}_{n-k}(\Real)}.
\]

\subsection{The Schubert decomposition}
Let $\{E_j:j=1,\ldots,n\}$ be a basis of $\mathbb{R}^n$. Then a $k$-dimensional vector subspace
is expressed by a set of vectors $\{f_i:i=1,\ldots,k\}$ with
\begin{equation}\label{fk}
f_i=\sum_{j=1}^n a_{i,j}E_j\qquad i=1,\ldots, k\,,
\end{equation}
With the $k\times n$ real matrix $A:=(a_{i,j})$,
\[
A=\begin{pmatrix}
a_{1,1} & a_{1,2} & \cdots  & a_{1,n} \\
\vdots  & \vdots   &  \ddots & \vdots \\
a_{k,1} & a_{k,2} & \cdots  & a_{k,n}
\end{pmatrix}
\in M_{k\times n}(\mathbb{R})\,,
\]
we write
\[    
(f_1,\ldots,f_k)^T=A\cdot (E_1,\ldots, E_n)^T\,.
\]
Here $M_{k\times n}(\Real)$ is the set of $k\times n$ matrices of the maximal rank $k$.
Each $A$ matrix gives a label of a point of $\Gr(k,n)$, and it can be put in a unique form
called the row reduced echelon form (RREF): For some $g\in {\rm GL}_k(\mathbb{R})$,
\[
g\cdot (f_1,\ldots,f_k)^T=g\cdot A\cdot (E_1,\ldots,E_n)^T\,,
\]
expresses the same $k$-dimensional subspace. Choosing an appropriate matrix $g$,
one can put $A$ in the RREF. This implies
\[
\Gr(k,n) \cong {\rm GL}_k(\Real)\backslash M_{k\times n}(\mathbb{R}).
\]
To each RREF we can associate a set of ordered integers $(\sigma_1,\ldots,\sigma_k)$, $\sigma_j$   indicating that the  $\sigma_j$-th column has a  pivot.  This parametrizes the cells   in a cell decomposition of
 $\Gr(k,n)$ which is called the Schubert decomposition,
\[
\Gr(k,n)=\bigsqcup_{1\le \sigma_1<\cdots<\sigma_k\le n}\,X(\sigma_1,\ldots,\sigma_k)\,,
\]
where $X(\sigma_1,\ldots,\sigma_k)$ is the Schubert cell defined by
the set of all matrices whose pivot ones are at $(\sigma_1,\ldots,\sigma_k)$ places.
The ordered set $(\sigma_1,\ldots,\sigma_k)$ is called the Schubert symbol
for the corresponding Schubert cell, and we sometime denote the cell simply by
the symbol.  
\begin{Example} \label{Gr24:example} The Grassmannian
$\Gr(2,4)$, the space consisting of two dimensional vector spaces
inside the four dimensional space $V(E_1,\ldots,E_4)$, has six Schubert cells, $X(i,j)$ with $1\le i<j\le 4$, which are given by
\[
\begin{array}{llll}
\displaystyle{X(1,2)=\left\{\begin{pmatrix} 
1 &0 & 0 & 0 \\
0 & 1 & 0& 0
\end{pmatrix}\right\}} \quad & 
X(1,3)=\left\{\begin{pmatrix} 
1 &0 & 0 & 0 \\
0 & * & 1& 0
\end{pmatrix}\right\} \\[3.0ex]
\displaystyle{X(2,3)=\left\{\begin{pmatrix} 
* &1 & 0 & 0 \\
* &0 & 1& 0
\end{pmatrix}\right\}} \quad & 
X(1,4)=\left\{\begin{pmatrix} 
1 &0 & 0 & 0 \\
0 & * & *& 1
\end{pmatrix}\right\} \\[3.0ex]
\displaystyle{X(2,4)=\left\{\begin{pmatrix} 
* &1 & 0 & 0 \\
* & 0 & *& 1
\end{pmatrix}\right\}} \quad & 
X(3,4)=\left\{\begin{pmatrix} 
*&* & 1& 0 \\
* & * & 0& 1
\end{pmatrix}\right\}
\end{array}
\]
Each $``*"$ indicates a free parameter for RREF, and the number of free parameters gives the dimension
of the cell.
  \end{Example}

The  dimension of the cell is given by 
\[
{\rm  dim}\, X(\sigma_1,\ldots,\sigma_k)=\sum_{j=1}^k(\sigma_j-j)\,,
\]
%where $Y_{\lambda}$ is the Young diagram associated to $X(\sigma_1,\ldots,\sigma_k)$ for
%the partition 
%$\lambda=(\lambda_1,\ldots,\lambda_k)$ with $\lambda_j=\sigma_{k-j+1}-(k-j+1)$, (note
%$\lambda_j\ge \lambda_{j+1}$ as the usual definition of the Young diagram,
%and each $\lambda_j$ expresses the number of boxes in the horizontal direction). 
Let $n_j$ be the total
number of cells of $j$ dimension. Then
the generating function $\sum_{j=0}^{k(n-k)}n_j q^{j}$ 
is given by counting the number of points in $\Gr(k,n,\mathbb{F})$ over a finite
field $\mathbb{F}=\mathbb{F}_q$, i.e.
\begin{equation}\label{qbinom}
\sum_{j=0}^{k(n-k)}n_j q^{j}=|\Gr(k,n,\mathbb{F}_q)|.
\end{equation}
The number $|\Gr(k,n,\mathbb{F}_q)|$ can be obtained from the complex
Grassmannian $\Gr(k,n,\Complex)$ restricted on $\mathbb{F}_q$, i.e.
\[
|\Gr(k,n,\mathbb{F}_q)|=\left[\begin{matrix} n \\ k \end{matrix}\right]_{q}=\frac{[n]_q!}{[k]_q!\,[n-k]_q!}\,,
\]
where $[n]_q!=[n]_q[n-1]_q\cdots[2]_q[1]_q$ with $[j]_q=(1-q^j)/(1-q)$, the $q$-analog of
the number $j$.
%and $\sum_{\lambda}$ indicates the sum over all Young subdiagrams of the rectangular
%diagram with the size $k\times (n-k)$.  
For example, we have
\[
|\Gr(2,4,\F_q)|=\left[\begin{matrix}4\\2\end{matrix}\right]_{q}=1+q+2q^2+q^3+q^4,
\]
(see Example \ref{Gr24:example}).

We also note that the polynomial $|\Gr(k,n,\mathbb{F}_q)|$ with $q=t^2$ is the Poincar\'e polynomial for the cohomology of $\Gr(k,n,\mathbb{C})$, denoted by $P_{(k,n)}^{\Complex}(t)$, that is,
\[
P_{(k,n)}^{\Complex}(t)=\sum_{j=0}^{k(n-k)}\,\beta_{2j}^{\Complex}t^{2j}\,,
\]
and each coefficient $\beta_{i}^{\Complex}$ gives the Betti number for $i$-cycles,
\[\beta^{\Complex}_i:=
 {\rm dim}\,H^{i}(\Gr(k,n,\mathbb{C}); \Real) =\left\{\begin{array}{lll}
 0 \qquad &{\rm if}\quad i={\rm odd}\\
n_j\qquad &{\rm if}\quad i=2j.
 \end{array}\right.
 \]
 
%%correction

\subsection{The Pl\"ucker embedding and the moment polytope}
In order to describe  the points of $\Gr(k,n)$, we can use the Pl\"ucker embedding,
\[\begin{array}{ccccc}
\Gr(k,n)  & \longrightarrow  & {P}\left( \bigwedge^k\,\mathbb{R}^n\right) \\[1.0ex]
[f_1,\ldots,f_k] &\longmapsto & {\mathbb R} f_1\wedge \cdots\wedge f_k
\end{array}
\]
where the wedge product can be expressed as
\begin{equation}\label{kframe}
f_1\wedge \cdots\wedge f_k=\sum_{1\le j_1<\cdots<j_k\le n}\xi(j_1,\ldots,j_k)\,E_{j_1}\wedge\cdots\wedge E_{j_k}\,,
\end{equation}
Here the coefficients $\xi(j_1,\ldots,j_k)$ are the Pl\"ucker coordinates given by the $k\times k$ minor
of the coefficient matrix $A$, i.e.
\begin{equation}\label{plucker}
\xi(j_1,\ldots,j_k)=\left| \left(a_{l,j_m}\right)_{1\le l,m\le k}\right|\,.
\end{equation}
Those coordinates satisfy the Pl\"ucker relations,
\begin{equation}\label{pluckerrelation}
\sum_{r=1}^{k+1}\,(-1)^{r}\,\xi(\alpha_1,\ldots,\check{\alpha}_r,\cdots,\alpha_{k+1})\,\xi(\alpha_r,\beta_1,\ldots,\beta_{k-1})=0,
\end{equation}
for any set of numbers $\{\alpha_1,\ldots,\alpha_{k+1},\beta_1,\ldots, \beta_{k-1}\}\subset\{1,\ldots,n\}$. Here the notation 
$\check{\alpha}_r$ implies the deletion of this number.

\begin{Example}\label{isotropic}
In the case of $\Gr(2,4)$, there is one Pl\"ucker relation,
\begin{equation}\label{PR24}
\xi(1,2)\xi(3,4)-\xi(1,3)\xi(2,4)+\xi(1,4)\xi(2,3)=0.
\end{equation}
Then $\Gr(2,4)$ can be expressed as
\[
\Gr(2,4)=\left\{[v]=\Real\sum_{1\le i,j\le 4}\xi(i,j)\,E_i\wedge E_j: \begin{array}{ccc}
\xi(i,j)~{\rm satisfy}~(\ref{PR24})\\[1.5ex]
\sum_{1\le i<j\le 4} |\xi(i,j)|^2\ne 0\end{array}\right\}.
\]
This variety can be also expressed as the following, which may be useful for counting
the $\F_q$ points of $\Gr(2,4)$.
Let $(x_1,x_2,x_3,y_1,y_2,y_3)$ be a new coordinate defined by
\[\begin{array}{llll}
\xi(1,2)=x_1-y_1, & \xi(1,3)=-x_2+y_2, & \xi(1,4)=x_3-y_3, \\[1.5ex]
\xi(3,4)=x_1+y_1, & \xi(2,4)= x_2+x_2, & \xi(2,3)=x_3+y_3.
\end{array}
\]
Then one can express $\Gr(2,4)$ as the set of {\it isotropic} vectors in $\Real P^5$,
\begin{equation}\label{GR24}
\Gr(2,4)=\left\{\Real (x_1,x_2,x_3,y_1,y_2,y_3):  \begin{array}{cccc}
x_1^2+x_2^2+x_3^3-y_1^2-y_2^2-y_3^2=0 \\[1.5ex]
x_1^2+x_2^2+x_3^2\ne 0 \end{array}
\right\}  \,\subset \Real P^5\,.
\end{equation}
Here the first condition is the Pl\"ucker relation (\ref{PR24}) and the second one shows
the maximal rank (i.e. at least one $\xi(i,j)\ne 0$).
\end{Example}

With the expression of the wedge product $f_1\wedge\cdots\wedge f_k$, one can define
the moment map $\mu:\Gr(k,n)\to \mathfrak{h}^*_{\mathbb{R}}$ (see \cite{gelfand:87}),
\begin{equation}\label{moment}
\mu(f_1\wedge\cdots\wedge f_k)=\frac{\sum_{1\le i_1<\cdots<i_k\le n}|\xi(\sigma_1,\ldots,\sigma_k)|^2\,(L_{\sigma_1}+\cdots+L_{\sigma_k})}
{\sum_{1\le \sigma_1<\cdots<\sigma_k\le n}|\xi(\sigma_1,\ldots,\sigma_k)|^2}
\end{equation}
where  $L_\sigma$ are the weights
of the standard representation $V$ of $\SL_n(\Real)$, and $\mathfrak{h}^*_{\Real}$ is
defined by
\[
\mathfrak{h}^*_{\mathbb{R}}={\rm Span}_{\mathbb{R}}\left\{L_1,\ldots,L_k:\sum_{i=1}^n L_i=0\right\}~\cong \Real^{n-1}\,.
\]
The image under the moment map is a convex polytope whose vertices
are marked by the elements of $\mathcal{S}^{(k)}_n$. In the representation theory, this polytope is a weight polytope
of the fundamental representation of ${\mathfrak{sl}}(n)$ on $\bigwedge^k V(E_1,\ldots,E_n)$. 
The set $\{E_{\sigma_1}\wedge\cdots\wedge E_{\sigma_k}: 1\le \sigma_1<\cdots<\sigma_k\le n\}$ then corresponds to a basis of weight vectors
of the representation $\bigwedge^kV$ (see for example \cite{fulton:91}).

\begin{figure}[t!]
\centering
\includegraphics[scale=0.35]{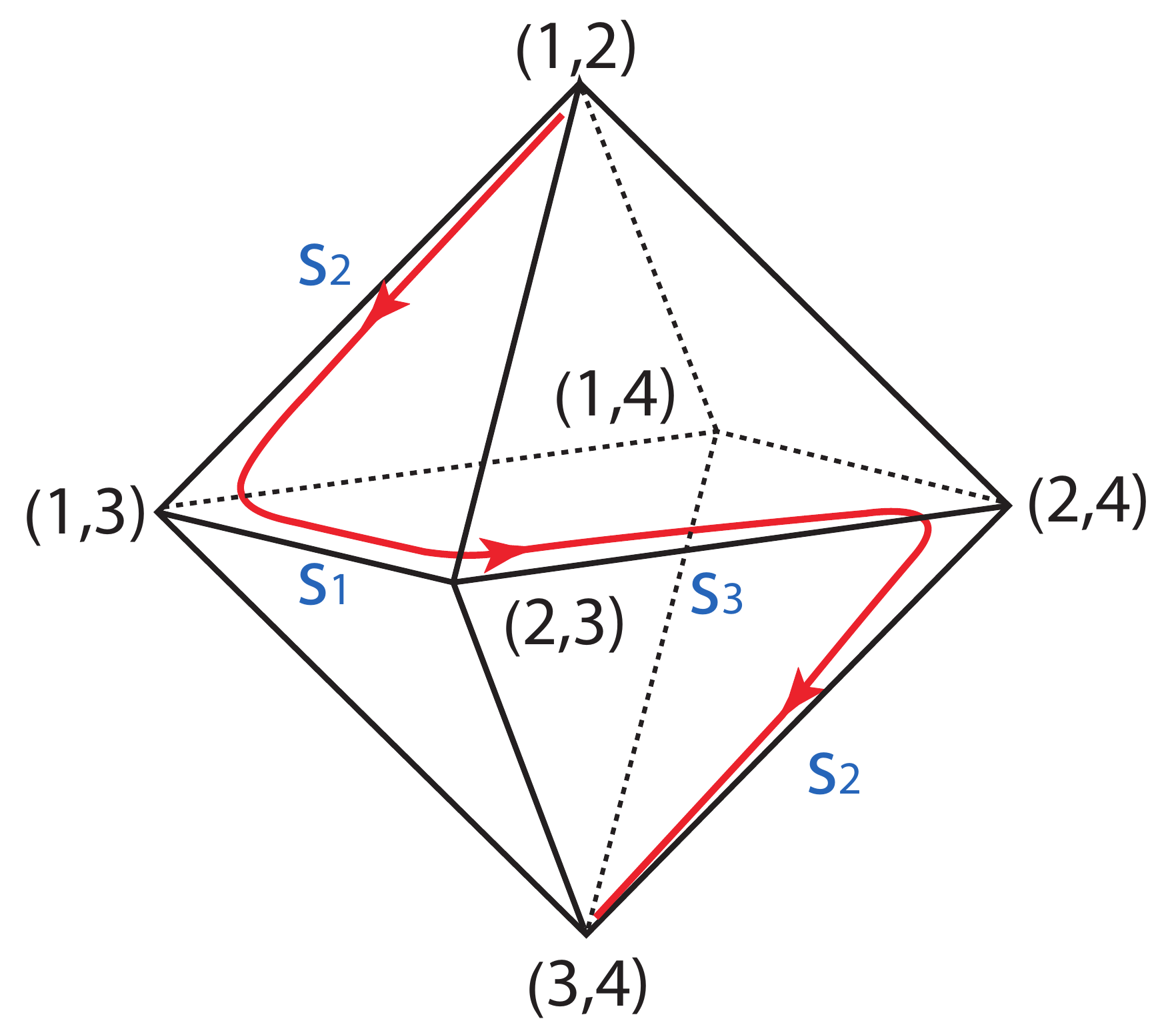}
\caption{The moment polytope of $\Gr(2,4)$.
Each weight is marked by $(\sigma_1,\sigma_2)$ with
$L=L_{\sigma_1}+L_{\sigma_2}$. The highest weight is given by
$(1,2)$, and all the weights are in the $\S_4^{(2)}$-orbit of highest weight.
The arrowed curve indicates a KP flow, which approximates the Weyl group action.
Each $s_j$ shows the simple reflection $j\leftrightarrow j+1$.}
\label{fig:Grassmann24}
\end{figure}

In Figure \ref{fig:Grassmann24},
we illustrate the moment polytope of $\Gr(2,4)$, which is an octahedron whose vertices are labeled
by the Schubert symbols $(\sigma_1,\sigma_2))$.
%%%%%%%%%%%%%%%%%%%%%%%%%%%%%%%%%%%%%%%%%%%%%%%%

%%%%%%%%%%%%%%%%%%%%%%%%%%%

\section{The KP flow on $\Gr(k,n)$} \label{kp}
In this section, we give a realization of $\Gr(k,n)$ in terms of finite dimensional solutions
of the KP equation (see, for example, \cite{CK:08, CK:09, K:10} for the recent development of the classification problem of
soliton solutions of the KP equation). The purpose of this section is to give a dynamical system, called a KP flow,
 on the moment polytope (\ref{moment}). In particular, we consider the flows which approximate Weyl group action
 along the weak Bruhat order (see section \ref{action}).  For example, a KP flow in $\Gr(2,4)$ is 
  illustrated in Figure \ref{fig:Grassmann24}, which gives the orbit $(1,2)\overset{s_2}{\to}(1,3)\overset{s_1}{\to}(2,3)\overset{s_3}{\to}(2,4)\overset{s_2}{\to}(3,4)$ with the simple reflections $s_j\in \mathcal{S}_4$, the
  symmetric group of permutations.  Those flows contain the information of the cohomology of the real
  variety $\Gr(k,n)$, which is the main motivation of the present paper.

\subsection{The KP flow on the moment polytope} \label{dynamical}
Let us start to fix the set of independent vectors $\{E_i: i=1,\ldots,n\}$ as 
given by the exponential functions of multi-variables $\mathbf{t}=(t_1,t_2,\ldots, t_{n-1})\in\Real^{n-1}$,
\[
E_i(\mathbf{t})=e^{\theta(\lambda_i;\mathbf{t})}\qquad{\rm with}\quad \theta(\lambda_i;\mathbf{t})=\sum_{r=1}^{n-1}\lambda_i^rt_r\,,
\]
where $\lambda_i\in\Real$ are all distinct, and we assume the ordering,
\[
\lambda_1~<\lambda_2~<~\cdots~<\lambda_n.
\]
With the set of functions $\{f_j: j=1,\ldots,k\}$ defined in (\ref{fk}), i.e.
$(f_1,\ldots,f_k)^T=A\,(E_1,\ldots,E_n)^T$ with the $k\times n$ coefficient matrix $A$, we define the following function, called the $\tau$-function, given by the Wronskian determinant of those
functions,
\begin{align}\label{tauA} 
\tau_A(\mathbf{t})&= {\rm Wr}(f_1,\ldots,f_k) 
  = \left|\begin{matrix}
  f_1 & \cdots & f_k \\
  \vdots &\ddots&\vdots\\
  f_1^{(k-1)}&\cdots&f_k^{(k-1)}
  \end{matrix}\right|\,,
\end{align}
Here the index $A$ indicates the matrix $A$ defined in (\ref{fk}). We then have:
\begin{Lemma}\label{BinetCauchy}
The $\tau$-function can be expressed in the form,
\[
\tau_A(\mathbf{t})=\sum_{1\le \sigma_1<\cdots<\sigma_k\le n}\xi(\sigma_1,\ldots,\sigma_k)\,E(\sigma_1,\ldots,\sigma_k;\mathbf{t})\,,
\]
where $\xi(\sigma_1,\ldots,\sigma_k)$ are the Pl\"ucker coordinates given by (\ref{plucker}), and
\[
E(\sigma_1,\ldots,\sigma_k;\mathbf{t})={\rm Wr}(E_{\sigma_1},\ldots,E_{\sigma_k})=\prod_{i<j}
(\lambda_{\sigma_j}-\lambda_{\sigma_i})\, \exp \theta(\sigma_1,\ldots,\sigma_k;\mathbf{t})\,,
\]
with  $\theta(\sigma_1,\ldots,\sigma_k;\mathbf{t}):=\sum_{i=1}^k\theta(\lambda_{\sigma_i};\mathbf{t})$.
\end{Lemma}
This Lemma is a direct consequence of the Binet-Cauchy theorem (see p.9 in \cite{G:59}).
With the formula for $f_1\wedge\cdots\wedge f_k$ in (\ref{kframe}), the $\tau_A$-function
gives a realization of the corresponding point on $\Gr(k,n)$ with the identification
$\tau_A\equiv c \tau_A$ for any nonzero constant $c$ (the projectivization). 
Namely, the Wronskian-map can be considered as the Pl\"ucker embedding,
${\rm Wr}:\Gr(k,n)\hookrightarrow \Real P^{\binom{n}{k}-1}$.
Here one should note that $\sum_{i=1}^k\lambda_{\sigma_i}$ should be distinct
for all $(\sigma_1,\ldots,\sigma_k)$. Then one can identify the set $\{E(\sigma_1,\ldots,\sigma_k;\mathbf{t})\}$ 
as the basis $\{E_{\sigma_1}\wedge \cdots\wedge E_{\sigma_k}\}$ of $\wedge^k\mathbb{R}^n$.
 
One can also define a map $\varphi_A:\Real^{n-1}\to \mathfrak{h}^*_{\Real}$ as the composite map,
$\varphi_A=\mu\circ \tau_A$, with
\[
\varphi_A(\mathbf{t})=\frac{\sum_{1\le \sigma_1<\cdots<\sigma_k\le n}|\xi(\sigma_1,\ldots,\sigma_k)
E(\sigma_1,\ldots,\sigma_k;\mathbf{t})|^2\,(L_{\sigma_1}+\cdots+L_{\sigma_k})}
{\sum_{1\le \sigma_1<\cdots<\sigma_k\le n}|\xi(\sigma_1,\ldots,\sigma_k)E(\sigma_1,\ldots,\sigma_k;\mathbf{t})|^2}.
\]
Then if $A$ is in the top cell, one can see that 
the closure of the set $\{\varphi(\mathbf{t}):\mathbf{t}\in\Real^{n-1}\}$ is a convex hull with
the vertices given by the weights $\{L_{\sigma_1}+\cdots+L_{\sigma_k}:1\le \sigma_1<\cdots<\sigma_k\le n\}$, denoted by Conv$(k,n)$, i.e.
\[
{\rm Conv}(k,n)=\overline{\left\{\varphi_A(\mathbf{t}):\mathbf{t}\in\Real^{n-1}\right\}}\,.
\]
This map defines a flow on the polytope, and we are interested in studying the zero set
of the $\tau_A$ function which is determined by a particular choice of the coefficient matrix $A$.
\begin{Example}
Let us consider the case $\Gr(2,4)$ to show how one can see the polytope in Figure \ref{fig:Grassmann24}
from the KP flow.  Figure \ref{fig:Ttype} illustrates the contour plots of the $t$-evolution of the solution 
$u(x,y,t)$.  Then as we explained in the introduction, the polytope can be seen from the dual graph of
the solution pattern as shown in the figure.  Here we take the matrix $A$ in the form,
\[
A=\begin{pmatrix} 
1 & 1 & 1 & 1 \\
\lambda_1 &\lambda_2 &\lambda_3 & \lambda_4
\end{pmatrix}.
\]
Note here that the $2\times 2$ minors of $A$ are all positive, i.e. $\xi(i,j)=\lambda_j-\lambda_i>0$ with
the ordering in $\lambda_j$'s, and this implies the $\tau$-function is positive definite and no blow-up
in the KP flow.  We then extend this example to the case with some signs in the Pl\"ucker coordinates, and show that the flow associated
with a particular choice of the signs can provide topological information of the real Grassmannians. 
%%%%%%%%%%%%%%%%%%%%%%%%%%%%%%%%
\begin{figure}[t!]
\centering
\includegraphics[scale=0.5]{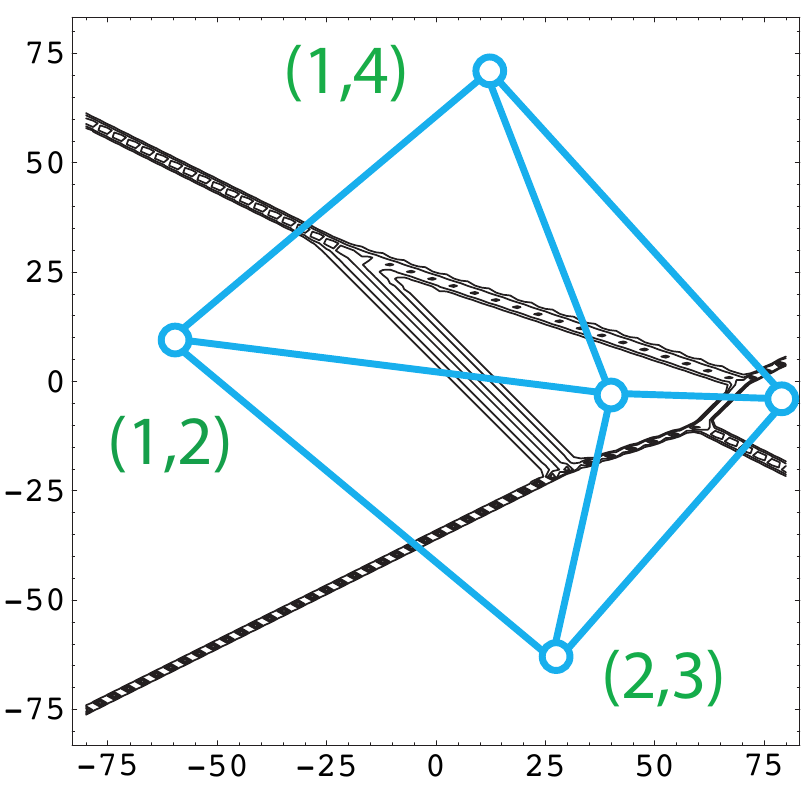} \hskip 0.5cm 
\includegraphics[scale=0.5]{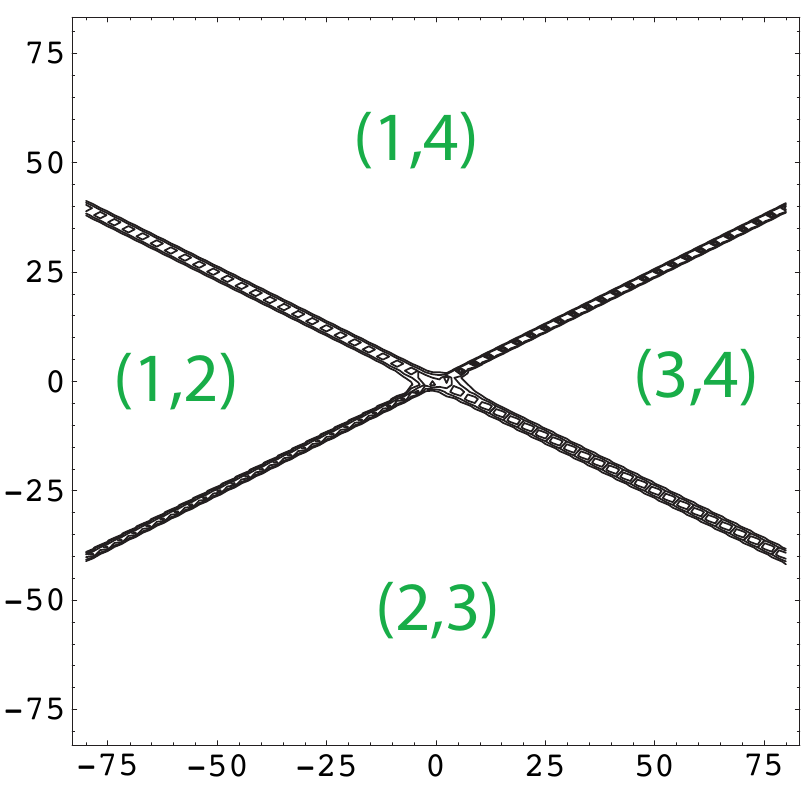} \hskip 0.5cm
\includegraphics[scale=0.5]{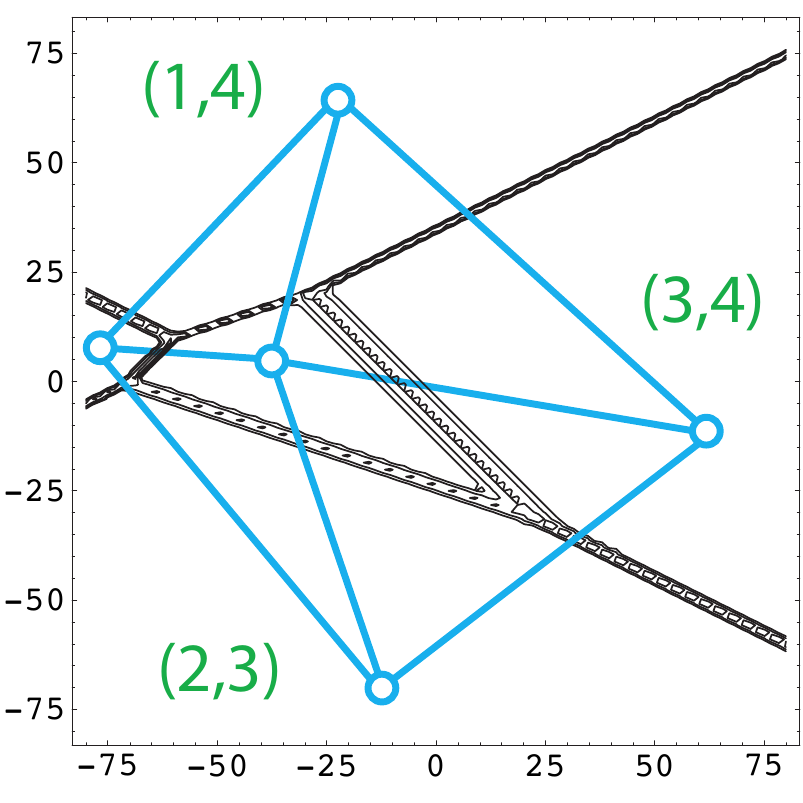}
\caption{Soliton solution of the KP equation on $\Gr(2,4)$.
The figures show the contour plot of the solution $u(x,y,t)$ for 
$t<0$ (left), $t=0$ (middle) and $t>0$ (right).
Each pair of numbers $(i,j)$ indicates the dominant exponential $E(i,j)$, i.e. the Pl\"ucker coordinates.
The bounded region in the left (right) figure corresponds to $E(1,3)$ ($E(2,4)$).  The light colored graph
in the figures shows the dual graph whose vertices represent those dominant exponentials, and
is the moment polytope (tetrahedron) shown in Figure \ref{fig:Grassmann24}.}
\label{fig:Ttype}
\end{figure}
%%%%%%%%%%%%%%%%%%%%%%%%
\end{Example}

\subsection{The KP equation and the Toda lattice}
As we explained in Introduction, the function 
\[
u(x,y,t)=2\frac{\partial^2}{\partial x^2}\ln\tau_A(x,y,t),
\]
satisfies the KP equation with the identification $t_1=x,t_2=y$ and $t_3=t$.
The higher times, $t_j$ for $j>3$, give the symmetry parameters of the KP equation,
and the set of all the flows parametrized by $t_j$'s forms the so-called KP hierarchy  (see for example \cite{MJD:00, CK:09, K:10}). 

In this paper, we consider a particular set of independent functions $\{f_i:i=1,\ldots,k\}$ such that
\[
f_i=\frac{\partial^{i-1}}{\partial t_1^{i-1}}f=: f^{(i-1)},\qquad {\rm with}\quad f=\sum_{j=1}^n\epsilon_jE_j,
\]
where we take $\epsilon_j\in\{\pm\}$, and the set of signs $(\epsilon_1,\ldots,\epsilon_n)$ is
referred to as the KP sign.
Note here that the $k\times n$ coefficient matrix $A$ is then given by 
\begin{equation}\label{hepsilon}
A=\begin{pmatrix}
1 &1& \cdots & 1 \\
\lambda_1 &\lambda_2 &\cdots & \lambda_n\\
\vdots &\vdots &\ddots &\vdots\\
\lambda_1^{k-1}&\lambda_2^{k-1}&\cdots &\lambda_n^{k-1}
\end{pmatrix} ~h_{\epsilon},
\end{equation}
where $h_{\epsilon}$ is the diagonal matrix whose entries are $\pm1$, i.e. $h_{\epsilon}:={\rm diag}(\epsilon_1,\ldots,\epsilon_n)$.  This matrix specifies the blow-ups of the KP flow, and plays an 
important role for our study (subsection \ref{signS}).
Because of the ordering $\lambda_1<\cdots<\lambda_n$, the sign of each Pl\"ucker coordinate
is given by
\begin{equation}\label{signP}
{\rm sign}(\xi(\sigma_1,\ldots,\sigma_k))=\prod_{j=1}^k\epsilon_{\sigma_j}=:\epsilon(\sigma_1,\ldots,\sigma_k),
\end{equation}
which will be used as the sign for the Schubert cell $X_{w_{\sigma}}=(\sigma_1,\ldots,\sigma_k)$ (see subsection \ref{signS}).

\begin{figure}[t!]
\centering
\includegraphics[scale=0.5]{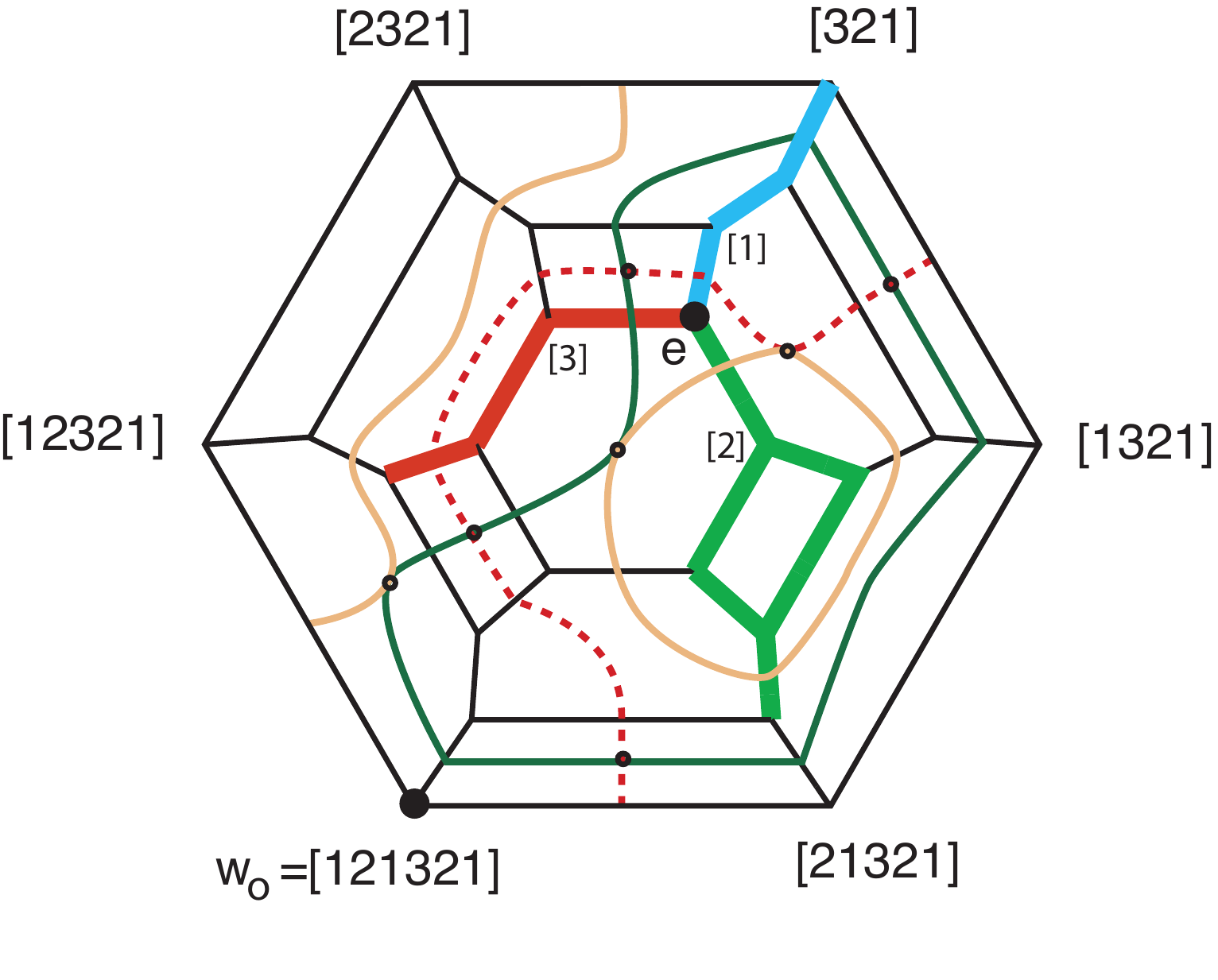}
\caption{The moment polytope of $\SL_4(\Real)$ Toda lattice hierarchy with the signs $\tilde{\epsilon}_j:={\rm sign}(a_j)=-$
for all $j=1,2,3$.
Each vertex is marked by the element of the symmetric group $\S_4$, which is denoted by
$[i_1i_2\ldots i_j]=s_{i_1}s_{i_2}\cdots s_{i_j}$
The three subgraphs indicated by the thick edges with different colors correspond to the incidence graphs of $\Gr(k,4)$.  Namely, the edges connecting with the vertices $\{e,[1],[21],[321]\}$ corresponds to $\Gr(1,4)$, the edges with $\{e,[2],[12],[32],[312],[2312]\}$ corresponds to $\Gr(2,4)$, and
those with $\{e,[3],[23],[123]\}$ for $\Gr(3,4)$.
 The zero sets of the $\tau$-functions are shown by the curves crossing the edges (see \cite{casian:06}).  For example,
the dotted curve shows $\tau_1=0$, the light color curve shows $\tau_2=0$, and 
the dark one shows $\tau_3=0$.}
\label{fig:G24}
\end{figure}

Without loss of generality, we choose the parameters $\lambda_j$ in $E_j$ 
with $\sum_{j=1}^n\lambda_j=0$ (since the matrix $L$ can be shifted by a constant in the diagonal, i.e. 
$L\to L+cI$ with the identity matrix $I$ does not change the equation).
Then the set of $\tau$-functions,
\begin{equation}\label{KPtau}
\tau_k:={\rm Wr}(f,f',\ldots, f^{(k-1)}), \qquad k=1,2,\ldots,n.
\end{equation}
gives the solution of the Toda lattice equation for $\SL_n(\Real)$, which is defined with
the $n\times n$ tri-diagonal matrix $L$,
\[
L=\begin{pmatrix}
b_1 & 1 & 0 &\cdots & 0 \\
a_1 & b_2 & 1 & \cdots & 0 \\
\vdots & \ddots &\ddots &\ddots&\vdots\\
0 & \cdots &\cdots & b_{n-1} & 1 \\
0 &\cdots &\cdots & a_{n-1}&b_n
\end{pmatrix}
\]
The Toda lattice equation is then expressed by the matrix equation (i.e. the Lax equation),
\[
\frac{dL}{dt}=[L,L_-],
\]
where $L_-$ represents the (strictly) lower triagular part of the matrix $L$.  The solution is then given by
\begin{equation}\label{Todaab}
a_j=\frac{\tau_{j+1}\tau_{j-1}}{\tau_j^2} \qquad b_j=\frac{d}{dt}\,\ln\frac{\tau_j}{\tau_{j-1}},
\end{equation}
with $\tau_0=1$.
The Toda lattice equation has several commuting flows, and the set of those flows forms the Toda
lattice hierarchy which can be defined as
\[
\frac{\partial L}{\partial t_j}=[L,L_-^j]\qquad {\rm with}\quad L^j_-:=(L^j)_-,\quad j=1,2,\ldots,n-1.
\]
Note that each $\tau_k$ then gives a solution of the KP equation on $\Gr(k,n)$ with the identifications
$t_1=x, t_2=y$ and $t_3=t$. The solution of the Toda lattice
hierarchy is given by the set $(\tau_1,\ldots,\tau_{n-1})$ (note $\tau_n$ is just a constant, since $\sum_{j=1}^n\lambda_j=0$).  Since each $\tau_k$  can be considered as a point of
$\Gr(k,n)$,  the solution of the Toda lattice equation defines a point of the flag manifold $\mathcal{B}:=\SL_n(\Real)/B$
with the Borel subgroup $B$.  This is just a consequence of the
diagonal embedding, denoted by $\iota$, of $\mathcal{B}$ into the product of the Grassmannians, i.e.
\[
\iota\quad:\quad\mathcal{B}\quad\hookrightarrow \quad \Gr(1,n)\times\Gr(2,n)\times\cdots\times\Gr(n-1,n).
\]
We then define a map $\varphi:\Real^{n-1}\to\mathfrak{h}^*_{\Real}$ as the composite map
$\varphi=\mu\circ \iota$,
\[
\varphi(\mathbf{t})=\sum_{k=1}^{n-1}\mu(\tau_k(\mathbf{t})),
\]
where $\mu:\Gr(k,n)\to\mathfrak{h}^*_{\Real}$ is the moment map defined in \eqref{moment}.
The moment polytope for the Toda lattice hierarchy is the permutohedron of the symmetric group $\mathcal{S}_n$,
whose vertices are labeled by the elements of $\S_n$.  In \cite{casian:06}, we used this fact, and
described the cohomology of the real flag variety in terms of the Toda flow and its blow-ups.  We basically follow the arguments given in that paper.
However we have found several new structures for the case of $\Gr(k,n)$, not just by a projection
$\pi:\mathcal{B}\to\Gr(k,n)$, and in fact, we found the explicit form of the Poincar\'e polynomials for $\Gr(k,n)$.
Figure \ref{fig:G24} illustrates the moment polytope (permutohedron) for $\SL_4(\Real)$ Toda lattice hierarchy. The figure shows the incidence graph of the real flag variety $\mathcal{B}=\SL_4(\Real)/B$
as shown in \cite{casian:06}.

%%%%%%%%%%%%%%%%%%%%%%%%%%%%%%%%%%%%%%%%%%%%%%%%
%%%%%%%%%%%%%%%%%%%%%%%%%%%
Each $\tau$-function gives the zero divisor on the polytope which corresponds to the singular solution of
the Toda lattice hierarchy.  Although one can recover  the incidence graph of $\Gr(k,4)$ from the Figure
as the subgraph of this flag picture, in this paper we give an alternative construction  based on the KP low which leads to explicit
cohomology calculations.  We do use the flag picture based on the Toda flow to prove our main theorem \ref{main1}.
Below we explain  the  links, similarities and differences between the Toda and the KP approaches with an example.

\subsection{Remark on the Toda lattice description} \label{remarktoda}

Here using the $\SL_4(\Real)$ Toda lattice,  we give a brief overview of the results  in  \cite{casian:06} on  the incidence graph of the real flag manifold $\mathcal{B}=\SL_4(\Real)/B$.
The cohomology of $\mathcal{B}$
then  arises from the Toda lattice with  the {\it Toda signs} $\tilde\epsilon_j={\rm sign}(a_j)$ (see Definition \ref{signos} below) attached to the vertex  $e$ given by  $\tilde\epsilon_j=-$ for all $j=1,2,3$.
Since the solution $(a_j,b_j)$ of the Toda lattice is given in the form \eqref{Todaab}, the $\tau_1$-function
has the form,
\[
\tau_1=E_1-E_2+E_3-E_4,
\]
that is, we have the KP signs $(\epsilon_1,\ldots,\epsilon_4)=(+,-,+,-)$ (the solution
corresponding to this choice of the signs has the maximum number of blow-ups from $t\to-\infty$ to $t\to\infty$, see \cite{KY:98}).  Note that from \eqref{Todaab}, the Toda signs are defined with the $\epsilon_j$'s in the $\tau_1$-function in the form,
\[
\tilde\epsilon_j:={\rm sign}(a_j)=\epsilon_j^{-1}\epsilon_{j+1}.
\]
This corresponds to the relation between the roots $\alpha$ and the fundamental weights $\omega$, i.e.
with the character $\chi$ of the group, they are defined by
$\tilde\epsilon=-\chi_{\alpha}(h)$ and $\epsilon=-\chi_{\omega}(h)$ for $h\in H$, the Cartan subgroup
\cite{CK:00} (see also subsection \ref{Todasigns}).
Then the sign change along the edge $w\to s_iw$ is given by (Proposition 5.1 in \cite{casian:06}, see also Proposition \ref{propsigns})
\[
\tilde\epsilon_j'=\tilde\epsilon_j\tilde\epsilon_i^{-C_{i,j}},
\]
where $C_{i,j}$ is the Cartan matrix for $\SL_4(\Real)$.
The co-boundaries  involving  two vertices $w$ and $s_iw$ are zero
when  one crosses a blow-up (singularity), that is, the edge crosses  $\tau_i=0$, (i.e. $\tilde\epsilon_j\tilde\epsilon'_j=-$). 

Then the   incidence graph of $\Gr(k,4)$  is  obtained as a subgraph of the incidence graph for the flag manifold. If we consider $\Gr(1,4)$ then the incidence graph is the subgraph containing the vertices
$\{ e, s_1, s_2s_1, s_3s_2s_1 \}$. In Figure \ref{fig:G24},  this is the subgraph 
connecting with the vertices $\{e,[1],[21],[321]\}$ which crosses curves of the form $\tau_i=0$ twice. 
That is,  along $s_1$ one crosses $\tau_1=0$, then along the edge from $[1]$ to $[21]$ none of the  $\tau_k=0$ are crossed, hence there is a non-zero co-boundary.  Finally  along the edge connecting
 $[21]$ and $[321]$ , $\tau_3$ becomes zero, i.e. the co-boundary is zero.  
 This gives the incidence graph giving the cohomology of
 $\Gr(1,4)$, that is,   we have
the graph, $e\to s_1\Rightarrow s_2s_1 \rightarrow s_3s_2s_1$.
In general, the incidence graph of the partial flag $G/P$ appears as a subgraph of the incidence
graph of the flag $G/B$ (see \cite{casian99} and Remark \ref{subgraph} below ).

 %labibi
 
 In this paper, we consider the cohomology of $\Gr(k,n)$ directly from the KP flow, that is, 
 we construct the incidence graph based on a single $\tau$-function, i.e. $\tau_k$.
 One should then note that for example, the incidence graph of $\Gr(1,4)$ from the KP flow is obtained by the choice of
 the $\tau$-function in \eqref{tauG14}, i.e. we have $(\epsilon_1,\ldots,\epsilon_4)=(+,-,-,+)$ as we discussed above.   Notice that 
 the corresponding Toda sign is $(\tilde\epsilon_1,\tilde\epsilon_2,\tilde\epsilon_3)=(-,+,-)$,
 which is {\it not} the same as the Toda signs used for the flag picture of Figure \ref {fig:G24}. 
 The  Toda sign $(-,+,-)$
 would give an incidence graph for cohomology  with twisted coefficients of the real flag manifold
 and this is not described by Figure \ref {fig:G24}. (See Figure 4 in \cite{CK:02} for the moment polytopes
 with different signs $(\tilde\epsilon_1,\tilde\epsilon_2,\tilde\epsilon_3)$.)

 In order to get the incidence graph of $\Gr(1,4)$ in terms of the KP flow, that is related to
 the Toda flow with sign $(-,+,-)$ and taking into account only the flow  as it crosses $\tau_1=0$, we then need to 
 consider  the quotient $\S_4/\langle s_2, s_3\rangle$ corresponding to the natural projection,
 $\pi:G/B\to G/P$ with the maximum parabolic $P$.
  Namely we are now in $\Gr(1,4)$ viewed as
 $\SL_4(\Real) / P$ and in  this quotient $\tau_1=0$ is crossed by that edge of the tetrahedron,
 and so  we recover the incidence graph of $\Gr(1,4)$.  This is
  illustrated in Figure \ref{fig:G14}.  Also note that the projection $\pi:G/B\to \Gr(2,4)$ gives the
  octahedron illustrated in Figures \ref{fig:Grassmann24} and \ref{fig:Ttype}.
 
 It turns out that  there are
 exactly two distinguished Toda signs (and two corresponding distinguished KP signs) which need to be used for  the cohomology of all the $\Gr(k,n)$
  to be expressed in terms of the KP flow.  For example,
 in the case $\Gr(k,4)$,  these are $(-,+,-)$ for $k=1,3$ and $(+,-,+)$ for $k=2$.  Note once more,  that the
Toda sign $(\tilde\epsilon_1,\tilde\epsilon_2,\tilde\epsilon_3)$ giving the incidence graph 
in \cite{casian:06} (cohomology with constant coefficients),  or in Figure \ref{fig:G24},   is not one of these two,  but rather it is associated to  $(-,-,-)$.
 
 A summary of the strategy in the  proof  of our main theorem  given in section \ref{teorema} can be expressed as follows.  The point is to check
 that the incidence graph of $\Gr(k,n)$, viewed as a subgraph of the incidence graph
 of the real flag manifold, derived from the Toda lattice and corresponding to the
 Toda signs $\tilde\epsilon_i=-$ for all $i$,  agrees with a graph defined in terms of the KP flow and one of the two distinguished
 KP signs.

%%%%%%%%%%%%%%%%%%%%%%%%%%%%%%%%%%%%%%%%%

\section{The set $\mathcal{S}^{(k)}_n$ and the  graph $\mathcal{G}_n^{(k)}$} \label{action}
In this section, we describe the algebraic structure of the Grassmannian $\Gr(k,n)$.
The set of Schubert cells $\{(\sigma_1,\ldots,\sigma_k):
1\le\sigma_1<\cdots <\sigma_k\le n\}$ for $\Gr(k,n)$ forms a partially ordered set ({\it poset}) with the 
Bruhat   order defined as follows:
Let $s_i\in \mathcal{S}_n$ be an adjacent transposition $s_i: i \to i+1$ for $i=1,\ldots,n-1$ for
the symmetric group $\mathcal{S}_n$. Then there exists $s_{\sigma_i}$-action on the cell $(\sigma_1,\ldots,\sigma_k)$,
if $\sigma_i+1<\sigma_{i+1}$ for some $i$ or $\sigma_k<n$. The action gives
\[
s_{\sigma_i}:  (\sigma_1,\ldots,\sigma_i,\ldots,\sigma_k)\longrightarrow (\sigma_1,\ldots,\sigma_{i-1},\sigma_i+1,\sigma_{i+1},\ldots,\sigma_k)\,.
\]
Then the arrow given by the action determines the {\it weak}  Bruhat  order between those cells.  With  action defined this way, each cell can be uniquely parametrized by a representative  of  minimal length in 
 the quotient $\S_n/\mathcal{P}_k$ (see \cite{bjorner:05}), that is, 
 \[ 
 \mathcal{S}^{(k)}_n:=\{{\rm the~reduced~words~of~mod}(\mathcal{P}_k)
 {\rm ~ending~in}~s_k\}.
 \]
 Here $\mathcal{P}_k$ is a Weyl
subgroup of $\mathcal{S}_n$ corresponding to the maximal parabolic subgroup,
and is generated by $\{s_1,\ldots, s_{k-1},s_{k+1},\ldots,s_{n-1}\}$, denoted as
\[
\mathcal{P}_k=\langle s_1,\ldots,\check{s_k},\ldots,s_{n-1}\rangle\cong \mathcal{S}_k\times\mathcal{S}_{n-k}.
\]
Note here that $\mathcal{P}_n=\S_n$ and $\S_n^{(n)}=\{e\}$.
In this sense, we identify an element
$w_{\bf{\sigma}}\in\mathcal{S}^{(k)}_n$ with the Schubert cell ${\sigma}=(\sigma_1,\ldots,\sigma_k)$.
This leads to the Schubert decomposition in terms of $\mathcal{S}_n^{(k)}$, that is, we have
\[
\Gr(k,n)=\bigsqcup_{w_{\sigma}\in\mathcal{S}_n^{(k)}} X_{w_{\sigma}},\quad{\rm with}\quad  X_{w_{\sigma}}=(\sigma_1,\ldots,\sigma_k).
\]
We also write $X_{w_{\sigma}}$ simply as $w_{\sigma}$, that is, $w_{\sigma}=(\sigma_1,\ldots,\sigma_k)$.
With the set $\mathcal{S}^{(k)}_n$, one should also note
\[
|\Gr(k,n,\mathbb{F}_q)|=\sum_{w\in\mathcal{S}^{(k)}_n}q^{l(w)}\,,
\]
where $l(w)$ is the length of $w$.

The set of Schubert cells with $s_i$-actions forms a graph consisting of the vertices given by the cells
and the edges given by the weak Bruhat order. We denote this graph by $\mathcal{G}^{(k)}_n$,
and call it the weak Bruhat graph of $\Gr(k,n)$.  Let us give few examples of $\mathcal{G}_n^{(k)}$.

\begin{Example} \label{ex1} We consider the cases with $n=4$ and $k=1,2,3$:

\begin{itemize}
\item[(1)] $\Gr(1,4)\cong \mathbb{R}P^3:$ The graph $\mathcal{G}_4^{(1)}$ is given by
\[\begin{matrix}
(1) & \mapright{s_1} & (2) & \mapright{s_2} & (3)&  \mapright{s_3} &(4)
\end{matrix}\,.
\]
Each cell $(k)$ can be parametrized by an element in $\S_4^{(1)}$,  the set of minimal length representative 
in the quotient $\mathcal{S}_4/\langle s_2,s_3\rangle$, that is, $(1)=e, (2)=s_1, (3)=s_2s_1, (4)=s_3s_2s_1$. \item[(2)]  $\Gr(2,4):$ The graph $\mathcal{G}_4^{(2)}$ is given by
\[
\begin{matrix}
 (1,2)   &\mapright{s_2}&      (1,3)             & \mapright{s_3}    &      (1,4)    \\
             &                         &\mapdown{s_1} &                          &        \mapdown{s_1}   \\
            &                          &      (2,3)          &  \mapright{s_3}    &     (2,4)    \\
            &                           &                        &                                &   \mapdown{s_2} \\
            &                            &                        &                                &       (3,4)
  \end{matrix}.
\]
Each cell $(i,j)$ is parametrized by a unique element of $\mathcal{S}^{(2)}_4$, i.e.
$(1,2)=e, (1,3)=s_2, (2,3)=s_1s_2, (1,4)=s_3s_2, (2,4)=s_1s_3s_2, (3,4)=s_2s_1s_3s_2$.
\item[(3)]  $\Gr(3,4)\cong\mathbb{R}P^3$:  The graph $\mathcal{G}_4^{(3)}$ is given by
\[
\begin{matrix}
(1,2,3) & \mapright{s_3} & (1,2,4)& \mapright{s_2} & (1,3,4)& \mapright{s_1} &(2,3,4)
\end{matrix}\,.
\]
The parametrization of each cell $(i,j,k)$ is given by $(1,2,3)=e,(1,2,4)=s_3, (1,3,4)=s_2s_3,
(2,3,4)=s_1s_2s_3$. Those are the elements of $\mathcal{S}^{(3)}_4$.
\end{itemize}

\end{Example}

\subsection{The decomposition of $\mathcal{S}^{(k)}_n$}

Here we consider a decomposition of
the set $\mathcal{S}^{(k)}_n$ into the set  $\{\mathcal{S}^{(k-1)}_{j-1}: j=k,k+1,\ldots,n\}$.
This decomposition is based on the following relation for the binomial coefficients,
\[
\binom{n}{k}=\binom{k-1}{k-1}+\binom{k}{k-1}+\cdots+\binom{n-1}{k-1}=\sum_{j=k}^n\binom{j-1}{k-1}\,.
\]
which is a direct consequence of the Pascal rule, $
\binom{n}{k}=\binom{n-1}{k}+\binom{n-1}{k-1}$.
Namely we have:
\begin{Proposition}\label{decompositionS}
The set $\mathcal{S}^{(k)}_n$ has a decomposition,
\[
\mathcal{S}^{(k)}_n=\mathcal{S}^{(k-1)}_{k-1}\cup\bigcup_{j=k}^{n-1}\mathcal{S}^{(k-1)}_{j}s_{j}\cdots s_k\,,
\]
where note $\mathcal{S}^{(k-1)}_{k-1}=\{e\}$.
\end{Proposition}

To prove
Proposition \ref{decompositionS},  we first state the following Lemma which is similar to the
Pascal rule,
\begin{Lemma}\label{pascal}
There is a decomposition,
\[
\mathcal{S}^{(k)}_n=\mathcal{S}^{(k)}_{n-1}\cup \mathcal{S}^{(k-1)}_{n-1}s_{n-1}s_{n-2}\cdots s_k\,.
\]
\end{Lemma}
\begin{Proof}
This corresponds to the decomposition of weight vectors,
\[
\bigwedge^k V(E_1,\ldots,E_n)=\bigwedge^kV(E_1,\ldots,E_{n-1})\oplus \left[\bigwedge^{k-1}V(E_1,\ldots,E_{n-1})\wedge E_n\right].
\]
It is then obvious that the weight vectors in the fast part are given by the orbit of $\mathcal{S}^{(k)}_{n-1}$
of the highest weight vector $E_1\wedge\cdots\wedge E_k$.

Now note that the element $w=s_{n-1}\cdots s_k$ maps the highest weight vector $E_1\wedge\cdots\wedge E_k$ to $E_1\wedge\cdots\wedge E_{k-1}\wedge E_n$. Then the weight vectors in $\wedge^{k-1}V(E_1,\ldots,E_{n-1})$ can be obtained by the orbit of $\mathcal{S}^{(k-1)}_{n-1}$ of the highest weight vector $E_1\wedge\cdots\wedge E_{k-1}$.
\end{Proof}

\begin{Remark} \label{decompositionx}  The decomposition in Proposition \ref{decompositionS} corresponds to  the following decomposition,
\[
\begin{array}{llll}
\displaystyle{\bigwedge^{k} V(E_1,\ldots ,E_n)}=\bigoplus_{j=k}^n\left[\bigwedge^{k-1}V(E_1,\ldots,E_{j-1})\wedge E_j\right].
\end{array}
\]
\end{Remark}

\begin{Example} \label{aa}
The set $\mathcal{S}_4^{(2)}$ can be written as  $\mathcal{S}_1^{(1)}\cup  \mathcal{S}_2^{(1)} s_2 \cup  \mathcal{S}_3^{(1)} s_3s_2$:
First we write the sets  $\mathcal{S}_m^{(k)}$ which are involved,

\begin{equation}
{\mathcal{S}_1^{(1)}}=  {\{ e \}  }, \qquad
{ \mathcal{S}_2^{(1)} }=  {\{ e, \,s_1 \} },\qquad
{\mathcal{S}_3^{(1)}}={\{ e,\, s_1, \,s_2s_1 \}}.
\end{equation}
Hence we have
\begin{equation}
{\mathcal{S}_1^{(1)}=}{\{e \} },\qquad
{  \mathcal{S}_2^{(1)} s_2=}{\{ s_2 },\,{s_1s_2\}} ,\qquad
{   \mathcal{S}_3^{(1)} s_3s_2 =} {\{s_3s_2 }, \,{s_1s_3s_2} ,\,{ s_2s_1s_3s_2 \}}.
\end{equation}
As in Remark \ref{decompositionx} we have a decomposition   
$$\bigwedge^2 V(E_1, E_2, E_3)=\left[\bigwedge^1 V(E_1)\wedge E_2\right] \oplus \left[\bigwedge^1 V(E_1,E_2)\wedge E_3\right] \oplus  \left[\bigwedge^1 V(E_1,E_2, E_3)\wedge E_4\right].$$ 
For example,
$ \bigwedge V(E_1,E_2, E_3)\wedge E_4$  corresponds to   ${ \mathcal{S}_3^{(1)}} s_3s_2$.  This is given explicitly by associating  $E_1\wedge E_4$, $E_2\wedge E_4$, $E_3\wedge E_4$  to the elements $s_3s_2, s_1s_3s_2, s_2s_1s_3s_2$ respectively (see also Example \ref{gr24} below).
\end{Example}

The decomposition of $\mathcal{S}_n^{(k)}$ provides an arrangement of the Schubert decomposition of $\Gr(k,n)$:
\begin{equation}\label{arrangement}
\Gr(k,n)=\bigsqcup_{j=k}^{n}\left(\bigsqcup_{1\le\sigma_1<\cdots<\sigma_{k-1}\le j-1}X(\sigma_1,\ldots,\sigma_{k-1},j)\right)\,.
\end{equation}
Each Schubert cell $X(\sigma_1,\ldots,\sigma_{k-1},j)$ can be labeled by an element
of $\mathcal{S}^{(k-1)}_{j}s_{j}\cdots s_k$.  For each set $\mathcal{S}^{(k-1)}_{j}$,
we can define the Bruhat graph $\mathcal{G}^{(k-1)}_{j}$.  We then have a decomposition of the graph $\mathcal{G}^{(k)}_n$,
\[
\mathcal{G}^{(k)}_n~=~\left[~\mathcal{G}^{(k-1)}_{k-1}~\mapright{s_k}~ \mathcal{G}^{(k-1)}_{k}~\mapright{s_{k+1}}~\cdots
\cdots ~\mapright{s_{n-1}}~\mathcal{G}^{(k-1)}_{n-1}~\right]\,.
\]
Here the $s_j$-action provides the edges connecting the vertices in $\mathcal{G}^{(k-1)}_{j-1}$ with the corresponding vertices in $\mathcal{G}^{(k-1)}_j$, i.e.
\[
\begin{matrix}
\G^{(k-1)}_{j-1}\ni (\sigma_1,\ldots,\sigma_{k-1},\fbox{j}) &\mapright{s_j}&(\sigma_1,\ldots,\sigma_{k-1},\fbox{j+1})\in
\G^{(k-1)}_{j},
\end{matrix}
\]
for all $1\le \sigma_1<\cdots<\sigma_{k-1}\le j-1$, and each index in the box is fixed.
\begin{Example} 
Consider the case of $\Gr(3,5)$. The graph $\mathcal{G}_5^{(3)}$ is given by
\[
\begin{matrix}
(1,2,\fbox{3}) &\overset{s_3}\longrightarrow & (1,2,\fbox{4}) & \mapright{s_4} & (1,2,\fbox{5}) &    & \\
            &                         & \downarrow&                         &\downarrow&    & \\
            &                         &  (1,3,\fbox{4}) & \overset{s_4}\longrightarrow & (1,3,\fbox{5}) &  \longrightarrow &(1,4,\fbox{5})\\
            &                        & \downarrow &                       &\downarrow &        &   \downarrow \\
            &                          & (2,3,\fbox{4}) &  \overset{s_4}\longrightarrow &  (2,3,\fbox{5}) &  \longrightarrow  &   (2,4,\fbox{5})   \\
            &                         &               &                            &      &        &\downarrow  \\
            &                         &               &                             &      &       & (3,4,\fbox{5})
\end{matrix}
\]
which has the structure of a sequence of smaller graphs,  i.e. the decomposition of the graph $\mathcal{G}_5^{(3)}$,          
\[
\mathcal{G}^{(3)}_5~= ~ \left[~\mathcal{G}_2^{(2)}~\overset{s_3}{\longrightarrow}~\mathcal{G}_3^{(2)}~
\overset{s_4}\longrightarrow~\mathcal{G}_4^{(2)}~ \right].
\]
\end{Example}

%%%%%%%%%%%%%%%%%%%%%%%%
\subsection{The decomposition of $\Gr(k,n)$} \label{grknbundle}
We now describe a decomposition of $\Gr(k,n)$ based
on the arrangement (\ref{arrangement}).
Let us first define $Y(j)$ for $j=k,k+1,\ldots,n$ as the unions of the Schubert cells in (\ref{arrangement}),
\[
Y(j):=\bigsqcup_{1\le\sigma_1<\cdots<\sigma_{k-1}\le j-1}X(\sigma_1,\ldots,\sigma_{k-1},j)\,.
\]
Then we explain that each $Y(j)$ is a manifold with the structure of a vector bundle. To do this,
let us fix vector spaces $\{V_0(j):j=k-1,k,\ldots,n\}$ satisfying
 \[
  V_0(k-1)\subset V_0(k) \subset \cdots \subset  V_0(n-1) \subset  V_0(n)\,=\, {\mathbb R}^n\,,
 \]
 with dim$(V_0(j))=j$.  Then $Y(j)$ can be described by
 \[
 Y(j)=\{ [V] : V\subset  V_0(j),  ~{\rm dim}( V\cap V_0(j-1))=k-1\} \qquad{\rm for}\quad j=k,k+1,\ldots,n.
 \]
Now one can see that  $Y(j)$ is a manifold  which has the structure of  a vector bundle ${\mathcal V}(j)$ with fibers of dimension $j-k$ over the Grassmanian  $\Gr(k-1, j-1)$.  We have well defined projections 
$\pi_j: Y(j) \to \Gr(k-1, j-1)$ which  are  given by $\pi_j([V])= [ V\cap V_0(j-1)]$.   This is because
$V\cap  V_0(j-1) \subset  V_0(j-1) $ is a vector space of dimension $k-1$ inside a fixed
vector space  which looks like $\Real^{j-1}$. The fibers of this projection can be computed using
explicit coordinates (see below).  These vector bundles determine the determinant line-bundles  over $\Gr(k-1,j-1)$,
\begin{equation}\label{thebundles}  
\bigwedge^{j-k} {\mathcal V}(j)={\mathcal E}(j).
\end{equation}

\begin{Remark}  \label{remark0} We have the following:
\begin{enumerate}

 \item [(i)]  The space  $Y(j)$ is a union  of cells corresponding to  $\mathcal{S}^{(k-1)}_{j-1}s_{j-1} \cdots s_{k}  \leftrightarrow   {\mathcal G}_{j-1}^{(k-1)}$
 
  \item [(ii)] Each element in $\mathcal{S}_{j-1}^{(k-1)}$ labels  a corresponding cell in $\Gr(k-1,j-1)$. 
 
  \item [(iii)]  The cell of lowest dimension in  $Y(j)$ corresponds to  $s_{j-1} \cdots s_{k}$. 
  
  \begin{enumerate}
  
  \item [$\bullet$] The union of all these cells $\{ s_{j-1} \cdots s_{k}:j=k+1,\cdots n-1  \} $ correspond to   the incidence graph of  $\Real P^{n-k}\cong\Gr(1,n-k+1)$.  
  
   \item [$\bullet$]  These  cells of  $\Real P^{n-k}$   are the ones associated to   the Schubert symbols
\[
(1,2,\ldots,{k-1},j) ,\quad  {\rm for}\quad j=k,\ldots,n\,.
\] 

 \item [$\bullet$]  These are also the elements $s_j\cdots s_k$  in the statement of Proposition \ref{decompositionS} and  in Example \ref{aa}  these are the cells labeled as $e,s_2, s_3s_2$. 

\end{enumerate}
  
    \item [(iv)] The fibers of the vector bundles ${\mathcal V}(j)$    with total space $Y(j)$ and base $\Gr(k-1, j-1)$ correspond to the cells of  this  $\Real P^{n-k}$.
 \end{enumerate}
\end{Remark}

In the following examples, we describe explicitly  the cells involved in each 
piece of this decomposition and the correspondence between  a cell decomposition of
$Y(j)$ and  $\bigwedge^{k-1} V(E_{1},\dots, E_{j-1})\wedge  E_{j}.$  The basis element
$E_{\sigma_1}\wedge \cdots\wedge E_{\sigma_{k-1}}\wedge E_j$ label  a cell with pivots 
at the positions $(\sigma_1,\ldots,\sigma_{k-1},j)$.

\begin{Example}  \label{gr23} Explicit decomposition  in the case of  $\Gr(2,3)\cong \Real P^2$. 
Let $V_0(2)\subset {\mathbb R}^3$ consist of vectors of the form $(x_1,x_2,0)$.
 We let $V$ denote a two  dimensional
vector space in  ${\mathbb R}^3$ and $[V]$ a corresponding point in $\Gr(2,3)$. We than have

\begin{align*}
Y(2)=&\{  [V_0(2)] \}=\left\{\begin{pmatrix}1 &0&0\\0&1&0\end{pmatrix}\right\}, \\[2.0ex]
Y(3)= & \{ [V]: V \subset   {\mathbb R}^3,  ~{\rm dim} (V\cap V_0(2))=1 \}\\[1.5ex]
    =&\left\{\begin{pmatrix}1&0&0\\0&*&1\end{pmatrix}\right\}\sqcup\left\{\begin{pmatrix}*&1&0\\ *&0&1\end{pmatrix}\right\}.\\
\end{align*}
Looking at the position of pivots  $Y(2)$ should be thought as corresponding to $\bigwedge^1 V(E_1)\wedge E_2$ and $Y(3)$ to the summand $\bigwedge^1 V(E_1,E_2)\wedge E_3$ in a decomposition of  
$\bigwedge^2 V(E_1,E_2, E_3)$.

    If we think of $\Gr(2,3) \cong {\mathbb R} P^2$  given as a sphere where antipodal points
are identified, then  this decomposition has a simple description. First  $Y(2)$ can be identified with the north and south poles. If we delete this pair of antipodal points we are left with a M\"obius band corresponding to $Y(3)$. This is a non-trivial  line bundle over the equator and in matrix notation,  $\pi_3$ projects a matrix   to its first row  so that
$\pi_3(Y(3))=\{(1~0~0)\}\sqcup\{(*~1~0)\}\cong \Real P^1$). The 
one dimensional fibers
are given by the second row in a matrix i.e.  $(0~*~1)$ and $(*~0~1)$ over $\{(1~0~0)\}$ and $\{(*~1~0)\}$, respectively.

The $\Real P^1$ in item (iii) of Remark  \ref{remark0} is a meridian. Note that the direction along meridians corresponds to the fibers of the fiber bundle over the equator giving the  M\"obius band.
 
\end{Example}

\begin{Example} \label{gr24}  Explicit decomposition  in the case of  $\Gr(2,4)$. 
We  once more emphasize the connection with a decomposition:

\begin{equation} \label{dd} 
\bigwedge^2 V(E_1, E_2, E_3, E_4)=\left[\bigwedge^1 V(E_1)\wedge E_2\right] \oplus \left[\bigwedge ^1V(E_1, E_2)\wedge E_3\right] \oplus  \left[\bigwedge^1 V(E_1, E_2, E_3)\wedge E_4\right]. \end{equation}

%We consider the matrix $h(+)$ acting on $\bigwedge^2 V(E_1, E_2, E_3, E_4)$. 
%Hence we have $(E_1, 1), (E_2, 1), (E_3,q), (E_4,q)$.

 Let $V_0(2)\subset V_0(3) \subset {\mathbb R}^4$ be two subspaces  of dimensions $2$ and $3$ respectively.  We let $V$ denote a two  dimensional
vector space in  ${\mathbb R}^4$ and $[V]$ a corresponding point in $\Gr(2,4)$. Set:
\begin{align*}
Y(2)=&\{  [V_0(2)] \}=\left\{\begin{pmatrix}1&0&0&0\\0&1&0&0\end{pmatrix}\right\},\\[2.0ex]
Y(3)=&\{ [V]: V \subset   V_0(3) , \,{\rm dim }(V\cap V_0(2))=1 \}\\[1.5ex]
       =&\left\{\begin{pmatrix}1&0&0&0\\0&*&1&0\end{pmatrix}\right\}\sqcup
       \left\{\begin{pmatrix} *&1&0&0\\ *&0&1&0\end{pmatrix}\right\},\\[2.0ex]
Y(4)=&\{ [V]: V\subset  {\mathbb R}^4,\,{\rm  dim}( V\cap V_0(3))=1 \}\\[1.5ex]
       =&\left\{\begin{pmatrix}1&0&0&0\\0&*&*&1\end{pmatrix}\right\}\sqcup\left\{
       \begin{pmatrix} *&1&0&0\\ *&0& * &1\end{pmatrix}\right\}\sqcup\left\{
       \begin{pmatrix} *& * & 1 & 0\\ * & * &0&1\end{pmatrix}\right\}.
\end{align*}
We have $Y(2)=\Gr(2,2)\cong \Gr(1,1)$, $Y(3)$ is a non-trivial line bundle over  $\Gr(1,2)$, and
$Y(4)$ is a non-trivial  ${\Real}^2$ bundle over $\Gr(1,3)$. The projections of
these bundles are $\pi_3 : Y(3)\to \Gr(1, 2)$ with $\pi_3([V])= [ V\cap V_0(2)]$, and 
$\pi_4 : Y(4)\to \Gr(1, 3)$ with $\pi_4([V])= [ V\cap V_0(3)]$. More explicitly, we have
\begin{align*}
\pi_3(Y(3))=&\{(1~0~0~0)\}\sqcup\{(*~1~0~0)\}=\Gr(1,2),\\[1.0ex]
\pi_4(Y(4))=&\{(1~0~0~0)\}\sqcup\{(*~1~0~0)\}\sqcup\{(*~*~1~0)\}=\Gr(1,3).
\end{align*}
The  fiber over   $\{(1~0~0~0)\}$ in $Y(3)$
is $(0~*~1~ 0)$, and
the fiber over   $\{(*~1~0~0)\}$ is $(*~0~1~0)$. This becomes a non-trivial line bundle over  the circle $\Gr(1,2)$.   The case of $Y(4)$ is described similarly and results in a  non-trivial vector bundle with fiber
of dimension $2$. 
The fiber over $\{(1~0~0~0)\}$ is $\{(0~*~*~1)\}$, the fiber over $\{(*~1~0~0)\}$ is $\{(*~0~*~1)\}$,
and the fiber over $\{(*~*~1~0)\}$ is $\{(*~*~0~1)\}$.
The fiber of  the vector bundle over $\Gr(1,3)$  is two dimensional because  of the two
parameters.
\end{Example}

%%%%%

%%%%%%%%%%%%%%%%%%%%%%%%%%%%%%%%%%%%%%%%%%%%%
\section{The graph $\G(k,n)$} \label{grafica} 

Here we define what will be shown to be  the  incidence graph for the cohomology of $\Gr(k,n)$ relative to Schubert cells. 
We first consider the case of
$\Gr(1,n)\cong \Real P^{n-1}$, and then introduce {\it signed}
Schubert cells to construct what will become the incidence graph.
 The goal of this section is to show that this  graph defined for
$\Gr(k,n)$ can be obtained from the results of $\Gr(1,n)$ and its dual $\Gr(n-1,n)$.

Let us recall that $\Gr(1,n)\cong\Real P^{n-1}$ has the Schubert decomposition,
\[
\Gr(1,n)=\bigsqcup_{1\le\sigma_1\le n}X(\sigma_1)\,,
\]
where $X(\sigma_1)$ is the Schubert cell given by $X(i)=\{(x_1,\ldots,x_{i-1},1,0,\ldots,0):x_j\in\Real\}\cong\Real^{i-1}$.
Then the cochain complex for $\Gr(1,n)$ is expressed by
\begin{equation}\label{incidenceRPn}
\mathbb{Z}[\langle(1)\rangle]~\overset{\delta_0}{\longrightarrow}~\mathbb{Z}[\langle(2)\rangle]~\overset{\delta_1}{\Longrightarrow}~\cdots\cdots~\overset{\delta_{2m-2}}{\longrightarrow}~
\mathbb{Z}[\langle(2m)\rangle]~\overset{\delta_{2m-1}}{\Longrightarrow}~\mathbb{Z}[\langle(2m+1)\rangle]~\overset{\delta_{2m}}{\longrightarrow}~\cdots\,.
\end{equation}
with $\langle (i)\rangle=\overline{X(i)}$ and the coboundary operator $\delta_i:\mathbb{Z}[\langle(i+1)\rangle]\to\mathbb{Z}[\langle(i+2)\rangle]$.  
Here the single arrow $\to$ indicates $0$ incidence number and the double arrow indicates
nonzero incidence number which is either $2$ or $-2$ \cite{koch}, i.e. 
\[
\delta_{i-1}\langle (i)\rangle=\left\{\begin{array}{lllll}
0 \quad & {\rm if}\quad& i={\rm odd} \\
\pm 2\,\langle(i+1)\rangle\quad &{\rm if}\quad &i={\rm even}
\end{array}\right.
\]
Then we recover the well-known formulas of the integral cohomology of $\Gr(1,n)\cong \Real P^{n-1}$ as follows:
\begin{itemize}
\item[(a)] For $n=2m+1$, the last arrow 
in (\ref{incidenceRPn}) is the double arrow, and we have
\[
H^k(\Real  P^{2m};\mathbb{Z})=\left\{\begin{array}{llll}
\mathbb {Z}\quad &{\rm if}\quad &k=0 \\
  0  \quad   &{\rm if} \quad & k={\rm odd}\le 2m-1 \\
 \mathbb{Z}_2\quad &{\rm if}\quad &k={\rm even}\le 2m
\end{array}\right.
\]
The Poincar\'e polynomial for $\Gr(1,2m+1)$ is then just $P_{(1,2m+1)}(t)=1$.
\item[(b)] For $n=2m+2$, the last arrow is the single one, and
\[
H^k(\Real  P^{2m+1};\mathbb{Z})=\left\{\begin{array}{llll}
\mathbb {Z}\quad &{\rm if}\quad &k=0 ~{\rm and~if}~ k=2m+1 \\
  0  \quad   &{\rm if} \quad & k={\rm odd}\ne 2m+1 \\
 \mathbb{Z}_2\quad &{\rm if}\quad &k={\rm even}\le 2m
\end{array}\right.
\]
The Poincar\'e polynomial for this case is $P_{(1,2m+2)}(t)=1+t^{2m+1}$.
\end{itemize}

%%%%%%%%%%%%%%%%%%%%%%%%%%%%%%%%%%%%%

\subsection{The signed Schubert cells}\label{signS}
In order to represent the incidence graph in terms of the Schubert cells, we introduce
the signed Schubert cells (see also subsection \ref{dynamical}). 
 Let us first define the signed vectors $e_j:=\epsilon_jE_j$ for $j=1,\ldots,n$ where
$E_j\in\mathbb{R}^n$ is the $j$-th standard basis vector, and $\epsilon_j\in\{\pm \}$.
Then the $k$-wedge product $e_{\sigma_1}\wedge\cdots\wedge e_{\sigma_k}$ has the sign
$\epsilon(\sigma_1,\ldots,\sigma_k)=\prod_{j=1}^k \epsilon_{\sigma_j}$ (see \eqref{signP}),  i.e.
\[
e_{\sigma_1}\wedge\ldots\wedge e_{\sigma_k}=\epsilon(\sigma_1,\ldots,\sigma_k)\,E_{\sigma_1}\wedge\ldots\wedge E_{\sigma_k}\,.
\]
We then define a (induced) sign of the Schubert cell $(\sigma_1,\ldots,\sigma_k)$ as the sign of the wedge product $e_{\sigma_1}\wedge\cdots\wedge e_{\sigma_k}$, i.e. 
$\epsilon(\sigma_1,\ldots,\sigma_k)$.  Thus we identify the signed cell $(\sigma_1,\ldots,\sigma_k)$
as $e_{\sigma_1}\wedge\cdots\wedge e_{\sigma_k}$. 

\begin{Remark}
With the signed bases $\{e_j:j=1,\ldots,n\}$, the $\tau$-function in \eqref{KPtau} can be written in the following form
with $f=\sum_{j=1}^k\epsilon_jE_j$,
\begin{align*}
\tau_k&={\rm Wr}(f,f',\ldots,f^{(k-1)})\\
&=\sum_{1\le \sigma_1<\cdots<\sigma_k\le n}|\xi(\sigma_1,\ldots,\sigma_k)|\epsilon(\sigma_1,\ldots,\sigma_k)E(\sigma_1,\ldots,\sigma_k)
\end{align*}
where $|\xi(\sigma_1,\ldots,\sigma_k)|$ is the $k\times k$ minor of the Vandermonde matrix in $A$.
\end{Remark}

We now define a graph ${\mathcal G}(k,n)$, which depends on the  choice of signed vectors $e_i$  and which is obtained from the weak  Bruhat graph whose edges indicate that
the cells connected by the edge have the same sign.  Suppose that we have $s_{\sigma_i}$-action on the cell $(\sigma_1,\ldots,\sigma_k)$,
i.e. $\sigma_i+1<\sigma_{i+1}$. Then if the cells connected by the action have the same sign, i.e. $\epsilon(\sigma_1,\ldots,\sigma_i,\cdots,\sigma_k)\,\epsilon(\sigma_1,\ldots,\sigma_i+1,\ldots,\sigma_k)=+$, we put the double edge for the Bruhat order, i.e.
\[
s_{\sigma_i}:(\sigma_1,\ldots,\sigma_i,\ldots,\sigma_k)\Longrightarrow (\sigma_1,\ldots,\sigma_i+1,\ldots,\sigma_k)\,.
\]
Otherwise (i.e. sign change), we keep the single edge of the Bruhat order. Thus this new graph is obtained
from the Bruhat graph by changing some of the edges to the double ones. The single edges correspond
to the cells which have the different signs, i.e. $\epsilon(\sigma)\epsilon(\sigma')=-$ for
the Bruhat ordered cells with $\sigma=(\sigma_1,\ldots,\sigma_k)$ and $\sigma'=(\sigma'_1,\ldots,\sigma'_k)$.

We also impose the signs for particular cells $(\sigma_1,\ldots,\sigma_k)$ as follows:
\begin{itemize}
\item[(a)] Assume $\epsilon_1=+$.
\item[(b)] Assign $\epsilon(1,\ldots,k-1,k+j)=(-)^{\lfloor\frac{j+1}{2}\rfloor}$ for $j=0,1,\ldots,n-k$.
\item[(c)] Assign $\epsilon(1,\ldots,k-j,k-j+2,\ldots,k+1)=(-)^{\lfloor\frac{j+1}{2}\rfloor}$ for $j=1,2,\ldots,k-1$.
\end{itemize}
The item (b) implies that we have the same sign pattern for $\mathbb{R}P^{n-k}=\Gr(1,n-k+1)$.  The item (c) then corresponds to the pattern of $\mathbb{R}P^{k}=\Gr(k,k+1)$. Those cells appear in the upper horizontal line and the
left vertical side of the Bruhat graph.
Then one can determine uniquely the signs $(\epsilon_1,\ldots,\epsilon_n)$ for a given pair of the numbers
$(k,n)$.
\begin{Proposition}
Following the sign assignment given above for $\Gr(k,n)$, the signs $\epsilon_j$ are given by
\[\left\{\begin{array}{lllll}
\epsilon_{2j+1}=(-)^j \quad &{\rm for}\quad j=0,1,\ldots, \lfloor\frac{n-1}{2}\rfloor, \\[1.0ex]
\epsilon_{2j+2}=(-)^{k+j}\quad  &{\rm for}\quad j=0,1,\ldots,\lfloor\frac{n-2}{2}\rfloor.
\end{array}\right.
\]
\end{Proposition}
\begin{Proof}
First we note that the sign choice in (b) leads to the condition,
\[
\epsilon_{k+j-1}\epsilon_{k+j}=(-)^j \qquad{\rm for}\quad j=1,2,\ldots, n-k.
\] 
Also from the choice (c), we have
\[
\epsilon_{k-j+1}\epsilon_{k-j+2}=(-)^j \qquad{\rm for}\quad j=1,2,\ldots,k.
\]
Combining those, we have
\[
\epsilon_j\epsilon_{j+1}=(-)^{k+j+1}\qquad {\rm for}\quad j=1,2,\ldots,n-1,
\]
With the condition (a), i.e. $\epsilon_1=+$, we obtain the desired formulae.
\end{Proof}
This proposition implies that for $\Gr(k,n)$, if $k$ is {\it odd}, we choose the signs,
\begin{equation}\label{signOdd}
(\epsilon_1,\epsilon_2,\epsilon_3,\epsilon_4,\ldots,\epsilon_n)=(+,-,-,+,\cdots, (-)^{\lfloor\frac{n}{2}\rfloor} ),
\end{equation}
and if $k$ is {\it even}, we choose
\begin{equation}\label{signEven}
(\epsilon_1,\epsilon_2,\epsilon_3,\epsilon_4,\ldots,\epsilon_n)=(+,+,-,-,\ldots,(-)^{\lfloor\frac{n-1}{2}\rfloor}).
\end{equation}
%alternative00
It will be useful to describe those sets of signs using the diagonal matrix $h_{\epsilon}$ defined in \eqref{hepsilon}.  Namely we define
 $h(-)$ and $h(+)$ as the $n\times n$ diagonal matrices $h_{\epsilon}$ corresponding to the sets of signs
 \eqref{signOdd} and \eqref{signEven}, respectively, i.e.
 \begin{align*}
 h(-):=&{\rm diag}\left(+,-,-,+,\ldots, (-)^{\lfloor \frac{n}{2}\rfloor}\right),\\
 h(+):=&{\rm diag}\left(+,+,-,-,\ldots,(-)^{\lfloor\frac{n-1}{2}\rfloor}\right).
 \end{align*}
 These act on $V(E_1,\cdots, E_n)$.   These
 matrices then belong to the group ${\rm SL}(n,\Real)^\pm$ and also act on $\bigwedge^k V(E_1,\cdots ,E_n)$.
 This gives each $E_{j_1}\wedge \cdots\wedge E_{j_k}$ (i.e. {\em each cell}) a sign, namely the corresponding eigenvalue,  
 \[
 h(\pm)~E_{\sigma_1}\wedge E_{\sigma_2}\wedge\cdots \wedge E_{\sigma_k}=\epsilon(\sigma_1,\sigma_2,\ldots,\sigma_k)~E_{\sigma_1}\wedge E_{\sigma_2}\wedge\cdots\wedge
 E_{\sigma_k}.
 \]
 The graph that has been constructed has vertices $(\sigma_1,\cdots, \sigma_k) \longleftrightarrow
 E_{\sigma_1}\wedge \cdots \wedge E_{\sigma_k}$.   For $k$ odd, we consider  the action $h(-)$, 
 and for $k$ even, the action of $h(+)$.  An edge $\Rightarrow$ exists between two vertices related by a simple reflection whenever the signs (eigenvalues of $h(\pm)$ ) agree for the two elements in $\bigwedge^k V(E_1,\cdots, E_n)$.
 
 \begin{Remark}
In the identification of the Schubert cell $(\sigma_1,\ldots,\sigma_k)$ with the wedge vector $E_{\sigma_1}\wedge\cdots\wedge E_{\sigma_k}$, the action between two cells in the weak Bruhat order
can be considered as the KP flow through the corresponding two dominant exponentials.
Then changing the sign $\epsilon(\sigma_1,\ldots,\sigma_k)$ is equivalent to having a zero in the $\tau$-function. That is, the KP flow has a singularity.  
\end{Remark}

 \begin{Notation} With the diagonal action of $h(\pm)$,  we can refer to this graph denoted by $\mathcal{G}(k,n)$ as the graph  associated to the action of $h(\pm)$ on $\bigwedge^k V(E_1,\cdots , E_n)$.  It will be shown that this graph is an incidence graph computing cohomology of $\Gr(k,n)$ relative to Schubert cells.
  \end{Notation}

We now state the main theorem.
 \begin{Theorem}\label{mainth1}The graph ${\mathcal G}(k,n)$   of $h(\pm)$ acting on $\bigwedge^kV(E_1,\cdots, E_n)$  agrees with the incidence graph of $\Gr(k,n)$.
 \end{Theorem}

We prove the Theorem in the section \ref{teorema}. Before closing this section, we give some lower dimensional examples.
\begin{Example} \label{ex2} Let us first consider $\Gr(2,4)$:
With the signs $(\epsilon_1,\ldots,\epsilon_4)=(+,+,-,-)$, we obtain
\[\begin{matrix} 
  (1,2)       &  \rightarrow &     (1,3)          &  \Rightarrow      &  (1,4)                    \\
                 &                              & \Downarrow &                                    &   \Downarrow       \\
                  &                             &      (2,3)          &   \Rightarrow     &      (2,4)          \\
                  &                             &                         &                                   &  \downarrow  \\
                  &                              &                       &                                     &       (3,4)
  \end{matrix}
\]
This is the incidence graph of $\Gr(2,4)$ and the nonzero incidence numbers are $\pm 2$.  This graph is, of course,  a subgraph of the incidence
graph for the real flag manifold  shown in  Figure \ref{fig:G24}. The integral cohomology $H^*(\Gr(2,4),\mathbb{Z})$ is then given by
\[
\begin{array}{llll}
H^0(\Gr(2,4),\mathbb{Z})&=\mathbb{Z},  \hskip 1.5cm   &H^3(\Gr(2,4),\mathbb{Z})&=\mathbb{Z}_2,\\
H^1(\Gr(2,4),\mathbb{Z})&=0,                                          &H^4(\Gr(2,4),\mathbb{Z})&=\mathbb{Z}, \\
H^2(\Gr(2,4),\mathbb{Z})&=\mathbb{Z}_2,
\end{array}
\]
It is interesting to note that the Betti numbers $\beta_0=1$ and $\beta_4=1$ are coming from the cells with $(1,2)$ and $(3,4)$, respectively.  The corresponding Young diagrams for those cells are given by
\[
(1,2)=\emptyset,\qquad (3,4)=\young[2,2][5].
\]
Here the Young diagram associated to the cell $(\sigma_1,\ldots,\sigma_k)$ is defined by 
$(\nu_1,\ldots,\nu_k)$ with $\nu_j=\sigma_{k-j+1}-(k-j+1)$ (note $\nu_j\ge\nu_{j+1}$
as the usual definition of the Young diagram, and each $\nu_j$ expresses the number of boxes in the row of the diagram).

\end{Example}

\begin{Remark}  \label{rpnminusk} Note that the top row of the graph corresponds to the elements

$$ E_1\wedge E_2~ \rightarrow~ E_1\wedge E_3~ \Rightarrow~ E_1\wedge E_4.$$ These terms $E_2, E_3, E_4$ already
appeared in the decomposition \eqref{dd} in Example \ref{gr24}  indicating the position of a pivot. They now  corresponds to
a portion of the graph of $\Real P^3$  with twisted coefficients, which agrees with the graph of  $\Real P^2$ with
constant coefficients,  and is  indicated with $\overbrace{\cdots}$,
$$E_1 ~\Rightarrow~ \overbrace {E_2~ \rightarrow~ E_3~ \Rightarrow~ E_4} .$$
In general the top row of the graph is one of the following depending on whether $k$ is odd or even:
\begin{enumerate}
\item [(a)]   If $k$ is odd, we have $h(-)$ action, and 
\[
E_1\wedge  \cdots \wedge E_{k-1} \wedge E_k \, \Rightarrow \, E_1 \wedge  \cdots  \wedge E_{k-1} \wedge E_{k+1}\,\rightarrow \, E_1\wedge\cdots  \wedge E_{k-1} \wedge  E_{k+2}\, \Rightarrow \,\cdots 
\]
\item [(b)]  If $k$ is even, we have $h(+)$ action, and 
\[ E_1\wedge  \cdots \wedge E_{k-1}  \wedge E_k \, \rightarrow\,  E_1\wedge  \cdots \wedge E_{k-1} \wedge E_{k+1}\,\Rightarrow\,  E_1\wedge\cdots \wedge E_{k-1} \wedge  E_{k+2} \,\rightarrow\, \cdots 
\]
\end{enumerate}
\end{Remark}

 We summarize the argument above as the following Lemma.

\begin{Lemma} \label {lemmatop}  We have the following. 
\begin{enumerate}
\item [(i)] The top row of the graph  ${\mathcal G}(k,n)$ associated to the action of $h(-)$ on $\bigwedge^kV(E_1,\cdots E_n)$ is  the incidence graph of $\Real  P^{n-k}$ with trivial coefficients if $k$ is odd and with twisted coefficients if $k$ is even.
\item [(ii)]The top row of the graph ${\mathcal G}(k,n)$ associated to the action of $h(+)$ on $\bigwedge^kV(E_1,\cdots E_n)$ is  the incidence graph of $\Real  P^{n-k}$ with twisted  coefficients if $k$ is even and with trivial  coefficients if $k$ is odd.
\end{enumerate}
\end{Lemma}

\begin{Example} We now consider $\Gr(3,6)$: The action is then given by $h(-)={\rm diag}(+,-,-,+,+,-)$, and
the graph is given by
\[\begin{matrix} 
  (1,2,3)_0 \quad  &     (1,2,4)_1        &  \quad\Rightarrow  \quad  &  (1,2,5)_2                & &       &  (1,2,6)_3               &          &               &     \\
                 & \Downarrow   &                           &  \Downarrow          &             &         & \Downarrow      &            &  			&                      \\
                  &     (1,3,4)_2          & \Rightarrow     &    (1,3,5)_3            & (1,4,5)_4      &      & (1,3,6)_4                &(1,4,6)_5 & \Rightarrow& (1,5,6)_6    \\
                   &     &                           &        &             &    &                                                 &                           & &                                      \\
                    &   (2,3,4)_3            &  \Rightarrow    &   (2,3,5)_4             &   (2,4,5)_5     &      & (2,3,6)_5              & (2,4,6)_6  & \Rightarrow & (2,5,6)_7\\
                      &                        &                            &                             &  \Downarrow &  &                      &\Downarrow&                &\Downarrow        \\
                         &                      &                           &                                &  (3,4,5)_6       &    &                       & (3,4,6)_7     &\Rightarrow & (3,5,6)_8   \\
                           &                 &                           &                                  &                   &      &                          &               &                      &          \\
                	&	&			&			 	&		      &   &		         &		&			&(4,5,6)_9
                  \end{matrix}
\]
Here the single arrows (corresponding to the zero incidence numbers) are all eliminated. The suffix in each cell shows the
dimension of the cell, that is, the dimension $d$ is given by $d=\sum_{j=1}^3(\sigma_j-j)$ for 
each cell of $(\sigma_1,\sigma_2,\sigma_3)$. The cohomology is then given by
\[\begin{array}{llllll}
H^0(\Gr(3,6),\mathbb{Z})=\mathbb{Z}, \hskip 3cm    & H^5(\Gr(3,6),\mathbb{Z})=\mathbb{Z}, \\[1.0ex]
H^1(\Gr(3,6),\mathbb{Z})=0,                                      & H^6(\Gr(3,6),\mathbb{Z})=\mathbb{Z}_2\oplus\mathbb{Z}_2, \\[1.0ex]
H^2(\Gr(3,6),\mathbb{Z})=\mathbb{Z}_2,                        &H^7(\Gr(3,6),\mathbb{Z})=\mathbb{Z}_2,   \\[1.0ex]
H^3(\Gr(3,6),\mathbb{Z})=\mathbb{Z}_2,                                                                     & H^8(\Gr(3,6),\mathbb{Z})=\mathbb{Z}_2, \\[1.0ex]
H^4(\Gr(3,6),\mathbb{Z})=\mathbb{Z}\oplus\mathbb{Z}_2\oplus\mathbb{Z}_2,    &H^9(\Gr(3,6),\mathbb{Z})=\mathbb{Z}. \\[1.0ex]
\end{array}
\]
Note here that the Betti numbers $\beta_0=1, \beta_4=1,\beta_5=1$ and $\beta_9=1$ are coming from
the Schubert cycles with $(1,2,3), (1,4,5), (2,3,6)$ and $(4,5,9)$, and the Young diagrams of those cycles
are given by
\[
(1,2,3)=\emptyset,\qquad (1,4,5)=\young[2,2][5],\qquad (2,3,6)=\young[3,1,1][5],\qquad (4,5,6)=\young[3,3,3][5].
\]
The Young diagram for $(4,5,6)$ may be considered to be a combination of $(1,4,5)$ and $(2,3,6)$.
This pattern is common for the case of $\Gr(k,n)$ with $k=$odd and $n=$even (see Section \ref{Ppolynomial})
\end{Example}

%%%%%%

\section{Proof of the main Theorem for the incidence graph} \label{teorema}

Here we begin to provide a roadmap of the main argument to prove the main theorem \ref{main1}
through a simple example.

\subsection{The case of $\Gr(2,4)$}\label{ex3} 
Let us describe the example of $\Gr(2,4)$:
First we recall that  the  graph $\G(2,4)$ is induced by the KP signs $(+,+,-,-)$ (i.e. $h(+)={\rm diag}(+,+,-,-)$) by  Lemma \ref  {lemmatop}  or Remark  \ref{rpnminusk}. We start with the decomposition, $\Gr(2,4)=Y(2)\cup Y(3)\cup Y(4)$, and
observe that if we proceed inductively (either on $n$ or on $k$), then pieces of the incidence  graph are already available. In this case we can assume  that we know the graphs  corresponding to the $Y(j) \leftrightsquigarrow \Gr(1,j-1)$ for $ j=2,3,4$ (those are the columns in the graph below).  We also know   the top row, i.e the subgraph  associated to a copy of   $\Real P^2$.  Thus we get the following,
\[\begin{matrix} 
  (1,2)       &  \longrightarrow &     (1,3)          &  \Longrightarrow      &  (1,4)                    \\
                 &                              & {\Downarrow}  &                                    &   {\Downarrow}       \\
                  &                             &      (2,3)          &   {?}     &      (2,4)          \\
                  &                             &                         &                                   &  {\downarrow} \\
                  &                              &                       &                                     &       (3,4)
  \end{matrix}
\]
where we are still missing the edge indicated with ``\,?\,''.
We  note two things:
\begin{enumerate}
\item [(a)]  The incidence  graphs associated to the columns $Y(j)$ correspond to {\em twisted } coefficients.
\item [(b)] To determine  the  missing edge, we need to show that the
arrows on the top row extend to the columns (are \lq\lq constant \rq\rq along the columns).  The missing edge is  $\Rightarrow$.
\end{enumerate}
We explain (a) and (b) in  the rest of this section, which provides the proof of Theorem \ref{main1}. However (a) requires some notation from \cite{casian:06} 
to encode the structure of $K$-equivariant local systems on a flag manifold.
Keeping in mind this simple example, the roadmap of the proof  consists in giving  a complete description of the
local systems associated to the line bundles  ${\mathcal E}(j)$.  The description takes place in the context of  $K$-equivariant line bundles  on the flag manifold ${\mathcal B}$  for $K={\rm O}_n(\Real)$  or 
$K={\rm SO}_n(\Real)$ with respect to (a).  After introducing some notation and  Proposition \ref{propsigns} below
on the sign change under the Weyl  group action,  the general situation of $\Gr(k,n)$ becomes an issue of bookkeeping.  The bookkeeping is done through 
the Toda  signs $\tilde \epsilon_i$  which were introduced in \cite{casian:06} and above  in subsection \ref {remarktoda}.  What determines the (twisted)  coefficients
in the incidence graph of Grassmanians along the columns, that is item (a),   is the structure of the vector bundle described in Example \ref {gr24}  in terms of the projections $\pi_3, \pi_4$ and
corresponding determinant  (line) bundles. The vector bundles,  roughly speaking, have fibers corresponding to the  cells of the  $\Real P^2$ along the top row of the graph $\G(2,4)$.

%labibi

% we revisit the Example of $\Gr(2,4)$ one last time in Example \ref{revisit}

\subsection{Some standard notation}
Let $G$ be a real split semi-simple Lie group associated to 
 the real Lie algebra ${\mathfrak g}$. For this paper the relevant cases will be $G=\SL_n(\Real)^{\pm}$ or $G=\SL_n(\Real)$ but some of the statements in this section apply to the more general situation.
We fix $H$ a  split Cartan subgroup of $G$  with Lie algebra ${\mathfrak h}$, $B=HN$ a Borel subgroup and $P$ a maximal parabolic subgroup containing $B$.   We let $K$ denote a maximal compact Lie subgroup with Lie algebra ${\mathfrak  k}$,  $T=K\cap H$  is a finite subgroup of $H$
(usually $K={\rm O}_n(\Real)$, or $K={\rm SO}_n(\Real)$ here and $T$ the diagonal $n\times n$ matrices  with entries $\pm 1$).

Let  $\{h_{\alpha_i},e_{\pm\alpha_i}\}$ is the Cartan-Chevalley basis of $\mathfrak g$ with the simple roots
$\Pi=\{\alpha_1,\cdots,\alpha_l\}$
which satisfy the relations,
\[
    [h_{\alpha_i} , h_{\alpha_j}] = 0, \quad
    [h_{\alpha_i}, e_{\pm \alpha_j}] = \pm C_{j,i}e_{\pm \alpha_j} \ , \quad
    [e_{\alpha_i} , e_{-\alpha_j}] = \delta_{i,j}h_{\alpha_j},
\]
where $(C_{i,j})$ is the $l\times l$ Cartan matrix of $\mathfrak g$.

We first review the computation of integral cohomology of  $G/B$ with  $K$-equivariant local coefficients: 
Let us recall that there is  a filtration by Bruhat cells with ${\mathcal B}_j:=\sqcup_{l(w)\le j}{NwB/B}$,
\[ 
\emptyset\, \subset \, {\mathcal B}_0\, \subset  \, {\mathcal B}_{1} \, \subset \, \cdots\, \subset \, {\mathcal B}_{l(w_o)}= {G/B}
\]
where $w_o$ indicates the longest element of the Weyl group $W$.
 We  have  coboundary maps,
$\delta:H^s({\mathcal B}_s, {\mathcal B}_{s-1}; {\mathbb Z}) \to 
H^{s+1}({\mathcal B}_{s+1}, {\mathcal B}_{s};{\mathbb Z})$,  and these define
a chain complex which computes the cohomology of $G/B$.

Recall that $\Gr(k,n)=G/P_k$ for a maximal parabolic $P_k$ and that the Weyl group of its Levi factor is denoted ${\mathcal P}_k$.   On  the level of cells parametrized by the Weyl group,
$W=\mathcal{S}_n$ in this case, there is a bijection, $\mathcal{S}^{(k)}_n\rightarrow W/{\mathcal P}_k$.

\begin{Notation} \label{flag} Flag manifold and the Bruhat cells $ {\mathcal B}_{w}$:
Consider the flag manifold for $G=\SL_n(\Real)$ or  $\SL_n(\Real)^{\pm}$,   $G/B$ with $B=HN$,  consisting of all real flags 
$\{ 0\subset V_1 \subset V_2 \subset \dots\subset V_n= \Real^n \}$.  Let   ${\mathcal B}_{w}$ denote the $N$ orbit 
  $N{w}B/B$.  Hence  $B_{w} \cong w N \cap N^-$.  
There is  a decomposition into Bruhat cells, i.e. into $N$ orbits,  $G/B=\bigsqcup_{w \in W} {\mathcal B}_{w} $.  
\end{Notation}

 We now recall that in the cases of $G=\SL_n(\Real)$ or $G=\SL_n(\Real)^{\pm}$  studied here,  $S^{(k)}_n $ consists of representatives of cosets  in $W/{\mathcal P}_k$ of minimal length. 

\begin{Remark} \label{subgraph}
 
The subspace consisting of  cells 
 $X(k,n):= \bigsqcup_{w\in S^{(k)}_n } {\mathcal B}_{w}$   is   {\em homeomorphic} to $G/P_k\cong \Gr(k,n)$ because the projection $\pi: G/B\to G/P_k$
is such that $\pi$ restricted to each $ {\mathcal B}_{w}$ is a bijection whenever ${w}\in S^{(k)}_n \subset W$.  The fibers of this projection are real flag
manifolds associated to the Levi factor and  $ {\mathcal B}_{w}$ intersected with this fiber  is one point, the lowest dimensional Bruhat cell in this 
flag manifold. On the level of cells this corresponds to restricting the quotient
$W\rightarrow W/ {\mathcal P}_k $ to the subset  $S^{(k)}_n$ which parametrizes $W/ {\mathcal P}_k$. This explains why the incidence graphs of  real Grassmanians
are found as subgraphs of the incidence graph of the real flag manifold.

\end{Remark}

\subsection{Connection of  the cohomology of $G/B$  with Hecke algebra operators}\label{connection}

 Here we give a quick summary of \cite {casian99} on the cohomology of $G/B$ as  reformulated
in \cite{casian:06}.

 For the purposes
of this paper we consider $K={\rm SO}_n(\Real)$, its complexification $K_{\mathbb C}$ as well as
$G=\SL_n(\Real)$ and $G_{\mathbb C}={\rm SL}_n({\mathbb C})$.  However this  can be done in the more general context of  \cite{lusztig83}. We 
consider the real flag manifold ${\mathcal B}=G/B$ and its complexification ${\mathcal B}_{\mathbb C}=G_{\mathbb C}/ B_{\mathbb C}$.
For example for $G=\SL_n(\Real)$,   ${\mathcal B}_{\mathbb C}$ consists of ${\mathbb C}P^1$. The real
flag manifold $\mathcal{B}$ is contained as a circle inside the open $K_{\mathbb C}$ orbit ${\mathbb C}\setminus \{ 0 \}$.

Recall that the Hecke algebra  ${\mathcal H}$ of the Weyl group defined in  \cite{kazhdan:80} is  a deformation of the usual group algebra of $W$.
As a set it is given by
 ${\mathcal H}={\mathbb Z}[q, q^{-1}]\otimes_{\mathbb Z} {\mathbb Z}[W]$, that is, the set 
  of formal linear combinations of elements in $W$  with coefficients in 
 ${\mathbb Z}[q, q^{-1}]$.  The multiplication is defined as in p. 189 of  \cite{kazhdan:80}.
The elements $w\in W$ are denoted by $T_w$ when viewed inside  ${\mathcal H}$
and we have $T_xT_y=T_{xy}$ when $l(xy)=l(x)+l(y)$
and  for any simple reflection $s_i$, we have $(T_{s_i}+1)(T_{s_i}-q)=0$. This replaces the equation $s_i^2=1$.

If ${\mathcal D}$ denotes the set of $K$-equivariant local systems on the  flag manifold 
${\mathcal B}_{\mathbb C}$  (dextended by zero from local systems on $K_{\mathbb C}$ orbits) then
${\mathbb Z}[q,q^{-1}] \otimes _{{\mathbb Z}}{\mathbb Z} [{\mathcal D}]$ becomes a module over the Hecke algebra.  The incidence graph
of the real flag manifold with local coefficients ${\mathcal L}_o$  can be described in terms of this  ${\mathcal H}$ action on local systems. If 
${\mathcal D}_o \subset {\mathcal D}$ consists  of  the elements in ${\mathcal D}$ which are supported on the open dense orbit, then the  incidence graph of $G/B$ with local coefficients 
 is given by the relations between the various elements
 $T_{\tau}^{-1}{\mathcal L}$.  Let  $w_o$ be  the longest element in $W$, and  consider
 $T_{w_o}^{-1}{\mathcal L}$. Suppose that ${\mathcal L}_o\in {\mathcal D}_o$ occurs  in
 the expression $T_{w_o}^{-1}{\mathcal L}$ , (or we set $q=1$  and obtain  ${\mathcal L}_o$) then
the incidence graph  that we are describing  corresponds to  $H^*(G/B; {\mathcal L}_o)$.   The vertices correspond to  $W$.  
We  associate  \lq\lq graded characters \rq\rq  to elements of $W$ as follows. We let
$\theta (e)= T_{w_o}^{-1}{\mathcal L}$, and then set $T_\sigma  \theta (e)=\theta (\sigma)$.  The element  $\theta (\sigma)$ 
corresponds to  cells parametrized by  $\sigma^{-1}=w$.   Each graded character $\theta (\sigma)$  corresponds to a local system ${\mathcal L}_w \in {\mathcal D}_o$. For example if we set $q=1$ then  $\theta (\sigma)$
reduces to  ${\mathcal L}_w$. Let  $q^R$  denote  the power of $q$ of ${\mathcal L}_e$
in   $\theta (e)$. We  let  $q^{\eta (w) }$ be the power of $q$ of  ${\mathcal L}_w$ in
 $q^R\theta (\sigma)$. We have  readjusted so that $\eta (e)=0$  and all the $\eta (w)$ are non-negative integers. 
 \begin{Remark}
 In \cite {casian:06},
 the number  $\eta (w)$ is   described in terms of the number of blow-ups of the Toda flow.
Namely, we count the number of singular points in the flow from the top cell (corresponding to the flow
for $t\ll 0$) to the Bruhat cell marked by $w\in W$ (see Figure \ref{fig:G24}).
\end{Remark}

Now the following is in the case when $G$ is  ${\mathbb R}$ split e.g. $G={\rm SL}_n({\mathbb R})$.  We can describe the edges of the incidence graph in the following way (after \cite {casian:06}):
We have an edge $w \Rightarrow w^\prime$
whenever $w \le w^\prime$ in the Bruhat order, $l(w^\prime)=l(w)+1$ and  $\eta(w)=\eta(w^\prime)$.

 \begin{Example}\label{zz1}
In the case of $G={\rm SL}_2({\mathbb R})$,  
 ${\mathcal D}= \{ {\mathcal C},{\mathcal L},  \delta_+, \delta_-  \}$ where  ${\mathcal C}$  denotes a trivial 
 sheaf  on ${\mathcal O}_o={\mathbb C}^*$, $\delta_{\pm }$ are sheaves supported on  the points $0$, 
 $\infty$ respectively,  and ${\mathcal L}$ is a non-trivial local system on ${\mathbb C}^*$. We have  
 $T_{s_1} {\mathcal C} = (q-2){\mathcal C} +(q-1)(\delta_-+\delta_+)$ and $T_{s_1}{\mathcal L}=-{\mathcal L}$.    
 
 In the Hecke algebra $T_{s_1}^{-1}= q^{-1}(T_{s_1} + (1-q))$. Hence by  applying 
  $T^{-1}_{s_1}$ to ${\mathcal C}$ we obtain  $-q^{-1} {\mathcal C} + q^{-1}(q-1)(\delta_++\delta_-)$. 
  We also have $ q^{-1}(T_{s_1}{\mathcal L} + (1-q){\mathcal L})= q^{-1}(-{\mathcal L}+ (1-q){\mathcal L})=-{\mathcal L}$.    Hence if we start with ${\mathcal C}$, then  $\theta (e)=-q^{-1} {\mathcal C} + q^{-1}(q-1)(\delta_++\delta_-)$, $\theta (s)={\mathcal C}$. Then (setting $q=1$ in $\theta (e)$) ${\mathcal L}_o={\mathcal C}$
  corresponding to constant coefficients. Since the power of $q$ associated to ${\mathcal C}$ is , respectively, $-1$ and $0$, or shifting to get non-negative integers $\eta(e)=0, \eta (s)=1$. Since
 $\eta(e)\not=\eta(s)$  there is no edge $\Rightarrow $ joining $e$ and $s$. This situation corresponds to the existence
 of one blow up in the Toda lattice.  We end up with and incidence graph containing to
  vertices $e$ and $s$ and no edges $\Rightarrow$. This gives the cohomology of $\Gr(1,2)$, a circle,
  with constant coefficients.  If we consider ${\mathcal L}$ then we end up with $e\Rightarrow s$
  the incidence graph of $\Gr(1,2)$ with local coefficients ${\mathcal L}$. This second case corresponds to
  an irreducible principal series module or to the case in which there are no blow-ups in the 
  Toda flow.
  
   Finally, we note the connection with the representation theory of  ${\mathfrak sl}(2;{\mathbb R})$.  By rewriting   $T^{-1}_{s_1}q^1{\mathcal C}$ as
  $-q^{-1/2}\hat C_{\mathcal C} +\hat C_{\delta_{+}} +\hat C_{\delta_{-}}$ with $\hat C_{\mathcal C}=q^{-1/2}( {\mathcal C}+\delta_{+}+\delta_{-})$,  $\hat C_{\delta_{\pm}}=\delta_{\pm}$,  we recover the weight filtration of the principal series representration containing the trivial representation
  as submodule  (replace $\hat C_{\mathcal C}$ with a trivial representation $C$ and $\hat C_{\delta_{\pm}}$ with two discrete series representations  $D_{\pm}$.  
   \end{Example}
 
We can now describe the incidence graph $\G(k,n)$ in terms of the description given above using  Hecke algebra operators. The incidence graph
is a graph with  vertices  $S_n^{(k)} \subset W$.  Only two $K_{\mathbb C}$ equivariant  local systems ${\mathcal L}(e)$ are considered. One trivial and one
non-trivial. Then  we have an edge $w \Rightarrow s_iw$
whenever $l( s_iw)=l(w)+1$ and  $\eta(w)=\eta( s_iw)$ (i.e. no sign change).

\subsection{The edge  in the incidence graph associated  to ${\mathcal B}_{w} \cup {\mathcal B}_{s_iw}$}

As in \cite{casian99} and \cite {casian:03} the edge in the incidence graph corresponding to   ${\mathcal B}_{w} \cup {\mathcal B}_{s_iw}$ 
is   encoded in the action
of the  Hecke algebra operator $T_s$ on the graded character $\theta (\sigma)$ , i.e.
whether $\eta (w)=\eta (sw)$ or  $\eta (w)\ne\eta (sw)$.   In terms of the action of the Weyl group on the set of signs 
of the Toda flow corresponding to $sgn(a_i)$ (introduced in  \cite{casian:06}  and below in Definition \ref{signos}),  there is an edge $\Rightarrow$, precisely when $s_i$ does not change signs.
If signs change,  we just indicate the relation between $w$ and $s_iw$ ion the Bruhat order with $\to$, which, does not count
as an edge in the graph, and which corresponds to crossing $\tau_i=0$, a blow-up in the Toda flow.

 %Since we are only studying a bundle over a circle: ${\mathcal B}_{w} \cup {\mathcal B}_{s_iw}$,  we only need  the restriction of this line bundle to 
%a circle corresponding to the simple reflection $s_i$  relative to $w(\Delta^+)$ i.e. relative to $ws_i$.

We recall that each simple root ${\alpha_i}$ gives rise to a Lie group
homomorphism $\Phi_{\alpha_i}: \SL_2(\Real)\to G$. We denote
by $z_i=\Phi_{\alpha_i}( \begin{pmatrix} -1&0\\
      0&-1  \end{pmatrix}       )$.  This element can be expressed as
      $z_i =\exp(\pi \sqrt{-1} h_{\alpha_i})$. We can also considered
      $\sigma (z_i)=\exp(\pi \sqrt{-1} \sigma h_{\alpha_i})$ associated to $\sigma s_i$ i.e. to $\sigma \alpha_i\in \sigma (\Delta^+)$.
      The sign of $\chi_w (  \sigma (z_i)   )=\pm 1$  then determines whether $w \rightarrow  s_i w$ (coboundary is zero) or $ w\Rightarrow s_i w$ (coboundary is $\pm 2$).

\subsection{$K$-equivariant local systems on $G/B$  and their  Toda signs}\label{Todasigns}

 Borel subgroups  $w(B)$ containing $H$ and corresponding to the positive  roots systems $w(\Delta^+)$   define points $x_w$ in the flag manifold $G/B$. These points are, in turn,  representatives of 
Bruhat cells. We follow the notation in \cite{casian99} which allows us to describe the action of the Weyl group on these local systems.  

\begin{Definition}

Consider a  $K$-equivariant local system ${\mathcal L}$  determined by  a character of $T$ given by $\chi({\mathcal L})$.  We consider the local system  ${\mathcal L}_{w}$
given by the $T$ character  $\chi({\mathcal L} ) \otimes \bigwedge^{l(w)} {\sigma }{\mathfrak n }\cap{ \mathfrak n}^*$ ($\sigma=w^{-1})$.  These correspond to $K_{\mathbb C}$-equivariant local systems supported on the open $K_{\mathbb C}$ orbit on ${\mathcal B}_{\mathbb C}$ and have already been described above in terms of the
Hecke algebra action.

\end{Definition}

Given $w \in W$ and the $K$-equivariant local system ${\mathcal L}_w$, we associate some signs which describe the local system from the perspective of simple reflections relative to $\sigma (\Delta^+)$.
As in  Definition 4.3 of   \cite{casian99}  to  the $K$-equivariant  local system ${\mathcal L}_{w}$ one associates a list of signs $\epsilon(\sigma, {\mathcal L})= (\tilde\epsilon_1, \cdots, \tilde\epsilon_l)$. These signs  describe  the local system relative to $x_w$ i.e. relative to simple roots in $\sigma (\Delta^+)$. Here $l$ the rank of the semi-simple Lie algebra. The signs keep track of triviality or non-triviality along directions corresponding to simple roots for $\sigma (\Delta^+)$.  A sign $-$ in the $i$th place means that along  a circle associated to $s_i$ the local system is trivial and a sign $+$ means that it is non-trivial.

\begin{Definition} \label{signos} \par
Let $\chi$ be a character of $T$ corresponding to the local system ${\mathcal L}$. Then ($w^{-1}=\sigma$) ${\mathcal L}_w$ corresponds to 
$\chi_w=\chi( {\mathcal L} ) \otimes \bigwedge^{l(w)} \sigma {\mathfrak n }\cap{ \mathfrak n}^*$.  We let $\epsilon(\sigma, {\mathcal L})=(\tilde\epsilon_1,\cdots, \tilde\epsilon_l)$  where
$\tilde\epsilon_i=-\chi_\sigma (z_i)$. This is $-$ if the local system is trivial along the direction
of the $i$th simple root relative to $\sigma \Delta^+$ and $+$ otherwise. For convenience we will refer to these signs as  {\em Toda signs}, 
since they have a clear interpretation in terms of the Toda flow. When $w=e$ we can refer to the Toda sign as {\em initial} Toda sign. This initial
Toda sign determines which local coefficients are being used in cohomology computations.

 \end{Definition}
 
 \begin{Example} Consider the flag manifold $\mathcal{B}=G/B$ of $G=\SL_3(\Real)$. 
 Recall the  four such characters, $\chi$ with
$\chi( \begin{pmatrix}
e_1 & 0 & 0\\
0& e_2& 0\\
0&0& e_3
\end{pmatrix} ) =e_i$ for $i=1,2,3$ and the trivial character.
Assume $i=1$.   We compute $\tilde\epsilon_i=-\chi_w(z_i)$ where $w=e$.
We have \[
z_1=\begin{pmatrix} -1 & 0 & 0\\ 
0 &-1 & 0\\
0&0&1
\end{pmatrix}\qquad {\rm and}\qquad z_2=\begin{pmatrix}  1 & 0 & 0\\ 
0 &-1 & 0\\
0&0&-1
\end{pmatrix}.
\]
If we evaluate   $-\chi_e (z_1)$ we get  $+$,  and if we evaluate $-\chi_e (z_2)$   
 we get  $-$.  Hence the two signs associated to ${\mathcal L}$  are $(+,-)$. 

We can now compute the correspondence:

$$\begin{matrix}
\epsilon (e,{ \mathcal L}(\tilde\epsilon_2))  &=& (+,+) \\
\epsilon (e, {\mathcal L}(\tilde\epsilon_3)) &=& (-,+) \\
\epsilon (e, {\mathcal L}(1) )&=& (-,-) \\
\epsilon (e, {\mathcal L}(\tilde\epsilon_1)) &=& (+,-) \\
\end{matrix}$$

 \end{Example}
% Note that I  the   ${\mathcal L}_\sigma$ in this paper would be denoted ${\mathcal L}_{\sigma^{-1}}$ in the inventiones paper
 
\subsection{Computation of the  signs $\epsilon(\sigma, {\mathcal L})$} \label{epsilon} 
 
We  first  write the  well-known action of the Weyl group on roots and coroots:
\begin{align*} 
%\label{cartan1}
s_i \alpha_j &= \alpha_j - C_{j,i} \alpha_i, \\[1.5ex]
% \label{cartan2}
s_i \check \alpha_j &= \check \alpha_j - C_{i,j} \check \alpha_i
\end{align*}
 In the first formula the matrix that is involved is the {\em  transpose}
of the Cartan matrix $(C_{i,j})$. In the second, using coroots, the matrix that appears
is the usual Cartan matrix.
Recall that each simple root $\alpha_j$  gives rise to a Lie group homomorphism
$\Phi_{\alpha_j}: \SL_2(\Real)\to G$ corresponding to an injection of ${\mathfrak{sl}}_2(\Real)$ 
into the Lie algebra.
The elements  $h_{\alpha_j}$  correspond to the $ \left(
\begin{matrix}
1 & 0   \\
0 & -1 \\
\end{matrix}\right)$ in  the  copy of ${\mathfrak{sl}}_2(\Real)$ associated to $\alpha_j$ and corresponds to
the coroot $\check \alpha_j$.  
The reason to introduce the $h_{\alpha_j}$ is that  we then can write  $\Phi_{\alpha_j} ( \left(
\begin{matrix}
-1 & 0   \\
0 & -1 \\
\end{matrix}\right))=\exp( \pi \sqrt{-1}\, h_{\alpha_j} )$.
 From the correspondence between $\check \alpha_i $ and $h_{\alpha_i}$, we have
\begin{equation} 
\label{cartan3}
s_i h_{\alpha_j} =h_{ \alpha_j} - C_{j,i} h_{\alpha_i}.
\end{equation}
The following will also be useful:
\begin{equation*}  \label{wsequation}
\sigma s_i \Delta^+ \cap \Delta^- = \sigma \Delta^+\cap \Delta^- \cup \{ -\sigma \alpha_i \}
\end{equation*}
which corresponds to  
 \begin{equation} \label{lie1equation}
 \sigma s_i{\mathfrak n}\cap   {\mathfrak n }^*=
 \sigma {\mathfrak n}\cap   {\mathfrak n }^*\oplus \Real X_{-\sigma\alpha_i}.
  \end{equation}
By applying  the inverse of $\sigma s_i$ we obtain 
 \begin{equation*} \label{lieequation}
  {\mathfrak n}\cap s_i \sigma ^{-1}  {\mathfrak n }^*=
 s_i {\mathfrak n}\cap \sigma ^{-1}  {\mathfrak n }^*\oplus \Real X_{\alpha_i}.
  \end{equation*}
 This version was used in  \cite{casian:06} but the 
$s_i$  on the right hand side of  the equation above was inadvertently left out  in the Proof of  Proposition 5.1 in \cite{casian:06}. We  then can correct the simple proof of  Proposition 5.1 in \cite{casian:06}:
\begin{Proposition}   \label{propsigns} For  $\epsilon(\sigma , {\mathcal L})=(\tilde\epsilon_1,\ldots,\tilde \epsilon_l)$,
we have $\epsilon(s_i\sigma ,{\mathcal L})=(\tilde\epsilon'_1,\ldots,\tilde \epsilon'_l)$ where
$\tilde\epsilon'_j =\tilde\epsilon_j\tilde\epsilon_i^{C_{i,j}}$. 
\end{Proposition}
\begin{Proof}
We set $\sigma=w^{-1}$ and  start with \eqref {lie1equation}:
 $ \sigma s_i{\mathfrak n}\cap   {\mathfrak n }^*=
 \sigma  {\mathfrak n}\cap   {\mathfrak n }^*\oplus \Real X_{-\sigma \alpha_i}$
where $X_{\alpha_i}$ is a  root vector associated to the simple root $\alpha_i$.  Therefore
\begin{equation} \label{correct1}
\bigwedge^{l(w)+1} \sigma s_i {\mathfrak n}\cap {\mathfrak n}^*=\bigwedge ^{l(w)} (\sigma {\mathfrak n}\cap {\mathfrak n}^*) \wedge  {\mathbb R} X_{-\sigma \alpha_i}
\end{equation}
In terms of characters of $T$ we then have :
$\tilde\epsilon^\prime _j = -\chi_{s_iw} (\exp( \pi \sqrt{-1}  ws_i  h_{\alpha_j} ))$ (by definition).  Using \eqref{correct1}  this becomes  $ -( \chi_w \chi_{-\sigma \alpha_i})  (\exp( \pi \sqrt{-1} h_{\sigma \alpha_j} )).$  We evaluate $- \chi_w  (\exp( \pi \sqrt{-1} \sigma s_i h_{\alpha_j}) )$ and using \eqref{cartan3}, this becomes\begin{align*}
&-\chi_w  (\exp( \pi \sqrt{-1}\sigma h_{ \alpha_j} ))(\chi_w (\exp( \pi \sqrt{-1}\sigma h_{\alpha_i}))^{C_{i,j}} \\
=&-\chi_w  (\exp( \pi \sqrt{-1}\sigma h_{ \alpha_j}) )(-1)^{C_{i,j}}(-\chi_w (\exp( \pi \sqrt{-1}\sigma h_{\alpha_i}))^{C_{i,j}}\\
=&~\tilde\epsilon_j \tilde \epsilon_i^{C_{i,j}} (-1)^{C_{i,j}}.
\end{align*}
Now we consider the additional factor  $ \chi_{-w\alpha_i}  (\exp( \pi \sqrt{-1} \sigma s_ih_{\alpha_j}) )$. This
is just 
\[
\exp(  \pi \sqrt{-1} \langle-\sigma  \alpha_i; \sigma h_{\alpha_j} \rangle)=\exp(  \pi \sqrt{-1} \langle- \alpha_i; h_{\alpha_j} \rangle) = (-1)^{C_{i,j}}.
\]
Hence we obtain $\tilde\epsilon_j \to \tilde\epsilon_j\tilde\epsilon_i^{C_{i,j}}$.
\end{Proof}

\begin{Example} In the case of $\SL_3(\Real)$ we have
 $\epsilon( \sigma , {\mathcal L})=( \tilde\epsilon_1, \tilde\epsilon_2)$ then
 $\epsilon(s_1\sigma,  {\mathcal L})= (\tilde\epsilon_1, \tilde\epsilon_1\tilde\epsilon_2)$ and
 $\epsilon(s_2\sigma,  {\mathcal L})= (\tilde\epsilon_1\tilde\epsilon_2, \tilde\epsilon_1)$.

In the general case of  $\SL_n(\Real)$, this has a very simple description and defines an action of the Weyl group on a set of signs,
\[
(\tilde\epsilon_1,\cdots, \tilde\epsilon_{j-1},\tilde\epsilon_j, \tilde\epsilon_{j+1}, \cdots,\tilde{\epsilon}_{n-1} )\quad \mapright{s_j}\quad
(\tilde\epsilon_1,\cdots, \tilde\epsilon_j \tilde\epsilon_{j-1},\tilde\epsilon_j, \tilde\epsilon_j \tilde\epsilon_{j+1} ,\cdots,\tilde{\epsilon}_{n-1}).
\]

\end{Example}

\begin{Example} \label{revisit}
 In the case of $\Gr(2,4)$ we consider the cells  $\sigma\in\{ e, s_2, s_3s_2\} $  which
correspond to a copy of  $\Real P^2$ (as in Remark \ref{remark0}  item (ii)) .   We compute the structure
of the line bundles ${\mathcal E}(j)$ defined in   \eqref {thebundles}  and associated to the vector bundles
 ${\mathcal V}(j)$.  As  in Remark \ref{remark0} the fibers of  ${\mathcal V}(j)$ correspond
 to the cells of  $\Real P^2$.

 We start  with   the  Toda sign $\epsilon( \sigma , { \mathcal L}(1) )=(-,-,-)$  and apply  these Weyl group elements $\sigma$  to  it. We obtain
\[
\begin{matrix}
\epsilon( e, { \mathcal L}(1) ) &\mapright{s_2} &  \epsilon( s_2, {\mathcal L}(1) ) &\mapright{s_3}  & 
\epsilon( s_3s_2,{ \mathcal L}(1) ) \\
{\updownarrow}& {}& {\updownarrow} & {} &  {\updownarrow} \\
(-,-,-) &\mapright{s_2}&  (+,-,+)& \mapright{s_3} & (+,-,+)\\
\end{matrix}
\]

The fact that a $+$ appears in the first spot in the second and third of the Toda signs discussed means
that in the incidence graph of $\Gr(2,4)$, {\em twisted} coefficients will have to be used along the columns
i.e. along  $Y(3)\cong \Gr(1,2)$ and $ Y(4)\cong \Gr(1,3)$.

\end{Example}

\begin{Remark} \label{spot}

The general case of $\Gr(k,n)$ is similar.  There are two local systems ${\mathcal  L}$ to consider, that is, the constant  one
 and a second twisted local system.  The signs that correspond are  $\epsilon( \sigma , { \mathcal L})=(\tilde\epsilon_1, \cdots ,\tilde\epsilon_n )$ 
where $\tilde\epsilon_i=-$  for  all $i$ in  the constant coefficient case, and  $\tilde\epsilon_i=-$ for all $i\not=k$,  $\tilde\epsilon_k=+$ for the twisted case.
We now consider the signs obtained by applying the elements $s_k, s_{k+1}, \cdots $.  For convenience we can regard 
$\Gr(k,n)$ as embedded inside $\Gr(k, n+1)$. This simply adds one more $-$ at the end e.g.  $\epsilon( \sigma , { \mathcal L}(1) )=(-,-,-,-)$ instead of
$(-,-,-)$.

Note that in the constant coefficient case  when $s_{k+2r}$  is applied, a new $+$ sign appears in the  $k+2r+1$ spot {\em for the first time}.
All the preceding signs $ s_{m} \cdots s_{k+1}s_k\epsilon( \sigma , { \mathcal L}(1))=(\tilde\epsilon_1^\prime , \tilde\epsilon_2^\prime , \cdots,\tilde\epsilon_n' )$  with  $m \le k+2r-1$, are such that $\tilde\epsilon^\prime_{k+2r+1}=-$. When elements in  $\mathcal{S}_{k+2r}^{(k-1)}  \subset  \langle s_1,\cdots ,s_{k+2r-1} \rangle$ are applied in this  constant coefficient case,  the $+$ sign in the $k+2r+1$ spot remains unchanged by the definition of the Weyl group action on signs.

In the twisted case when $s_{k+2r-1}$ is applied,  a new $+$ sign appears in the  $k+2r$ spot. When  elements in  $\mathcal{S}_{k+2r-1}^{(k-1)} \subset  \langle s_1,\cdots s_{k+2r-1} \rangle$
this $+$ sign remains unchanged.

\end{Remark}

\begin{Example} For example, we have for the local system $\mathcal{L}(1)$,

$$\begin{matrix}
\epsilon( e, { \mathcal L}(1) ) &\longrightarrow &  \epsilon( s_2, {\mathcal L}(1) ) &\longrightarrow  & 
\epsilon( s_3s_2,{ \mathcal L}(1) ) \\
{\updownarrow}& {}& {\updownarrow} & {} &  {\updownarrow} \\
(-,-,-,-) &\mapright{s_2}&  (+,-,+,-)& \mapright{s_3} & (+,-,+,-)\\
\end{matrix}$$
and for $\mathcal{L}^*$,
$$\begin{matrix}
\epsilon( e, { \mathcal L}^*) &\longrightarrow &  \epsilon( s_2, {\mathcal L}^* ) &\longrightarrow  & 
\epsilon( s_3s_2,{ \mathcal L}^* ) \\
{\updownarrow}& {}& {\updownarrow} & {} &  {\updownarrow} \\
(-,+,-,-) &\mapright{s_2}&  (-,+,-,-)& \mapright{s_3} & (-,-,-,+)\\
\end{matrix}$$

\end{Example}

\begin{Remark} \label{involved}

The general structure of all graphs involved can be loosely described by a diagram (see Example \ref {ex3}):

$$\begin{matrix}

{j=k}& {} & {j=k+1}& {} &  {j=k+2} & {} &  {j=k+3}\\
%{\updownarrow}& {} & {\updownarrow }&  {}&  {\updownarrow}&  {} &  {\circ}   \\

 {\circ}& \mapright{s_k}  & {\circ}&  \mapright{s_{k+1}} &  {\circ} & \mapright{s_{k+2}} &  {\circ}  \\
 {\downarrow}& {} & {\downarrow }&  {}&  {\downarrow}&  {} &  { \downarrow }   \\
{\S_{k-1}^{(k-1)}}& {} & {\S_{k}^{(k-1) }}&  {}&  {\S_{k+1}^{(k-1) } }&  {} &   {\S_{k+2}^{(k-1) } }   \\
{\downarrow}& {} & {\downarrow }&  {}&  {\downarrow}&  {} &  { \downarrow }   \\

 {\{e\}}& {}  & { \langle s_1,\cdots, s_{k-1} \rangle}& {}&  {\langle s_1,\cdots, s_{k} \rangle} & {} &  {\langle s_1,\cdots, s_{k+1} \rangle}   \\
 
 %{S_{j-1}^{(k-1)}\subset}&

\end{matrix}$$
Here $\circ$  may correspond to $(\sigma_1,\cdots \sigma_k)  \leftrightsquigarrow  E_{\sigma_1}\wedge \cdots \wedge E_{\sigma_k}$ 
so that  edges $\Rightarrow$  are  defined  by  agreement of eigenvalues of the matrix $h(\pm)$
on $\wedge^kV(E_1,\ldots,E_n)$.  They can also
 be replaced with  $l$ tuples   $\epsilon(\sigma,{\mathcal L})=(\tilde\epsilon_1,\cdots,\tilde\epsilon_{l})$, 
and the edges $\Rightarrow$  then correspond to non-zero co-boundaries.

We then note that for the two possible types of  graphs being considered:
\begin{enumerate}
\item [(i)]  The columns  associated to $j$ (e.g. to $Y(j)$) depend only on the simple reflections $\{ s_1,\cdots, s_{j-2} \}$
\item [(ii)] Any edges $\Rightarrow$ involve vertices in  adjacent columns.
\item [(iii)]  Any  edge $\Rightarrow$  between vertices labeled $j, j+1$  along the
top column (corresponding to  the incidence graph of  $\Real P^{n-k}$)  induces  edges $\Rightarrow$ connecting the cells
 in one column (vertices in  ${\mathcal Gr}(k-1,j-1)$) to the cells in the  next
( corresponding vertices in  ${\mathcal Gr}(k-1,j)$).
\end {enumerate}
\end{Remark}

Note that  in the case of the the graph  ${\mathcal G}(k,n)$  (ii) is true by definition.  We then briefly  discuss (ii) in  the case 
of the incidence graph of $\Gr(k,n)$. Using Definition 3.3 (d)  of  \cite{casian:06} which gives $\Rightarrow$, it is known that whenever there is a coboundary
$\Rightarrow$ relating $\sigma$ and $\sigma^\prime$, the corresponding signs must agree: $\epsilon(\sigma, {\mathcal L})=\epsilon(\sigma^\prime, {\mathcal L})$.
We consider first  the case in which we start with $\tilde \epsilon_i =-$ for all $i$.  We note (Remark \ref{spot}) that when  $s_{j}$  is applied to this Toda sign, either there is a $+$ 
in the $j$ th spot or a new $+$ sign appears in the  $j+1$ spot 
{\em for the first time} and this sign remains along the column  since it is unaffected by the action of $s_1,\cdots s_{j-2}$. 
Therefore the  $\Rightarrow$ can only involve adjacent  ${\mathcal G}_j^{(k-1)}$ i.e. the cells associated to
$Y(j)$ and $Y(j+1)$. The argument in the other case is similar. This establishes  (ii).

Again (iii) is true in the case of  ${\mathcal G}(k,n)$  by definition. In the case of the incidence graph of $\Gr(k,n)$ the connection  between vertices in the top row labeled $j,j+1$
  is  given  by the   $j$th simple  reflection. In terms of signs $\epsilon(\sigma,{\mathcal L})$ it  then depends only on the $j$th sign. However this sign doesn't change along
the $j$th column  (along $Y(j)$).  Therefore an edge $\Rightarrow$ on the top row extends to the entire column.

%\end{Remark}

\subsection{Proof of Theorem} 

The two distinguished $n$-tuples of signs are encoded in two diagonal matrices $h(-)$ and $h(+)$ acting on $V(E_1,\cdots, E_n)$ and all the spaces $\bigwedge^k V(E_1,\cdots E_n)$. An edge $\Rightarrow$ corresponds to the agreement of eigenvalues. The main theorem \ref{main1} can be restated as follows:

 \begin{Theorem}\label{mainth}The graph ${\mathcal G}(k,n)$   of $h(\pm)$ acting on $\bigwedge^kV(E_1,\cdots, E_n)$  agrees with the incidence graph of $\Gr(k,n)$.   More precisely, we have
 the following:
 \begin{enumerate}
 
\item [(a)]  If $k$ is odd and the  choice of signs corresponds to the matrix $h(-)$ $($respectively $h(+))$ then
 ${\mathcal G}(k,n)$ agrees with the incidence graph with constant coefficients $($respectively  non trivial local coefficients$)$.  
 
\item [(b)]  If $k$ is even  and the  choice of signs corresponds to the matrix $h(+)$ $($respectively $h(-))$ then ${\mathcal G}(k,n)$ agrees with the incidence graph with constant coefficients $($respectively  non trivial local  coefficients$)$.  
  
 \end{enumerate}
 
  \end{Theorem}
  
  \begin{Proof}
  We proceed by induction on $k$. We then assume  the statements (a) and (b) of the Theorem for smaller values  of $k$, the cases $k=1$ are already clear. 

There are several similar cases which are handled in identical manner within an inductive proof. We then focus, for simplicity of exposition,  on (b) and  the case when $k$ is even and consider  cohomology of 
$\Gr(k,n)$ with constant coefficients.  Hence we consider the  KP sign $(+,+, -,-,+,\cdots )$ corresponding to $h(+)$
as in Lemma \ref {lemmatop}.  This is the sign that guarantees that the top row of the graph of the action of $h(\pm)$ on 
$\bigwedge^k V(E_1,\cdots, E_n)$ is the incidence graph of $\Real P^{n-k}$ with constant coefficients (Lemma \ref {lemmatop}).

The  two graphs that we are comparing have the same set of vertices and by Remark  \ref{involved} ii) edges involve adjacent columns only.  We need to  show that  there is a non-zero co-boundary  in the incidence graph (given in terms of Toda signs), that is, 
 $ \sigma \Rightarrow s_i\sigma$  appears, exactly when in the graph ${\mathcal G}(k,n)$ (defined in terms of the KP signs, i.e.  edges crossing $\tau_k=0$)  there is a corresponding  edge. 

We start with the decomposition  $\Gr(k,n)=Y(k)\cup Y(k+1)\cup \cdots  Y(j) \cdots $ and as in the  example of $Gr(2,4)$  in subsection  \ref{ex3} observe  that, by induction,  pieces of the incidence  graph are already available and agree with corresponding pieces of the graph $\G(k,n)$.  We can then assume  that we know the graphs  corresponding to the $Y(j) \leftrightsquigarrow \Gr(k-1,j-1), j=k, k+1, \cdots $. We also know   the top row, i.e the subgraph  associated to a copy of   $\Real P^{n-k}$.

We now deal with  issue (a) brought up in  the $Gr(2,4)$ example in subsection
\ref{ex3}.  The spaces $Y(j)$, have some  {\em additional structure  of a vector bundle}   ${\mathcal V}(j)$ over  $\Gr(k-1,j-1)$ (subsection \ref{grknbundle}). Then the  line bundle ${\mathcal E}(j)$  defined in (\ref{thebundles}),  
 $\mathcal{E}(j):=\bigwedge^{j-k}{\mathcal V} (j)$,  determines the local system that we need over $Y(j)$. 
 It turns out that   we need to consider the incidence graphs of the  $Y(j) \leftrightsquigarrow \Gr(1,j-1), j=2,3,\ldots,n$ with {\em twisted} coefficients.   This follows  by noting that the  first  $j-1$ signs of  $ \epsilon( \sigma, { \mathcal L}(1) ) $  on the top row always include one $+$.
 
 Applying   $\sigma = s_k,  s_{k+1}, s_{k+2}, s_{k+3},\cdots  $  as occurring  in Proposition \ref{decompositionS} to the Toda  sign $(-, \cdots ,-)$  we have in the first step
 \[
 \{\tilde\epsilon_i=-~(\forall i)\}\quad\overset{s_j}{\longrightarrow}\quad\{\tilde\epsilon'_{j-1}=\tilde\epsilon'_{j+1}=+,~
 \tilde\epsilon'_i=\tilde\epsilon_i ~(i\ne j\pm1)\} 
 \]
Note here that the $+$ sign in the $(j-1)$-th place does not disappear when additional  simple reflections are applied. Hence along the columns
we have twisted coefficients.
 All that remains is to show that the arrows $\Rightarrow$ along the top row are  \lq\lq constant\rq\rq  along the columns  between cells in $Y(j)$ and cells in $Y(j+1)$ related by $s_j$. 
 This is Remark \ref{involved} (ii), (iii).
\end{Proof}

%correction

\begin{Remark}  Recall that each space $Y(j)$ corresponds to a Grasmannian: 
$$Y(j) \leftrightsquigarrow \Gr(k-1,j-1), \quad {\rm for}\quad  j=k, k+1, \cdots, n. $$ 
For the purpose of  describing
the incidence graphs it is harmless to replace $ \Gr(k-1,j-1)$ instead of  $Y(j)$.  In order  to keep track of the two possible incidence
graphs on  $ \Gr(k-1,j-1)$  (constant and twisted coefficients) we write  $  \mathcal{G}(k-1,j-1)^*$ when twisted coefficients
are involved. This way  $  \mathcal{G}(k-1,j-1)$  or $  \mathcal{G}(k-1,j-1)^*$  corresponds to $Y(j)$ at the level of incidence graphs.

 We use  this  shorthand notation to indicate the structure of the incidence graphs in Theorem \ref{mainth}:
 $  \mathcal{G}(k-1,r) \rightarrow   \mathcal{G}(k-1,r+1)  $  or  $  \mathcal{G}(k-1,r) \Rightarrow  \mathcal{G}(k-1,r+1) $ or
 $   \mathcal{G}(k-1,r)* \rightarrow  \mathcal{G}(k-1,r+1)^* $  or  $   \mathcal{G}(k-1,r)* \Rightarrow  \mathcal{G}(k-1,r+1)^* $. 
 For example the incidence graph of  $\Gr(4,9)$ with local coefficients becomes
 $$ \mathcal{G}(3,3) \Rightarrow  \mathcal{G}(3,4) \rightarrow  \mathcal{G}(3,5) \Rightarrow  \mathcal{G}(3,6)  \rightarrow  \mathcal{G}(3,7)  \Rightarrow  { \mathcal{G}(3,8)}.$$
 With constant coefficients this is indicated as
$$ \mathcal{G}(3,3)^* \rightarrow  \mathcal{G}(3,4)^*  \Rightarrow  \mathcal{G}(3,5)^* \rightarrow  \mathcal{G}(3,6)  \Rightarrow  \mathcal{G}(3,7)^*  \rightarrow  { \mathcal{G}(3,8)^*}.$$
 This can be further decomposed. For example $\Gr(3,4)^*$  can  be described as: $ \mathcal{G}(2,2) \Rightarrow  \mathcal{G}(2,3)$. Then $ \mathcal{G}(2,3)$
 corresponds to $ \mathcal{G}(1,1)^* \rightarrow  \mathcal{G}(1,2)^*$. 
\end{Remark}

With this notation, Theorem  \ref{mainth} can be restated as follows:

\begin{Theorem} \label{main2}  We have the following.

\begin{enumerate}
\item  [(i)]    The incidence graph of $\Gr(k,n)$ with constant coefficients is:
$$ \mathcal{G}(k-1,k-1)^* \rightarrow  \mathcal{G}(k-1,k)^* \Rightarrow  \mathcal{G}(k-1,k+1)^* \rightarrow \cdots   { \mathcal{G}(k-1,n-1)^*}.$$

\item  [(ii)]  The incidence graph of $\Gr(k,n)$ with  twisted  coefficients is:

$$ \mathcal{G}(k-1,k-1) \Rightarrow  \mathcal{G}(k-1,k)\rightarrow  \mathcal{G}(k-1,k+1) \Rightarrow \cdots   { \mathcal{G}(k-1,n-1)}.$$
\end{enumerate}
\end{Theorem}

This theorem leads to the well-known statement on the orientability.

\begin{Corollary} 
The real Grassmannian $\Gr(k,n)$ is orientable if and only if  $n$ is even.
\end{Corollary}

\begin{Proof}  We first show that if $n$ is even then $\Gr(k,n)$ is orientable.  We proceed by induction on $n=2j$. The  case of  $j=1$ corresponds to  $\Gr(1,2)$ which is a circle.

Using Theorem \ref{main2} we can write the incidence graph of $\Gr(k,n)$ as
$$\mathcal{G}(k-1, k-1)^* ~\to~   \mathcal{G}(k-1, k)^* ~ \Rightarrow ~\cdots  \cdots~  \mathcal{G}(k-1, n-1)^*.$$
At the same time the incidence graph of  $\Gr(k-1, n-1)$ with twisted coefficients  can be written as:
$$\mathcal{G}(k-2, k-2) ~\Rightarrow~  \mathcal{G}(k-2, k-1)^*  ~\rightarrow ~\cdots  \cdots ~ \mathcal{G}(k-2, n-2).$$
Therefore the top element  (determining orientability)  corresponds to the top  vertex  of the graph  $\mathcal{G}(k-2, n-2)$. Since $n$ is even, by induction the  top vertex  from the graph  $\mathcal{G}(k-2, n-2)$ contributes to cohomology and  $\Gr(k,n)$ is orientable.

 In the same way, one can easily see that if $n$ is odd then $\Gr(k,n)$ is not orientable. Note that we have  $\Gr(1,3)\cong\Gr(2,3)\cong{\mathbb R}P^2$.
\end{Proof}

\begin{Remark}
The homology group of the Grassmannian $\Gr(k,n)$ can be found by  Poincar\'e duality, if the variety is
orientable  ($n=$ even), i.e.
\[
H_{j}(\Gr(k,n),\mathbb{Z})=H^{k(n-k)-j}(\Gr(k,n),\mathbb{Z}) \qquad {\rm if}\quad n={\rm ~even}.
\]

In the non-orientable case ($n=$ odd), the homology group can be obtained using the incidence graph of
$\Gr(k,n)$ with twisted coefficients $\mathcal{L}^*$, that is, the graph $\G(k,n)^*$.  This is
a consequence of a generalization of Poincar\'e duality (IX.4 and VI.3 in \cite{iversen}, and
it is sometimes called the Poincar\'e-Lefschetz duality)  and that the so-called orientation sheaf  (V.3 in \cite {iversen})
${\mathcal Or}_X$ with $X=\Gr(k,n)$, corresponds, in the non-orientable cases,  to the  the twisted coefficients $\mathcal{L}^*.$

For example, consider $\Gr(2,5)$.
The graph $\G(2,5)^*$ is given by
\[
\begin{matrix}
(1,2) & \Rightarrow & (1,3) & \rightarrow & (1,4) & \Rightarrow & (1,5)\\
         &                &  \downarrow &        & \downarrow &      & \downarrow \\
         &                  &(2,3)           &\rightarrow  & (2,4) &\Rightarrow & (2,5) \\
         &                 &                      &                     &\Downarrow &          & \Downarrow\\ 
         &                  &                     &                      &   (3,4) & \Rightarrow & (3,5) \\
         &                 &                       &                     &             &                      & \downarrow \\
         &               &                        &                     &              &                      &  (4,5)
\end{matrix}
\]
which is obtained by taking the action of $h(-)$ instead of $h(+)$ for the constant coefficients.
The cohomology obtained from this graph is 
\[\begin{array}{llllll}
H^0(\Gr(2,5),\mathcal{L}^*)=0, \hskip 1.5cm    & H^4(\Gr(2,5),\mathcal{L}^*)=\mathbb{Z}_2, \\[1.0ex]
H^1(\Gr(2,5),\mathcal{L}^*)=\mathbb{Z}_2,                                      & H^5(\Gr(2,5),\mathcal{L}^*)=\mathbb{Z}_2, \\[1.0ex]
H^2(\Gr(2,5),\mathcal{L}^*)=\mathbb{Z},                        &H^6(\Gr(2,5),\mathcal{L}^*)=\mathbb{Z}.   \\[1.0ex]
H^3(\Gr(2,5),\mathcal{L}^*)=\mathbb{Z}_2,                        & {}                                         
\end{array}
\]
Then the homology group $H_*(\Gr(2,5),\mathbb{Z})$ can be obtained as $
H_j(\Gr(2,5),\mathbb{Z})=H^{6-j}(\Gr(2,5),\mathcal{L}^*)$.  In general, we have
\[
H_j(\Gr(k,n),\mathbb{Z})=H^{k(n-k)-j}(\Gr(k,n),\mathcal{L}^*)\qquad {\rm if}\quad n={\rm ~odd}.
\]
The homology groups of $\Gr(k,n)$ have been computed explicitly for lower dimensional examples in
\cite{J:79} (see Table IV in this reference).

\end{Remark}

%%%%%

\section{The Poincar\'e  polynomials} \label{recursive} 

  From the incidence graphs $\G(k,n)$ (trivial coefficients) and $\G(k,n)^*$ (twisted coefficients) 
  constructed in the previous section, 
we here consider the Poincar\'e  polynomials  $P_{(k,n)}(t)$ and 
$P^*_{(k,n)}(t)$ for those graphs.
Since our formulas for the Poincar\'e polynomials will be expressed in terms of $q$-analog of the binominal
coefficients $\left[\begin{matrix}m\\ j\end{matrix}\right]_q$, let us first list some of their properties:
\begin{enumerate}\label{pp}
\item [(q1)]  Pascal's triangle formula:  $\left[\begin{matrix} m \\ j \end{matrix}\right]_{q}= \left[\begin{matrix} m-1 \\ j \end{matrix}\right]_{q}+ q^{m-j}  \left[\begin{matrix} m-1 \\ j-1 \end{matrix}\right]_{q}.$

\item [(q2)] The polynomial  $\left[\begin{matrix} m \\ j \end{matrix}\right]_{q}$ is of degree $j(m-j)$.

\item [(q3)] Poincar\'e duality:  
 $
\left[\begin{matrix} m \\ j \end{matrix}\right]_{q} =\sum_{j=0}^{j(m-j)}\,n_j q^j=\sum_{j=0}^{j(m-j)}n_jq^{j(m-j)-j}
=q^{j(m-j)}\left[\begin{matrix}m\\j\end{matrix}\right]_{q^{-1}}.$
\end{enumerate}

\subsection{The polynomials $P^*_{(k,n)}(t)$}
We first note the following Lemma of the recursion relations of the polynomials.

\begin{Lemma} \label{pascal1*}The polynomials $P^*_{(k,n)}(t)$ satisfy
\begin{enumerate}
\item [(a$^*$)]  $P^*_{(2j,2m+1)}(t)=P^*_{(2j,2m-1)}(t)+ t^{4(m-j)+2} P^*_{(2j-2, 2m-1)}(t)$,
\item [(b$^*$)]  $P^*_{(2j+1,2m)}(t)=P^*_{(2j+1,2m-2)}(t)+ t^{4(m-j)-2} P^*_{(2j-1, 2m-2)}(t)$,
\item [(c$^*$)] $P^*_{(2j, 2m)}(t)=P^*_{(2j,2m-1)}(t)+ t^{2(m-j)} P_{(2j-1, 2m-1)}(t)$,
\item [(d$^*$)] $P^*_{(2j+1, 2m+1)}(t)=P^*_{(2j+1,2m)}(t)+ t^{2(m-j)} P_{(2j, 2m)}(t)$.
\end{enumerate}
\end{Lemma}
Before giving the proof of this Lemma, notice that this already  gives the following explicit formulas.
\begin{Proposition}  \label{prop1} The Poincar\'e polynomials $P^*_{(k,n)}(t)$ have the following
explicit form for $(k,n)=(2j+1,2m)$ or $(k,n)=(2j,2m+1)$:
\begin{enumerate}
\item[(i$^*$)]  $P^*_{(2j+1,2m)}(t)=0$,
\item[(ii$^*$)] $ P^*_{(2j,2m+1)}(t)=t^{2j}\left[\begin{matrix} m \\ j \end{matrix}\right]_{t^4}$.
\end{enumerate}
\end{Proposition}

\begin{Proof} We prove  (i$^*$)  by induction on  $m$.  For
$m=1$, we have the case of $\Gr(1,2)$ with local coefficients and $P^*_{(1,2)}(t)=0$.
Then by induction using (a) in Lemma \ref{pascal1*}, we obtain (i$^*$). 
 
To prove (ii$^*$), we also use induction on  $m$.  Again the case of $m=1$ is simple. We first note that
the relation (a) in Lemma \ref{pascal1*} gives
\begin{align*}
 P^*_{(2j,2m+1)}&=P^*_{(2j,2m-1)}(t)+ t^{4(m-j)+2} P^*_{(2j-2, 2m-1)}(t)\\
 &=
 t^{2j}\left[\begin{matrix} m-1 \\ j \end{matrix}\right]_{t^4}+t^{4(m-j)+2}t^{2(j-1)}\left[\begin{matrix} m-1 \\ j-1 \end{matrix}\right]_{t^4}.
 \end{align*}
We then use the property (q1) above to obtain
 $  t^{2j}  (\left[\begin{matrix} m-1 \\ j \end{matrix}\right]_{t^4}+t^{4(m-j)}\left[\begin{matrix} m-1 \\ j-1 \end{matrix}\right]_{t^4} )=  t^{2j}   \left[\begin{matrix} m \\ j \end{matrix}\right]_{t^4} $.
\end{Proof}

We now prove Lemma \ref{pascal1*}:

\begin{Proof} 
\begin{enumerate} 

\item [(a$^*$)] 
The incidence graph  $ \mathcal{G}(2j,2m+1)^*$  (twisted coefficients)  has the following description
corresponding to the decomposition in terms of the spaces $Y(j)$ (Theorem \ref{main2}).
\begin{align*}
&\{  \mathcal{G}(2j-1,2j-1) \Rightarrow   \mathcal{G}(2j-1,2j) \rightarrow  \cdots  \rightarrow   \mathcal{G}(2j-1,2m-3) \Rightarrow  \mathcal{G}(2j-1,2m-2) \}  \\
& \rightarrow ~  \mathcal{G}(2j-1,2m-1)~  \Rightarrow ~  { \mathcal{G}(2j-1,2m)}.
\end{align*}
The portion inside of $\{ \cdots \}$  corresponds to  $ \mathcal{G}(2j,2m-1)^*.$  The lowest degree of terms
in the last graph denoted by  ${ \mathcal{G}(2j-1,2m)}$  is   $2m-(2j-1)=2(m-j)+1.$
The incidence graph $ \mathcal{G}(2j-1,2m)$   is in turn given by
 \begin{align*}
 &\{ \mathcal{G}(2j-2,2j-2)^* \rightarrow   \mathcal{G}(2j-2,2j-1)^* \Rightarrow \cdots 
 \rightarrow  \mathcal{G}(2j-2,2m-3)^*  \Rightarrow  \mathcal{G}(2j-2,2m-2) \}^* \\
&  \rightarrow ~  \mathcal{G}(2j-2,2m-1)^*.
 \end{align*}
 Here the part  $\{ \cdots  \}$ is  $  \mathcal{G}(2j-1,2m-1) $. The lowest degree of terms 
 associated to the graph indicated by  $ \mathcal{G}(2j-2,2m-1)^*$ is 
 $2m-1-(2j-2)=2(m-j)+1$. 
 Altogether, taking into account the lowest degree  in ${ \mathcal{G}(2j-1,2m)}$,  there is a degree shift of:
 $2(m-j)+1 + 2(m-j)+1=4(m-j)+2 $.
 Thus the incidence graph  can be represented as follows:
 \[
 \begin{matrix}
  \mathcal{G}(2j,2m-1)^*& \rightarrow  &  \mathcal{G}(2j-1,2m-1)&{} &{} \\
 {}  &    {} &  {\Downarrow}&{}&{} \\
   {}  &    {} &  { \mathcal{G}(2j-1,2m-1)}&{\rightarrow} &{ \mathcal{G}(2j-2,2m-1)^*} \\
     \end{matrix}
\]     
The vertical $\Downarrow$ causes the  cancellation of all cohomology associated to 
 $\Gr(2j-1,2m-1)$.  Also the horizontal $\rightarrow$ correspond to multiplication by zero. We
 obtain   $P^*_{(2j,2m+1)}(t)=P^*_{(2j,2m-1)}(t)+ t^{N} P^*_{(2j-2, 2m-1)}(t)$ and
 $N$  is  the dimension shift  which was already computed and is given by  $4(m-j)+2 $.
 
 %%second case
 
 \item [(b$^*$)]  The incidence graph  $ \mathcal{G}(2j+1,2m)^*$  (twisted coefficients)  has the following description
(Theorem \ref{main2}).
\begin{align*}
&\{  \mathcal{G}(2j,2j) \Rightarrow   \mathcal{G}(2j,2j+1) \rightarrow  \cdots   \rightarrow   \mathcal{G}(2j,2m-4) \Rightarrow  \mathcal{G}(2j,2m-3) \}  \\
&  \rightarrow ~  \mathcal{G}(2j,2m-2) ~ \Rightarrow~   { \mathcal{G}(2j,2m-1)}.
\end{align*}  
The part $\{ \cdots \}$ corresponds to $ \mathcal{G}(2j+1, 2m-2)^*.$ The  lowest degree of the terms in the chain  complex corresponding to   ${ \mathcal{G}(2j,2m-1)}$ is 
$2m-1-2j= 2(m-j)-1$.
The incidence graph $ \mathcal{G}(2j,2m-1)$   is in turn given by
\begin{align*}
&\{ \mathcal{G}(2j-1,2j-1)^* \Rightarrow   \mathcal{G}(2j-1,2j)^* \rightarrow \cdots 
 \Rightarrow  \mathcal{G}(2j-1,2m-4)^*   \rightarrow  \mathcal{G}(2j-1,2m-3)^* \} \\
 & \Rightarrow ~  \mathcal{G}(2j-1,2m-2)^*.
 \end{align*}
 The part $\{ \cdots \}$ corresponds to $ \mathcal{G}(2j, 2m-1)$.  The lowest degree of the terms in  $ \mathcal{G}(2j-1,2m-2)^*$ is $2m-2-2j+1=2(m-j)-1$.  Altogether we have a degree shift of  $ 2(m-j)-1+ 2(m-j)-1=4(m-j)-2$.

 We thus have a graph that can
 be indicated by the diagram:
\[
 \begin{matrix}
  \mathcal{G}(2j+1,2m-2)^*& \rightarrow  &  \mathcal{G}(2j,2m-1)&{} &{} \\
 
 {}  &    {} &  {\Downarrow}&{}&{} \\
 
  {}  &    {} &  { \mathcal{G}(2j,2m-1)}&{\rightarrow} &{ \mathcal{G}(2j-1,2m-2)^*} \\
 \end{matrix} .
 \]
 This gives the  recursive formula.  The value of $r$ is  $4(m-j)-2$.
  
  \item [(c$^*$)]  The incidence graph  $ \mathcal{G}(2j,2m)^*$    has the following description:
\begin{align*}
&\{  \mathcal{G}(2j-1,2j-1) \Rightarrow   \mathcal{G}(2j-1,2j) \rightarrow  \cdots \rightarrow  \mathcal{G}(2j-1,2m-3)   \Rightarrow   \mathcal{G}(2j-1,2m-2)  \}  \\
& \rightarrow ~  { \mathcal{G}(2j-1,2m-1)}.
\end{align*}
This then can be represented as:
\[
 \mathcal{G}(2j,2m-1)^* \rightarrow   { \mathcal{G}(2j-1,2m-1)}.
 \]
This gives rise to  (c$^*$).

\item [(d$^*$)]  The incidence graph  $ \mathcal{G}(2j+1,2m+1)^*$    has the following description:
\begin{align*}
&\{  \mathcal{G}(2j,2j) \Rightarrow   \mathcal{G}(2j,2j+1) \rightarrow  \cdots  \rightarrow  \mathcal{G}(2j,2m-2)    \Rightarrow   \mathcal{G}(2j,2m-1)\}  \\
& \rightarrow  ~ { \mathcal{G}(2j,2m)}.
\end{align*}
 This then can be represented as
\[
\mathcal{G}(2j+1,2m)^* \rightarrow   { \mathcal{G}(2j,2m)}.
\]
 From here (d$^*$) follows.
  \end{enumerate}
  This completes the proof.
\end{Proof}

Similarly, we can derive the formulas for  $P^*_{(2j, 2m)}$ and  $P^*_{(2+1j, 2m+1)}(t)$.

\begin{Proposition}  \label{prop2} We have
\begin{enumerate}

%\item  $P^*_{(2j+1,2m)}=0$ 
%\item $ P^*_{(2j,2m+1)}=t^{2j}\left[\begin{matrix} m \\ j \end{matrix}\right]_{t^4}$

\item[(iii$^*$)]  $P^*_{(2j+1, 2m+1)}(t)=t^{2(m-j)}\left[\begin{matrix} m \\ j \end{matrix}\right]_{t^4}$,

\item[(iv$^*$)]   $P^*_{(2j, 2m)}(t)=(t^{2j} + t^{2(m-j)} ) \left[\begin{matrix} m \\ j \end{matrix}\right]_{t^4}$.
\end{enumerate}

\end{Proposition}

\subsection{The polynomials $P_{(k,n)}(t)$}
Let us note the following Lemma for the additional recursion relations of the Poincar\'e polynomials.
\begin{Lemma} \label{pascal1} 
The polynomials $P_{(k,n)}(t)$  satisfy
\begin{enumerate}
\item [(a)]   $P_{(2j,2m+1)}(t)=P_{(2j,2m)}(t)+ t^{2(m-j)+1} P^*_{(2j-1, 2m)}(t)$,
\item [(b)]   $P_{(2j+1,2m)}(t)=P_{(2j+1,2m-1)}(t)+ t^{2(m-j)-1} P^*_{(2j, 2m-1)}(t)$,
\item [(c)]  $P_{(2j,2m )}(t)=P_{(2j,2(m-1))}(t)+ t^{4(m-j)} P_{(2j-2, 2(m-1))}(t)$,
\item [(d)]  $P_{(2j+1,2m+1)}(t)=P_{(2j+1,2m-1)}(t)+ t^{4(m-j)} P_{(2j-1, 2m-1)}(t)$.
\end{enumerate}
\end{Lemma}

Before going through the proof,  notice that this gives rise to explicit formulas for the polynomials $P_{(k,n)}(t)$.

\begin{Theorem}\label{PoincarePoly} We have the explicit form of the Poincar\'e polynomials $P_{(k,n)}(t)$,
 \begin{enumerate}
\item[(i)] $P_{(2j,2m)}(t)=P_{(2j,2m+1)}(t)=P_{(2j+1,2m+1)}(t)=\left[\begin{matrix} m \\ j \end{matrix}\right]_{t^4}$

\item[(ii)]  $P_{(2j+1,2m)}(t)=(1+t^{2m-1})\left[\begin{matrix} m-1 \\ j \end{matrix}\right]_{t^4}$

\end{enumerate}

\end{Theorem}

\begin{Proof} We use (c) of Lemma \ref{pascal1} and induction on $m$. We then 
can rewrite this recursive formula as
\[
\left[\begin{matrix} m-1 \\ j \end{matrix}\right]_{t^4} +t^{4(m-j)} \left[\begin{matrix} m-1 \\ j-1 \end{matrix}\right]_{t^4}.
\]
By the property (q1) above,  this
is just  $\left[\begin{matrix} m \\ j \end{matrix}\right]_{t^4}$.
Using  (a) of Lemma \ref{pascal1} and  (i$^*$) of Proposition \ref {prop1} we have  $P_{(2j,2m+1)}(t)=P_{(2j,2m)}(t)+ t^{2(m-j)+1} P^*_{(2j-1, 2m)}(t)=P_{(2j,2m)}(t)$.
Using (d) of Lemma \ref{pascal1} and induction on $m$, we have   $P_{(2j+1,2m+1)}(t)=P_{(2j+1,2m-1)}(t)+ t^{4(m-j)} P_{(2j, 2m-1)}(t)$. This gives
\[
\left[\begin{matrix} m-1 \\ j \end{matrix}\right]_{t^4} +t^{4(m-j)} \left[\begin{matrix} m-1 \\ j-1 \end{matrix}\right]_{t^4},
\]
 which again corresponds to  $\left[\begin{matrix} m \\ j \end{matrix}\right]_{t^4}$.

Finally we use (b) of Lemma \ref{pascal1} and (ii$^*$) of Proposition \ref {prop1} to obtain 
\[
\left[\begin{matrix} m-1 \\ j \end{matrix}\right]_{t^4} + t^{2m-1}\left[\begin{matrix} m-1 \\ j \end{matrix}\right]_{t^4},
\]
which gives (2).
\end{Proof}

We now prove Lemma \ref{pascal1}.
\begin{Proof}
\begin{enumerate}

\item [(a)]  We have a graph of the form,
\begin{align*}
&\{  \mathcal{G}(2j-1,2j-1)^* \rightarrow   \mathcal{G}(2j-1,2j)^* \Rightarrow  \cdots   \rightarrow  \mathcal{G}(2j-1,2m-2)^*  \Rightarrow   \mathcal{G}(2j-1,2m-1)^*  \} \\
& \rightarrow  ~ { \mathcal{G}(2j-1,2m)^*}.
\end{align*}
The part $\{ \cdots \}$  corresponds to  $ \mathcal{G}(2j,2m).$   
The lowest degree of the terms in the last graph denoted by  ${ \mathcal{G}(2j-1,2m-1)}$  is   $2m-2j+1=2(m-j)+1.$ We then have a graph which can be summarized as follows:
\[
\begin{matrix}
 \mathcal{G}(2j,2m)& \rightarrow  &  \mathcal{G}(2j-1,2m)^*
 \end{matrix} .
 \]
  From here  $P_{(2j,2m+1)}=P_{(2j,2m)}(t)+ t^r P^*_{(2j-1, 2m)}(t)$ follows and $r=2(m-j)+1$.

\item [(b)]  We have:
\begin{align*}
& \{  \mathcal{G}(2j,2j)^* \rightarrow   \mathcal{G}(2j,2j)^* \Rightarrow  \cdots \rightarrow  \mathcal{G}(2j,2m-3)^*   \Rightarrow   \mathcal{G}(2j,2m-2)^*  \} \\
& \rightarrow  ~ { \mathcal{G}(2j,2m-1)^*}.
\end{align*}
The portion inside $\{ \cdots \}$  corresponds to  $ \mathcal{G}(2j+1,2m-1)^*.$
We then end up with a diagram that can be summarized as: 
\[ \begin{matrix}
 \mathcal{G}(2j+1,2m-1)& \rightarrow  &  \mathcal{G}(2j,2m-1)^*
  \end{matrix} .
\]  
The lowest degree in the chain complex corresponding to  $ \mathcal{G}(2j,2m-1)^*$ is then: $2(m-j)-1$ 

%plus $2(m-j)$. Hence  $4(m-j)-1$.
%%%

\item [(c)]   The incidence graph  $ \mathcal{G}(2j,2m)$   has the following description
(Theorem \ref{main2}).
\begin{align*}
&\{  \mathcal{G}(2j-1,2j-1)^* \rightarrow   \mathcal{G}(2j-1,2j)^* \Rightarrow  \cdots  \rightarrow   \mathcal{G}(2j-1,2m-4)^* \Rightarrow  \mathcal{G}(2j-1,2m-3)^* \}  \\
&  \rightarrow  ~ \mathcal{G}(2j-1,2m-2)^* ~ \Rightarrow ~  { \mathcal{G}(2j-1,2m-1)^*}.
\end{align*}
The portion inside of $\{ \cdots \}$  corresponds to  $ \mathcal{G}(2j,2m-2).$  
The lowest degree of the terms in the last graph denoted by  ${ \mathcal{G}(2j-1,2m-1)}$  is   $2m-2j=2(m-j).$   The incidence graph $ \mathcal{G}(2j-1,2m-1)^*$   is in turn given by
\begin{align*}
&\{ \mathcal{G}(2j-2,2j-2) \Rightarrow   \mathcal{G}(2j-2,2j-1) \rightarrow \cdots
 \rightarrow  \mathcal{G}(2j-2,2m-4) \Rightarrow  \mathcal{G}(2j-2,2m-3) \} \\
&  \rightarrow   \mathcal{G}(2j-2,2m-2).
\end{align*}
Here the part $\{ \cdots  \}$ is now  $  \mathcal{G}(2j-1,2m-2)^* $. 
 The lowest degree of terms 
 associated to the graph indicated by  $ \mathcal{G}(2j-2,2m-2)^*$ is 
 $2m-2-(2j-2)=2(m-j)$. 
 Altogether, taking into account the lowest degree  in ${ \mathcal{G}(2j-2,2m-2)}$,  there is a degree shift of:
$4(m-j)$.
We thus have a graph that can
 be indicated by the diagram:
 \[
  \begin{matrix}
 \mathcal{G}(2j,2m-2)& \rightarrow  &  \mathcal{G}(2j-1,2m-2)^*&{} &{} \\
 
 {}  &    {} &  {\Downarrow}&{}&{} \\
 
  {}  &    {} &  { \mathcal{G}(2j-1,2m-2)^*}&{\rightarrow} &{ \mathcal{G}(2j-2,2m-2)}  \\
 \end{matrix} 
 \]
  From here we obtain the recursive formula.

\item [(d)]  We start with the corresponding  graph that can be represented by a diagram:
\[
 \begin{matrix}
 \mathcal{G}(2j+1,2m-1)& \rightarrow  &  \mathcal{G}(2j,2m-1)^*&{} &{} \\
 
 {}  &    {} &  {\Downarrow}&{}&{} \\
 
  {}  &    {} &  { \mathcal{G}(2j,2m-1)^*}&{\rightarrow} &{ \mathcal{G}(2j-1,2m-1)} \\
 \end{matrix} 
 \]
The lowest degree in the piece corresponding to  ${ \mathcal{G}(2j-1,2m-1)} $ is $4(m-j)$. This becomes the
formula $P_{(2j+1,2m+1)}(t)=P_{(2j+1,2m-1)}(t)+ t^{4(m-j)} P_{(2j, 2m-1)}(t).$
\end{enumerate}
This completes the proof of Lemma \ref{pascal1}.
\end{Proof}
\begin{Remark}
We note the similarity in the formulae of the Poincar\'e polynomials for the real and complex Grassmannians.
That is, in the case of (i) in Theorem \ref{PoincarePoly}, we have
\[
P_{(k,n)}(t)=P^{\Complex}_{(\lfloor{k}/{2}\rfloor,\lfloor{n}/{2}\rfloor)}(t^2).
\]
Also, in particular, if we take the limit $n\to\infty$ for the case $|t|<1$, we have for both cases (i) and (ii)
\[
P_{(k,\infty)}(t)=P^{\Complex}_{(\lfloor{k}/{2}\rfloor,\infty)}=\prod_{j=1}^{\lfloor k/2\rfloor}\frac{1}{1-t^{4j}}.
\]
This is a consequence of the structure of the cohomology ring of the Grassmannians
in terms of the characteristic classes (see for example \cite{BT:82, MS:74}):
It is well known that the real cohomology ring of the complex Grassmannian $\Gr(k,n,\Complex)$ can be described by
\[
H^*(\Gr(k,n,\Complex),\Real)\cong \frac{\Real[c_1,\ldots,c_{n-k},\bar{c}_1,\ldots,\bar{c}_k]}{\{c\cdot  \bar{c}=1\}},
\]
where  $c:=1+c_1+\cdots +c_{n-k}$ and $\bar{c}_i=1+\bar{c}_1+\cdots+\bar{c}_k$ with the Chern classes $c_j$'s  defined by
\[
c_j\in H^{2j}(\Gr(k,n,\Complex),\Real).
\]
The Poincar\'e polynomial $P_{(k,n)}^{\Complex}(t)$ is then given by
\[
P_{(k,n)}^{\Complex}(t)=\left[\begin{matrix}n\\k\end{matrix}\right]_{t^2}=\prod_{j=0}^{k-1}\frac{1-t^{2(n-j)}}{1-t^{2(k-j)}}.
\]
In particular, for the classifying space $BU(k)$ as the infinite Grassmannian $\Gr(k,n,\Complex)$ with $n=\infty$, the cohomology ring is given by
\[
H^*(BU(k),\Real)\cong \Real[c_1,c_2,\ldots, c_k], \qquad {\rm with}\quad c_j\in H^{2j}(BU(k),\Real).
\]
and the Poincar\'e polynomial becomes the series given by
\[
P^{\Complex}_{(k,\infty)}(t)=\prod_{j=1}^{k}\frac{1}{1-t^{2j}}.
\]

For the classifying space $BO(k)$ as the infinite Grassmannian $\Gr(k,\infty)$, the cohomology ring is
known to be
\[
H^*(BO(k),\Real)\cong \Real[p_1,p_2,\ldots,p_{\lfloor{k}/{2}\rfloor}],
\]
where the generators of the ring are given by the Potrjagin classes $p_j\in H^{4j}(\Gr(k,n),\Real)$.
The Poincar\'e series of $BO(k)$ is
\[
P_{(k,\infty)}(t)=\prod_{j=1}^{\lfloor k/2\rfloor}\frac{1}{1-t^{4j}}=\lim_{n\to\infty}P_{(k,n)}(t),
\]
where the limit of course make sense for $|t|<1$.   
\end{Remark}

%%%%%%%%%%%%%%%%%%%%%%%%%%%%%%%%%%%%%%%%%%%%%%%%%%%%%%%%%%%%

\section{The $\mathbb{F}_q$ points on $\Gr(k,n)$ and the Poincar\'e polynomials}\label{Ppolynomial}

Here we first introduce the {\it weighted} Schubert cell where the weight is given by the $q^{\eta(w_{\sigma})}$ (and $q^{\eta(w_{\sigma})^*}$) defined in subsection \ref{connection}.
Then we define a polynomial $p_{(k,n)}(q)$ (and $p^*_{(k,n)}(q)$) from the incidence graph
$\G(k,n)$ (and $\G(k,n)^*$) based on the weights of the Schubert cells.
It turns out that the polynomial $p_{(k,n)}(q)$ is related to the number of points on
$\Gr(k,n)$ over the finite field $\mathbb{F}_q$, and also related to the Poincar\'e polynomials
$P_{(k,n)}(t)$ found in the previous section. The point here is that the notion of the weighted Schubert
cells gives a simple method to compute the Poincar\'e polynomials and the $\mathbb{F}_q$ points on
$\Gr(k,n)$.

\begin{Remark}
As in the case of the Toda-flow for the real flag variety discussed in \cite{casian:06}, those polynomials
$p_{(k,n)}(q)$ can be computed by counting the number of blow-ups along the KP flow.
\end{Remark}

\subsection{The $q$-weighted Schubert cells}   \label {power} 
Let us first recall that  the Schubert cell
$(\sigma_1,\ldots,\sigma_k)$ can be identified as an element   $w_{\sigma}$ in
${\mathcal{S}}^{(k)}_n$. Hence  the vertices of the incidence graph 
$\G(k,n)$ correspond to certain minimal length Weyl group representatives $w_{\sigma}$.  As shown in
subsection \ref{connection}, given  ${\mathcal L}$  local system on $\Gr(k,n)$ 
certain  powers of $q$ can be associated   to each  vertex $w_{\sigma}$.  Since there are only two local systems two consider, we can simplify the notation of 
\cite{casian:06}, and  just denote by $q^{\eta (w_{\sigma})}$ the power of $q$ assigned in the  case when  ${\mathcal L}$ is constant and $q^{\eta (w_{\sigma})^*}$ the power of $q$ in the  case
when ${\mathcal L}$ is twisted.  We then associate powers $q^{\eta (w_{\sigma})}$ to each vertex $w_{\sigma}$ of $\mathcal{G}(k,n)$ and  powers $q^{\eta (w_{\sigma})^*}$ to each vertex $w_{\sigma}$ of  $\mathcal{G}{(k,n)}^*$. We call the Schubert cells $X_{w_{\sigma}}$ with those powers $q^{\eta(w)}$ 
the {\it weighted} Schubert cells, denoted by $(w_{\sigma},q^{\eta(w_{\sigma})})$ for each $w_{\sigma}\in\S_n^{(k)}$.

Here we consider only the $\eta(w)$ (the $\eta(w)^*$ can be treated in a similar way).
Let us first define the following set of weighted vectors.
\begin{itemize}
\item[(a)] For odd $k$, 
\[
\left\{ e_j(q)=(E_j; q^{\lfloor\frac{j}{2}\rfloor}): j=1,\ldots,n\right\}\,.
\]
(i.e. the weights are assigned as $(1,q,q,q^2,q^2,\ldots,q^{\lfloor\frac{n}{2}\rfloor})$),
\item[(b)]  For even $k$,
\[
\left\{e_j(q)=(E_j; q^{\lfloor\frac{j-1}{2}\rfloor}): j=1,\ldots, n\right\}\,.
\]
(i.e. the weights are assigned as $(1,1,q,q,\ldots, q^{\lfloor\frac{n-1}{2}\rfloor})$).
\end{itemize}
This is a $q$-deformation of the signed vector $e_j$, and
the sign for each vector $E_j$ is given by setting $q=-1$.  Then we can find 
 the explicit form of the $\eta(w_{\sigma})$.
\begin{Lemma} The function $\eta(w_{\sigma})$ with the representation
$w_{\sigma}=(\sigma_1,\ldots,\sigma_k)$ is given by
\begin{align*}
\eta(w_{\sigma})&=\sum_{j=1}^k\lfloor\frac{\sigma_j-j}{2}\rfloor \quad \quad\quad{\rm if}\quad k={\rm odd},\\
\eta(w_{\sigma})&=\sum_{j=1}^k\lfloor\frac{\sigma_j-j-1}{2}\rfloor \quad\quad{\rm if}\quad k={\rm even}.
\end{align*}
\end{Lemma}
To show this, we note that the $\eta(w_{\sigma})$ satisfies the following conditions which
uniquely determine $\eta(w_{\sigma})$ for given Schubert cell $w_{\sigma}=(\sigma_1,\ldots,\sigma_k)$.
\begin{enumerate}
\item [(i)]  To the top cell $e=(1,2,\ldots,k)$ we associate $\eta (e)=0$.
\item [(ii)]  If two Schubert cells $w$ and $w'$ are joined by $\Rightarrow$, i.e.
 $w \Rightarrow w^\prime$, then $\eta (w)=\eta (w^\prime)$.
\item [(iii)]   If two Schubert cells $w$ and $w'$ are  joined by $\rightarrow$ (not $\Rightarrow$), then $\eta (w^\prime)= \eta (w)+1$.
\end{enumerate}
We then define a  polynomial $p_{(k,n)}(q)$ as an  alternating sum,
\[
p_{(k,n)}(q):=\sum_{w\in\S_n^{(k)}} (-1)^{\ell (w) } q^{ \eta(w)} 
\]

We now offer a direct construction of the polynomials  $p_{(k,n)}(q)$ which leads to their
direct calculation.  An alternative way to proceed is through arguments  similar to those
leading to the recursive formulas for Poincar\'e polynomials.
We have the following Theorem:
\begin{Theorem} \label{polynomials}
The polynomials $p_{(k,n)}(q)$ take the forms,
\begin{align*}
&p_{(2j,2m)}(q)=p_{(2j,2m+1)}(q)=p_{(2j+1,2m+1)}(q)=\left[\begin{matrix}m\\j\end{matrix}\right]_{q^2}\\
&p_{(2j+1,2m)}(q)=(1-q^m)\left[\begin{matrix}m-1\\j\end{matrix}\right]_{q^2}\,.
\end{align*}
\end{Theorem}
\begin{Proof}
Let us first consider the case $\Gr(2j,2m)$. In this case we have
\[
\left\{(E_1;1),(E_2;1),(E_3,q),(E_4;q),\ldots,(E_{2m-1};q^{m-1}),(E_{2m};q^{m-1})\right\}\,.
\]
We compute $p_{(2j,2m)}(q)$ by computing all wedge products $(\sigma_1,\cdots, {\sigma_{2j}})\leftrightarrow e_{\sigma_1}(q)\wedge\cdots\wedge e_{\sigma_{2j}}(q)$. It is immediate that if a cell 
contains only one term from a pair $\{(E_{2i-1};q^{i-1}),(E_{2i};q^{i-1})\}$, then this cell is canceled by
the cell containing the same terms except the term from the pair replaced by the other one.
Thus the cells which contribute to $p_{(2j,2m)}(q)$ are given by the wedge products
containing pairs $\{(E_{2i-1};q^{i-1}),(E_{2i},q^{i-1})\}$ for some $i$, i.e.
\[
\left(E_{2a_1-1}\wedge E_{2a_1}\wedge \cdots\wedge E_{2a_j-1}\wedge E_{2a_j}; q^{j(j-1)+4|Y(a)|}\right)\,.
\]
Here $Y(a)$ represents the Young diagram corresponding to $(a_1,\ldots,a_j)$, and $4$ implies that
each box in $Y(a)$ consists of $2\times 2$ cube 
of the boxes for the original Young diagram $Y_{\lambda}$ with
\[
\lambda=\left(2(a_1-1)+1,2(a_1-1)+2,\ldots, 2(a_j-j)+2j-1,2(a_j-j)+2j\right).
\]
 Since the generating function of
the number of cells in $\Gr(j,m)$ is given by (\ref{qbinom}), we obtain the formula $p_{(2j,2m)}(q)$.

Now let us consider the case $\Gr(2j,2m+1)$. In this case we have,
\[
\left\{(E_1;1),(E_2;1),\ldots,(E_{2m-1};q^{m-1}),(E_{2m};q^{m-1}),(E_{2m+1};q^m)\right\}\,.
\]
It is easy to see that if a cell contains $E_{2m+1}$ in $2j$-wedge product, then this cell has
no contribution in $p_{(2j,2m+1)}(q)$.
Then the situation is the same as the case $\Gr(2j,2m+1)$. A similar argument can be applied for 
the case of $\Gr(2k+1,2n+1)$ (in this case, all the terms containing $E_1$ vanish).

Finally we consider the case $\Gr(2j+1,2m)$. We start with
\[
\left\{(E_1;1),(E_2;q),(E_3;q),\ldots,(E_{2m-2};q^{m-1}),(E_{2m-1};q^{m-1}),(E_{2m};q^m)\right\}\,.
\]
We note that the cells containing both $E_1$ and $E_{2m}$ have no contribution.
There are two types of cells which contribute; (i) those consisting of $E_1$ and pairs
$\{E_{2i},E_{2i+1}\}$ of the same degree $q^i$, and (ii) those consisting of $E_{2n}$ and
pairs of the same degrees. The first case (i), i.e. without $E_{2m}$, gives the same polynomial
as in the case $\Gr(2j+1,2m-1)$. The second case (ii), i.e. without $E_1$, we have the same 
polynomial times $-q^m$ due to the degree of $E_{2m}$ and $(-1)^{l(w)}=(-1)^{2m-1}$.
We thus obtain the result for $p_{(2j+1,2m)}(q)$.
\end{Proof}
\begin{Remark}
The polynomials $p_{(k,n)}(q)$ in Theorem \ref{polynomials} are related to the Poincar\'e polynomials $P_{(k,n)}(t)$ found in Theorem \ref{PoincarePoly}.
\begin{itemize}
\item[(a)] For $\Gr(k,n)$ with $(k,n)=(2j,2m)$ or $(2j,2m+1)$ or $(2j+1,2m+1)$, the Poincar\'e
polynomial is given by $P_{(k,n)}(t)=p_{(k,n)}(q=t^2)$, i.e.
\[
P_{(k,n)}(t)=\left[\begin{matrix}m\\ j\end{matrix}\right]_{t^4}=\prod_{i=0}^{j-1}\frac{1-t^{4(m-i)}}{1-t^{4(j-i)}}.
\]
For example, $P_{(2,4)}(t)=P_{(2,5)}(t)=P_{(3,5)}(t)=1+t^4$.
\item[(b)] For $\Gr(k,n)$ with $(k,n)=(2j+1,2m)$, $P_{(k,n)}(t)$ is given by
\[
P_{(k,n)}(t)=(1+t^{2m-1})\,\left[\begin{matrix} m-1\\ j\end{matrix}\right]_{t^4}.
\]
Namely, the factor $(1-q^m)$ in the polynomial $p_{(k,n)}(q)$ is replaced by $(1+t^{2m-1})$ and $\left[\begin{matrix}{m-1}\\{j}\end{matrix}\right]_{q^2}$
with $q=t^2$.
For example, $p_{(3,6)}(q)=(1-q^3)\left[\begin{matrix}2\\1\end{matrix}\right]_{q^2}$ and $P_{(3,6)}(t)=(1+t^5)(1+t^4)=1+t^4+t^5+t^9$.
\end{itemize}
\end{Remark}

\begin{Example}  Following the constructive arguments in the proof of Theorem \ref{polynomials}, 
we directly find the polynomials $p_{(k,n)}(q)$ for the cases $\Gr(4,8)$ and $\Gr(5,12)$:
\begin{itemize}
\item[(a)] For $\Gr(4,8)$, we have the following cells which contribute the polynomial $p_{(4,8)}(q)$,
\[\begin{array}{lllll}
(1,2,3,4)=e,\\
(1,2,5,6)=s_4s_5s_3s_4,\\
(1,2,7,8)=s_6s_7s_5s_6s_4s_5s_3s_4,   \\
(3,4,5,6)=s_2s_3s_1s_2s_4s_5s_3s_4,    \\
(3,4,7,8)=s_2s_3s_1s_2s_6s_7s_5s_6s_4s_5s_3s_4, \\
(5,6,7,8)=s_4s_5s_3s_4s_2s_3s_1s_2s_6s_7s_5s_6s_4s_5s_3s_4
\end{array}
\]
Here the cells are represented by the elements of $\mathcal{S}^{(4)}_8$.
Then the polynomial is given by
\[
p_{(4,8)}(q)=1+q^2+2q^4+q^6+q^8=\left[\begin{matrix}4\\2\end{matrix}\right]_{q^2}\,.
\]
The Young diagrams of those cells are given by
\begin{align*}
&(1,2,3,4)=\emptyset,\qquad (1,2,5,6)=\young[2,2][4],\qquad (1,2,7,8)=\young[4,4][4] \\
&(3,4,5,6)=\young[2,2,2,2][4],\qquad (3,4,7,8)=\young[4,4,2,2][4],\qquad (5,6,7,8)=\young[4,4,4,4][4]
\end{align*}
Note here that the box $\young[2,2][4]$ gives the unit for those diagrams and each diagram represents the
Pontrjagin class $p_j\in H^{4j}(\Gr(4,8),\Real)$ for $j=1,2$, and this explains the relation,
\[
p_{(4,8)}(q)=|\Gr(2,4,\mathbb{F}_{q^2})|.
\]
Then the cohomology ring may be expressed by
\[
H^*(\Gr(4,8),\Real)\cong \frac{\Real[p_1,p_2,\bar p_1,\bar p_2]}{\{p\cdot \bar p=1\}},
\]
where $p=1+p_1+p_2$ and $\bar p=1+\bar p_1+\bar p_2$.
\item[(b)] For $\Gr(5,12)$, we have the contributing cells containing either $E_1$ or $E_{12}$,
\[
\begin{array}{llllllll}
(1,2,3,4,5)      \qquad       \quad&    &\quad        &   (2,3,4,5,12)    \\
(1,2,3,6,7)                       &     &                   &  (2,3,6,7,12)     \\
(1,2,3,8,9)                        &   &                  &   (2,3,8,9,12)     \\
(1,2,3,10,11)                  &   &            &  (2,3,10,11,12) \\
(1,4,5,6,7)                        &[ 1~ \to~12]&                  &   (4,5,6,7,12)   \\
(1,4,5,8,9)                      &  &                    &   (4,5,8,9,12)   \\
(1,4,5,10,11)                   &&                  &    (4,5,10,11,12)  \\
(1,6,7,8,9)                       & &                   &    (6,7,8,9,12)    \\
(1,6,7,10,11)                    & &                 &  (6,7,10,11,12)\\
(1,8,9,10,11)                       & &              &    (8,9,10,11,12)
\end{array}
\]
The cells in the left hand side gives $1+q^2+2q^4+2q^6+2q^8+q^{10}+q^{12}=(1-q^{10})(1-q^8)/(1-q^4)(1-q^2)$, and those in the right hand side gives $-q^6$ times the same polynomial, i.e.
\[
p_{(5,12)}(q)=(1-q^6)\left[\begin{matrix}5\\2\end{matrix}\right]_{q^2}\,.
\]
The Young diagrams for the cells in the left column are given by
\begin{align*}
&\emptyset\qquad\young[2,2][4]\qquad \young[4,4][4]\qquad\young[6,6][4]\qquad\young[2,2,2,2][4]
\qquad\young[4,4,2,2][4]\qquad\young[4,4,4,4][4]\qquad\young[6,6,4,4][4]\qquad\young[6,6,6,6][4]
\end{align*}
The Young diagrams of the cells in the right column are given by the above ones combined with the first cell
in the right column,
\[
(2,3,4,5,12)=\young[7,1,1,1,1][4].
\]
That is, we have, from the second cell, 
\begin{align*}
\qquad\young[7,3,3,1,1][4]\quad \young[7,5,5,1,1][4]\quad\young[7,7,7,1,1][4]\quad\young[7,3,3,3,3][4]
\quad\young[7,5,5,3,3][4]\quad\young[7,7,7,3,3][4]\quad\young[7,5,5,5,5][4]\quad\young[7,7,7,5,5][4]
\quad\young[7,7,7,7,7][4]
\end{align*}
Note again that those cells are expressed by the diagrams $\young[2,2][4]$ and $\young[7,1,1,1,1][4]$.
Each diagram of divisible by $\young[2,2][4]$ represents the Pontrjagin class $p_j\in H^{4j}(\Gr(5,12),\Real)$, and the hook diagram may correspond to an extra element, say $r$, of degree 11 with the property $r^2=0$.  It then may be natural to conjecture that the cohomology ring of $\Gr(5,12)$ has the structure,
\[
H^*(\Gr(5,12),\Real)\cong\frac{\Real[p_1,p_2,p_3,\bar{p}_1,\bar{p}_2,r]}{
\{p\cdot \bar{p}=1,~r^2\}},
\]
where $p=1+p_1+p_2+p_3,$ and $\bar{p}=1+\bar{p}_1+\bar{p}_2$.
\end{itemize}
\end{Example}

 The proof of the following Proposition is similar to  the calculation of  $P^*_{(2j,2m+1)}(t)$, 
 $P^*_{(2j+1,2m+1)}(t)$
and is omitted.

\begin{Proposition} \label{star} We have:

\begin{enumerate}

\item[]  $p^*_{(2j,2m+1)}(q)=q^{j}\left[\begin{matrix} m \\ j \end{matrix}\right]_{q^2}$

\item[]  $p^*_{(2j+1,2m+1)}(q)=q^{m-j}\left[\begin{matrix} m \\ j \end{matrix}\right]_{q^2}$

\end{enumerate}

\end{Proposition}

%%%%%%%%%%%%%%%%%%%%%%%%%%%%%%%%%%%%%%%%%%%%

\subsection{The number of $\mathbb{F}_q$-points on $\Gr(k,n)$}\label{Fqpoints}
The main goal of this section is to show that
the $\mathbb{F}_q$-points on the Grassmannian $\Gr(k,n)$ is given by
\[
|\Gr(k,n)_{\mathbb{F}_q}|= q^r\,|p_{(k,n)}(q)|
\]
where $r$ is given by $r=k(n-k)-{\rm deg}\,(p_{(k,n)}(q))$.

 In order to calculate 
$|\Gr(k,n)_{\mathbb{F}_q}|$, 
we first introduce the complexification $\Gr(k,n)_{\Complex}$ of the real Grassmannian $\Gr(k,n)$,
which is not the complex Grassmannian $\Gr(k,n,\mathbb{C})$, but rather, a Zariski open subset
of $\Gr(k,n,\Complex)$ having the same homotopy type as $\Gr(k,n)$. So, for instance the complexification of ${\mathbb R}P^1=\Gr(1,2)$ is not $\Gr(1, 2, \Complex)$
which is ${\mathbb C} P^1$, but rather,  the set of all lines in ${\mathbb C}^2=\{(x,y)\}$  
such that $x^2+y^2\not=0$. 
This is an open subset of   ${\mathbb C}P^1$ which has the homotopy type of a circle.

\begin{Definition}
 We fix a real vector subspace of dimension $k$,  $V\subset {\mathbb R}^n$ and its complexification
$V_{\mathbb{C}}=V \oplus_{\Real} \sqrt{-1}\, V$.  The group ${\rm SO}_n(\Real)$ acts transitively on   $\Gr(k,n)$ but the orbit ${\rm SO}_n({\mathbb C})\cdot V_{\mathbb{C}}$
of the vector space $V_{\mathbb{C}}$ does not exhaust all  of $\Gr(k,n,\mathbb C)$.  We define this orbit  as the {\it complexification} of $\Gr(k,n)$, and denoted by $\Gr(k,n)_{\mathbb C}$.   The real Grassmanian  $\Gr(k,n)$ can then be described as $K/ L\cap K$, where $K={\rm SO}_n(\Real)$  and
the Levi factor $L$ of a maximal parabolic subgroup
containing the Borel subgroup $B$  of upper triangular matrices in $\SL_n({\mathbb R})$ and 
the complexification is  $K({\mathbb C})/ L({\mathbb C}) \cap K({\mathbb C})=\SO_n({\mathbb C})\cdot V_{\mathbb{C}}$  which has the same homotopy type as  $\Gr(k,n)$. 

This naturally extends to other semi-simple Lie groups.  In brief, if $G=KAN$ (Iwasawa decomposition)
and $P=LN$ is a maximal parabolic with Levi factor $L$ containing a maximally split Cartan  subgroup $H=TA$ , then we have $G/P\cong K/K\cap L$ and the complexification will be $K({\mathbb C})/ L({\mathbb C})\cap K({\mathbb C})$.
\end{Definition}

\begin{Example} \label{rpone} Consider ${\Real}P^{n-1}$.  If we fix $V$ a one dimensional vector subspace of ${\Real}^{n}$.
Then  the $\SO(n,{\Complex})$ orbit  containing $V_{\mathbb{C}}$ can be described as the set of all lines in $\Complex^n$, denoted by ${\mathbb C}(x_1,\cdots x_{n})$,
such that $x_1^2+\cdots x_{n}^2\not=0$. This is an open subset  $\Gr(1,n)_{\mathbb C}$ of  
the Grassmanian $\Gr(1,n, {\mathbb C}) ={\mathbb C}P^{n-1}$.
\end{Example}

\begin{Example}  \label{gr2} Let us consider $\Gr(2,4)$. Using the Pl\"ucker  embedding   as in Example \ref{isotropic}, the complexification can be explicitly  described as follows.
$$\Gr(2,4)_{\Complex} = \left\{ {\mathbb C} (z_1, z_2, z_3, w_1, w_2, w_3) : 
\begin{array}{lll}
z_1^2+z_2^2+z_3^2-w_1^2-w_2^2-w_3^2=0, \\
z_1^2+z_2^2+z_3^2 \not=0
\end{array}\right\}.$$   This is an open subset of  $\Gr(2,4,{\mathbb C})$.
This formula is useful for counting the $\mathbb{F}_q$-points on $\Gr(k,n)$ (see below).
\end{Example}

We consider $\overline{\mathbb{F}}_q$ an algebraic closure of a field $\mathbb{F}_q$ with $q$ elements, and consider Grassmannian varieties $\Gr(k,n,\overline{\mathbb{F}}_q)$.
As a set of points this is the set of all the $k$-dimensional subspaces of $\overline{\mathbb{F}}^n_q$.

However the $\Gr(k,n,\overline{\mathbb{F}}_q)$ is the $\overline{\mathbb{F}}_q$-analogue of
the standard complex Grassmannian variety $\Gr(k,n,\mathbb{C})$.  Here we are interested instead in the $\overline{\mathbb{F}}_q$-counterparts of the complexifications $\Gr(k,n)_{\mathbb{C}}$ of
real Grassmannians in counting their $\mathbb{F}_q$-points.

Assume that  $q$ is a power of a prime number $p\not=2$  such  that  in ${\mathbb F}_q$ the 
 polynomial  $x^2+1$ is reducible i.e. $\sqrt{-1} \in  {\mathbb F}_q$.
 If $p$ has the form $p=4k+1$ with $k$ an integer,  this will be the case e.g. $q=5$. 
 We then consider the analogue of $\Gr(k,n)$ over $\mathbb{F}_q$, that is, the Zariski open
 subset $\Gr(k,n)_{\overline{\mathbb{F}}_q}$ of $\Gr(k,n,\overline{\mathbb{F}}_q)$ and then
 the corresponding  ${\mathbb F}_q$ points. Then we have the following results
 for $|S^n(\mathbb{F}_q)|$, the number of $\mathbb{F}_q$-points on $S^n$:
 \begin{itemize}
 \item[(a)] For $n=2m-1$, we have
 \[
 |S^{2m-1}(\mathbb{F}_q)|=q^{m-1}(q^m-1)\,.
 \]
 \item[(b)] For $n=2m$, we have
 \[
 |S^{2m}(\mathbb{F}_q)|=q^m(q^m+1)\,.
 \]
 \end{itemize}
 This can be obtained as follows: Let us first consider the case $n=1$, i.e.
 \[
 S^1(\mathbb{F}_q)=\{(x,y)\in\mathbb{F}_q^2:x^2+y^2=1\}\,.
 \]
 Then using the formulae for the stereographic projection; $x=\frac{2u}{u^2+1},
 y=\frac{u^2-1}{u^2+1}$ with $y\ne 1$ and $\{u\in \mathbb{F}_q:u^2+1\ne0\}$. 
 Since $\sqrt{-1}\in\mathbb{F}_q$,
 we have $2$ points in $\{u^2+1=0\}$. Counting the point $(0,1)$, the north pole, we have
 $|S^{1}(\mathbb{F}_q)|=q-2+1=q-1$. Now consider the case $n=2$, we have
 $x=\frac{2u_1}{u_1^2+u_2^2+1}, y=\frac{2u_2}{u_1^2+u_2^2+1}, z=\frac{u_1^2+u_2^2-1}{u_1^2+u_2^2+1}$ with $z\ne 1$ and $\{(u_1,u_2)\in \mathbb{F}_q^2: u_1^2+u_2^2+1\ne 0\}$.
 This gives $q^2-(q-1)$ points (note $(q-1)$ is the number of points in $u_1^2+u_2^2+1=0$).
 We now add the points of the north pole $(x,y,1)$ with $x^2+y^2=0$. This gives $2(q-1)+1$,
 where $2(q-1)$ for $x=\pm\sqrt{-1}y\ne 0$ and $1$ for $(0,0,1)$. Then we have
 $|S^2(\mathbb{F}_q)|=q^2-(q-1)+2(q-1)+1=q(q+1)$. Using the induction, one can show
 that the number of points in the north pole is given by
 \begin{equation}\label{Fpoints}\left\{\begin{array}{llll}
 |\{(x_1,\ldots, x_{2m-1})\in \mathbb{F}_q^{2m-1}:x_1^2+\cdots+x_{2m-1}^2=0\}|=q^{2m-2},\\[2.0ex]
 |\{(x_1,\ldots,x_{2m})\in \mathbb{F}_q^{2m}:x_1^2+\cdots+x_{2m}^2=0\}|=q^{2m-1}+q^m-q^{m-1}.
 \end{array}\right.
 \end{equation}
 Then one can obtain the above formulae for $|S^n(\mathbb{F}_q)|$.

\begin{Example} \label{example2}  We consider the analogue of ${\mathbb R}P^{n-1}$ over ${\mathbb F}_q$,  that is, 
the Zariski open subset $\Gr(1,n)_{{\mathbb F}_q}$ of $\Gr(1,n, {\mathbb F}_q)$.  Recall that we have
the Schubert decomposition,
\[
\Complex P^{n-1}=\Gr(1,n,\Complex )=\{(*,\ldots,*,1)\}\sqcup\{(*,\ldots,*,1,0)\}\sqcup\cdots\sqcup\{(1,0,\ldots,0)\}\,,
\]
and we have  $|\Gr(1,n,\mathbb{F}_q)|=[n]_q=1+q+\cdots+q^{n-1}$. To find $|\Gr(1,n)_{\mathbb{F}_q}|$, the number of $\mathbb{F}_q$ points on $\Real P^{n-1}$,   we have to remove the points in the Schubert cells given by
\[
\{(x_1,\ldots,x_k,1,0,\ldots, 0)\in \mathbb{F}_q^n: x_1^2+\cdots+x_k^2+1=0\}\qquad k=1,\ldots,n-1.
\]
Since $\sqrt{-1}\in \mathbb{F}_q$, those sets are equivalent to $S^{k-1}$ on $\mathbb{F}_q$
whose points can be counted from the formulae above. Thus we get
\begin{itemize}
\item[(a)] For $n=2m$, 
\[
|\Gr(1,2m)_{\mathbb{F}_q}|=q^{m-1}(q^m-1)\,.
\]
\item[(b)] For $n=2m+1$,
\[
|\Gr(1,2m+1)_{\mathbb{F}_q}|=q^{2m}\,.
\]
\end{itemize}
As was shown in Proposition \ref{polynomials}, those polynomials are related to
$p_{(1,n)}(q)$, i.e. $p_{(1,2m)}(q)=1-q^m$ and $p_{(1,2m+1)}(q)=1$, and we have the form,
\[
|\Gr(1,n)_{\mathbb{F}_q}|=q^{r}|p_{(1,n)}(q)|.
\]
\end{Example}

\begin{Example} \label{example4} We now consider $\Gr(2,4)$.   As in Example  \ref{isotropic},  
 $\Gr(2,4)$ becomes, via the Pl\"ucker embedding,  the set of one dimensional isotropic vector spaces in ${\mathbb R}^6$.  We 
then  have, as in Example  \ref{gr2},  the following description of the $\bar{\mathbb F}_q$ points of the variety we are studying
\[
\Gr(2,4)_{\bar {\mathbb F}_q} = \{ (x_1, x_2, x_3, y_1, y_2, y_3) \in \bar {\mathbb F}_q^6 : x_1^2+x_2^2+x_3^2-y_1^2-y_2^2-y_3^2=0,
x_1^2+x_2^2+x_3^2 \not=0 \}.
\]
Since $\sqrt{-1}\in {\mathbb F}_q$, we can transform the equation
$x_1^2+x_2^2+x_3^2-y_1^2-y_2^2-y_3^2=0$ giving {\it isotropy} into  $x_1^2+x_2^2+x_3^2+y_1^2+y_2^2+y_3^2=0$. The number of  solutions of $x_1^2+x_2^2+x_3^2+y_1^2+y_2^2+y_3^2=0$ is 
given by $q^5+q^3-q^2$ (see (\ref{Fpoints})). We now must subtract those points for which $x_1^2+x_2^2+x_3^2=0$
and  $y_1^2+y_2^2+y_3^2=0$.  From (\ref{Fpoints}), one can see that there are $q^2$ solutions of  $x_1^2+x_2^2+x_3^2=0$. So with an additional
$q^2$ solutions of   $y_1^2+y_2^2+y_3^2=0$, we end up with  $q^5-q^4+q^3-q^2=q^2(q-1)(q^2+1)$.
We now divide by $q-1$ to count the projectivization and obtain a total of $q^2(1+q^2)$ points.
\end{Example}

\begin{Example}
Using the formulae for $|S^{n}(\mathbb{F}_q)|$ for the case with $\sqrt{-1}\in\mathbb{F}_q$, one can write an explicit   formula giving the number of $\mathbb{F}_q$ points
of finite Chevalley group $\SO_n(\mathbb{F}_q)$:
First recall that $\SO_{n+1}(\Real)/\SO_n(\Real)\cong S^n$. Then one expects:
\[
|\SO_n(\mathbb{F}_q)|=\prod_{k=1}^n|S^{n-k}(\mathbb{F}_q)|\,,
\]
which corresponds  to the results \cite{carter:93}:
\begin{itemize}
\item[(a)] For $n=2m$,
\[
|\SO_{2m}(\mathbb{F}_q)|=2q^{m(m-1)}(q^2-1)(q^4-1)\cdots (q^{2m-2}-1)(q^m-1)\,.
\]
\item[(b)] For $n=2m+1$,
\[
|\SO_{2m+1}(\mathbb{F}_q)|=2q^{m^2}(q^2-1)(q^4-1)\cdots (q^{2m}-1)\,.
\]
\end{itemize}
In \cite{casian:06}, we show that 
those polynomials $p(q)=(q^2-1)(q^4-1)\cdots (q^{2m-2}-1)(q^m-1)$ and $p(q)=(q^2-1)(q^4-1)\cdots (q^{2m}-1)$ are related to the cohomology of the real flag variety $\SL_n(\Real)/B$ through
the singular solutions (blow-ups) of the Toda lattice.
\end{Example}

We assume that $q$ is a power of a prime $p$ , $m$ is relatively prime to $p$ and $p\not=2$.  Moreover, as before
we assume that $\sqrt{-1}\in \mathbb{F}_q$.  As  in \cite {casian:06}    one has the following.
\begin{Proposition} The cohomology   $H^*(\Gr(k,n)_{\overline{\mathbb{F}}_q}; \overline {\mathbb Q}_m)$ has Frobenius
eigenvalues of the form $q^i$.

\end{Proposition}

\begin{Proof}  This  follows from the arguments in \cite {casian:06} with almost no change.  As can be seen from  \cite {casian99} or more directly by Remark \ref {subgraph} above,
 the incidence graph of  $\Gr(k,n)$ is a subgraph of the flag manifold.  
 The Frobenius eigenvalues are still computed in terms of Hecke algebra operators as in Section 5 of   \cite {casian:06} but the  Weyl group elements in the expressions $T_{w^{-1}}^{-1}$ are restricted to a smaller subset corresponding to representatives of cells in a Grassmanian. This is also implicit  in Section 9 of  \cite {casian99} but notation is more convoluted.

One may also consider the spectral sequence associated to the fibration  
$$ (K\cap L)(\overline  {\mathbb F}_q)~ \to   ~(K)(\overline  {\mathbb F}_q) ~ \to ~(K/L\cap K)(\overline  {\mathbb F}_q)$$
This reduces the argument to the cases of $ K\cap L$ and $K$ over  $\overline  {\mathbb F}_q$.  The case of $K$ and $L\cap K$  is just the case discussed in   \cite {casian:06}  of the flag manifold.
\end{Proof}

We now recall that from the Lefschetz fixed point formula for $Fr$ the alternating sum is given by
\[
\sum_{s}(-1)^s {\rm Tr}\left((Fr)_*|_{H_c^s( \Gr(k,n)_{\overline{\mathbb{F}}_q}
;\bar {\mathbb Q}_m)}\right)=| \Gr(k,n)_{{\mathbb{F}}_q} |
\]
By  Poincar\'e duality, $H_c^{2k(n-k) -s}(  \Gr(k,n)_{\overline{\mathbb{F}}_q} ; \bar {\mathbb Q}_m)$ is the dual of
$H^{s}( \Gr(k,n)_{\overline{\mathbb{F}}_q} ;  {\mathbb V})$ where  ${\mathbb V}$ is the dual of a constant sheaf.  In the orientable cases we can replace  ${\mathbb V}$ with $\bar {\mathbb Q}_m$ i.e. constant coefficients.

\subsection {Frobenius eigenvalues calculation}

This calculation is based on section 9 of   \cite {casian99}   or  on   Section 6.1   \cite {casian:06}   which is restricted to real split cases and has simpler notation.    Frobenius eigenvalues $q^i$   increase along a graph $\G(k,n)$  or  $\mathcal{G} {(k,n)}^*$  as described in  (i),(ii), (iii) in Subsection  \ref{power}.

In terms of the action of the Hecke algebra in Lemma 3.5 of   \cite{lusztig83} , one must apply 
$T_{s+i}^{-1}=q^{-1} (T_{s_i}+(1-q))$.  Now $\Rightarrow$ corresponds to  case (e) and the Frobenius eigenvalue does not
change. The case of $\to$ corresponds to  (d2) and the Frobenius eigenvalue is multiplied by $q$.

By construction of the polynomial  $p_{(k,n)}(q)$  we have 
\[
\sum_{s}(-1)^s {\rm Tr}\left((Fr)_*|_{H^s( \Gr(k,n)_{\overline{\mathbb{F}}_q};
\bar {\mathbb Q}_m)}\right)=p_{(k,n)}(q)
\]
This leads to the following Proposition,
\begin{Proposition}   We have
\[
| \Gr(k,n)_{{\mathbb{F}}_q} |=\pm q^{k(n-k)-D}p_{(k,n)}(q),
\]
with $D=\rm{deg}(p_{(k,n)}(q))$.
\end{Proposition}

\begin{Proof}  We consider the cases with $n=$even (orientable case) and $n=$odd (non-orientable case):
\begin{enumerate}

\item Assume that $n$ is even.
We use Poincar\'e duality, $H_c^{2k(n-k)-s}(\Gr(k,n)_{\overline{\mathbb{F}}_q} ; \bar {\mathbb Q}_m)$ is the dual of $H^{s}( \Gr(k,n)_{\overline{\mathbb{F}}_q} ;  \bar {\mathbb Q}_m )$.  This corresponds  to considering the polynomial  $ q^{k(n-k)}p_{(k,n)}(q^{-1}).$
We now used the formulas  listed at the beginning  of Section \ref {pp}.  We note that  if $D=\rm{deg}(p_{(k,n)}(q))$ then, according to these formulas $q^Dp_{(k,n)}(q^{-1})=\pm p_{(k,n)}(q)$. Hence we obtain
\[
\sum_{s}(-1)^s {\rm Tr}\left((Fr)_*|_{H_c^s( \Gr(k,n)_{\overline{\mathbb{F}}_q};
;\bar {\mathbb Q}_m)}\right)=| \Gr(k,n)_{{\mathbb {F}}_q} |=\pm q^{k(n-k)-D}p_{(k,n)}(q).
\]

\item Assume that $n$ is odd.
By  Poincar\'e duality, $H_c^{2k(n-k)-s}(  \Gr(k,n)_{\overline{\mathbb{F}}_q} ; \bar {\mathbb Q}_m)$ is the dual of
$H^{s}( \Gr(k,n)_{\overline{\mathbb{F}}_q} ;  {\mathcal L} )$ where ${\mathcal L}$ is a twisted local system. 
We then have $q^{k(n-k)}p^*_{(k,n)}(q^{-1}).$  We then use Proposition \ref{star}.  We obtain 
$q^{k(n-k)}q^rp_{(k,n)}(q^{-1})$ where $r=-j$ or $r=m-j$ if $k=2j, 2j+1$ and $n=2m+1$.
\end{enumerate}
This completes the proof.
\end{Proof}

%%%%%%%%%%%%%%%%%%%%%%%%%%%%%%%%%%%%%%%%%%%%%%

\bibliographystyle{amsalpha}

\end{document}